\documentclass{amsart}
\usepackage{style} 
\usepackage{cleveref}
\setmargins{margin=2.5cm}
\usetikzlibrary{external}

\title[optimal convergence rates for parametrized variational problems]{On optimal convergence rates for parametrized variational problems with applications to nonlocal models}
\author{Zhaolong Han}
\address{Department of Mathematics, University of California, San Diego, CA 92093, United States} 
\email{zhhan@ucsd.edu}
\author{James Scott}
\address{Department of Mathematics and Statistics, Auburn University, AL 36849, United States}
\email{james.m.scott@auburn.edu}
\author{Xiaochuan Tian}
\address{Department of Mathematics, University of California, San Diego, CA 92093, United States} 
\email{xctian@ucsd.edu}
\thanks{Z.~Han and X.~Tian were supported in part by NSF grant DMS-2240180 and the Alfred P.~Sloan Fellowship. J.~Scott was supported in part by the NSF grant DMS-2606364. The early stages of this project were conducted during a SQuaRE at the American Institute for Mathematics (AIM) titled {\it Variational methods for multiscale and nonlinear nonlocal models with applications to peridynamics}. The authors thank AIM for providing a supportive and mathematically rich environment.}
\begin{document}

\begin{abstract}
  We develop an abstract variational framework for convergence rates of
parametrized problems and their Galerkin approximations, with applications
to nonlocal diffusion. The main objective is to quantify modeling and
discretization errors in terms of the nonlocal horizon, the mesh size,
the Sobolev regularity of the local solution, and the compatibility of
boundary data. We establish error estimates in energy and $L^2$ norms
and develop an adaptation of the Aubin--Nitsche argument that uses
regularity of the limiting local dual problem. This approach yields
improved $L^2$ estimates without requiring additional regularity of the
nonlocal dual solutions, including for integrable kernels. We apply the
framework to two constant-horizon nonlocal Laplacians with exterior
boundary data and to heterogeneous localization models with classical
Dirichlet, Neumann, and Robin conditions. The resulting estimates provide
quantitative guarantees of asymptotic compatibility and improve existing
error bounds under weaker regularity assumptions. Numerical experiments
illustrate the convergence rates and examine the sharpness of the estimates.
\end{abstract}



\maketitle

%\tableofcontents

\section{Introduction}\label{sec:intro}
Nonlocal models describe spatial interactions through integral operators and
provide a setting for diffusion and material behavior involving interactions
over a finite distance
\cite{silling2000reformulation,silling2010peridynamic,du2012analysis,
du2019nonlocal}.
They have seen greatly increased use in recent years in the modeling and simulation of fracture and damage via peridynamics
\cite{silling2010peridynamic,ha2010studies,hu2012peridynamic},
anomalous diffusion in complex systems \cite{metzler2000random},
and non-Fickian transport in heterogeneous porous media
\cite{benson2000application,berkowitz2006modeling}.
Numerical implementations for such models include meshfree particle methods,
finite element methods, and finite difference and spectral methods
\cite{silling2005meshfree,du2012analysis,delia2020numerical}.
When the interaction radius, or horizon, \(\del\) tends to zero, suitably
normalized nonlocal operators approach classical differential operators.
For numerical approximation, this limit introduces a second small parameter
alongside the mesh size \(h\). A fundamental question is determining the dependence on
both parameters of the error between a
discrete nonlocal solution and the corresponding local solution. On a bounded domain, the answer also depends on how the
boundary data are prescribed and on the regularity of the local solution.
This paper develops a variational approach to these error estimates and
applies it to three classes of nonlocal elliptic-type problems.

The theory of asymptotically compatible schemes provides a framework for
studying the simultaneous limits \(\del\to0\) and \(h\to0\).
Tian and Du \cite{tian2014asymptotically,tian2020asymptotically} established a framework for schemes that converge to the nonlocal solution
for a fixed horizon and recover the corresponding local solution as the
horizon and mesh size tend to zero independently; such schemes are called
asymptotically compatible (AC). For the nonlocal diffusion problems considered there,
finite element spaces containing continuous piecewise linear functions
give AC schemes, whereas convergence of piecewise constant approximations to the local
solution is established under the scale restriction \(h=o(\del)\). Du, Scott, and Tian
\cite{du2026asymptotically} extended this perspective to nonlinear
variational problems through \(\Gamma\)-convergence, including models with
heterogeneous localization and parameter-dependent domains. These results
establish convergence under broad assumptions without specifying a rate in the horizon and mesh size. %
We build on these efforts here, and establish a framework of sufficient conditions for rates of convergence in both parameters.

Quantitative estimates involving both the horizon and mesh size have
been developed for several numerical methods.
Quadrature-based finite difference and reproducing kernel collocation
schemes admit \(O(\del^2+h^2)\) estimates under \(C^4\) regularity of the
local solution on an extended domain and suitable exterior data
\cite{du2019asymptotically,leng2021asymptotically}.
Related analyses of meshfree quadrature and shifted polynomial
collocation also rely on Taylor expansions and sufficient classical regularity assumptions
\cite{chen2023analysis,trask2019asymptotically}.
For piecewise linear Galerkin approximations of volume-constrained
diffusion problems, Du, Xie, Yin, and Zhang \cite{du2024error} assumed
a \(C^4\) local solution and extension and obtained $
    O(\del^{(3+\mu)/2}+ \del^{-1} h^2)$ error bounds
     for $\mu=0,1$ 
in both the nonlocal energy and \(L^2\) norms for exact spherical
interaction neighborhoods.
These bounds assume that the prescribed exterior values approximate
a smooth extension of the local solution with maximum norm error
\(O(\del^{2+\mu})\).
Their analysis also accounts for additional errors from polygonal
approximations of the interaction neighborhoods.
For a nonlocal Neumann formulation,
Meng and Shi \cite{meng2025asymptotically} proved an \(L^2\) bound of
\(O(\del + h^k)\) for degree-\(k\) finite elements on shape-regular meshes assuming
\(H^{\max\{k+1,3\}}\) regularity of the local solution.

For a fixed horizon on a bounded domain, a common formulation prescribes
the solution on 
an exterior 
collar
\cite{du2012analysis}. The values in this collar contribute to the
nonlocal equation near the boundary, so their agreement with the local
solution influences the convergence rate. In particular, constant extension of the prescribed Dirichlet values
into the collar can limit convergence to first order in \(\del\), even
when the nonlocal operator is second-order consistent in the interior.
For smooth solutions of one-dimensional problems, incorporating the
boundary derivatives of the local solution into the exterior data
recovers an \(O(\del^2)\) convergence rate in the maximum norm
\cite{du2019uniform}.
Since these derivatives are generally unavailable from the prescribed
Dirichlet data alone, Lee and Du \cite{lee2022second} developed a boundary
treatment using nonlocal gradients to implement linear extrapolation
and proved second-order convergence under their kernel and regularity
assumptions.
An alternative approach constructs volume constraints from auxiliary
local problems and yields second-order convergence in the nonlocal
energy and \(L^2\) norms \cite{delia2020physically}.
These results show that recovering the second-order interior accuracy
on bounded domains requires a suitable treatment of the boundary
interactions.

Heterogeneous localization offers another way to treat boundaries: the
interaction radius varies with position and vanishes at the boundary,
allowing classical boundary conditions to be imposed without an exterior
collar. Trace theorems and variational formulations for such models have
been developed \cite{tao2019nonlocal,tian2017trace}. Scott and Du \cite{scott2024nonlocal,scott2025nonlocal} established
a systematic theory of the associated function spaces and variational
problems, including boundary-localized convolution operators,
Green's identities, and regularity results. For a one-dimensional problem with homogeneous local Dirichlet
conditions, Chang, Liang, and Qiao \cite{chang2026finite} analyzed
a conforming piecewise linear finite element approximation.
Their interaction radius vanishes quadratically at the boundary.
Under \(H^4\) regularity of the local solution and additional
smoothness assumptions on the localization profile, they proved
second-order operator consistency and obtained \(O(\del^2+h)\)
error bounds in both the nonlocal energy and \(L^2\) norms.

A remaining open question is ascertaining optimal error estimates in \(\del\) and \(h\)
under prescribed Sobolev regularity of the local solution. In particular, estimates based on fourth-order classical or Sobolev
regularity do not explain the convergence rates observed in practice
for solutions with lower regularity, nor do they give a general account
of how these rates depend on the compatibility of the boundary data. We address this question through a variational analysis of the nonlocal
problems and their Galerkin approximations.
For a large class of nonlocal problems, we obtain optimal estimates in \(\del\) and \(h\) with
prescribed Sobolev regularity of the local solution and compatibility of the nonlocal boundary data.

A second unresolved question concerns the improvement of the convergence rates when using the \(L^2\) norm instead of the energy norm.
For fractional-type kernels, regularity estimates for the nonlocal
dual problem can support an Aubin--Nitsche argument
\cite{burkovska2019regularity}.
For more general kernels, including integrable kernels, the nonlocal
operator need not have an elliptic smoothing property, so the Aubin--Nitsche argument cannot apply directly.
Nevertheless, computations exhibit higher \(L^2\) rates as the horizon
vanishes, and their theoretical justification has remained open
\cite{chang2026finite, du2024error}.
The second contribution of this paper is an Aubin--Nitsche-type
analysis that establishes this improvement by exploiting elliptic regularity
of the limiting local problem.
We evaluate the consistency residual on continuous and discrete dual
solutions and obtain the required estimates by comparison with the
local elliptic dual solution.
Thus improved \(L^2\) estimates follow as \(\del\to0\) even when no
corresponding higher regularity is available for the nonlocal dual
solutions.

We first apply the framework to two constant-horizon models: the usual
integral nonlocal Laplacian, called type~(I) \cite{du2012analysis}, and the composition of a
half-space nonlocal gradient and divergence, called type~(II)
\cite{han2024compactnessa,han2023nonlocal}. Both prescribe nonlocal boundary data on an exterior interaction layer.
For \(m\ge1\), let \(\Om\) be a \(C^{m,1}\) domain and assume that
the local solution and the nonlocal boundary data have \(H^{m+1}\)
regularity on their respective domains, with matching traces of the
function values and normal derivatives through order \(m-1\) on
\(\partial\Om\).
For the continuous problem, we obtain \(O(\del^{\min\{2,m-1/2\}})\) convergence in the nonlocal
energy norm and \(O(\del^{\min\{2,m\}})\) convergence in the \(L^2\) norm.
These estimates quantify how additional regularity and boundary compatibility improve
the convergence rate.
In particular, $O(\del^2)$ convergence in the \(L^2\) norm holds under \(H^3\)
regularity and matching traces of the function values and first normal
derivatives, extending the earlier \(C^4\)-based estimates to solutions
with lower regularity \cite[Theorem~2.2]{du2024error}.
For conforming Galerkin approximations with continuous piecewise
polynomials of degree \(k\ge1\), the same regularity and compatibility
assumptions give error bounds of order
\[
    O\bigl(\del^{\min\{2,m-1/2\}}+h^{\min\{m,k\}}\bigr)
    \quad\text{in the nonlocal energy norm},
\]
and
\[
    O\bigl(\del^{\min\{2,m\}}
      +(\del^{1/2}+h)h^{\min\{m,k\}}\bigr)
    \quad\text{in }L^2(\Om).
\]
For \(m=1\), \(H^2\) regularity and matching traces of the function values give
\(O(\del^{1/2}+h)\) energy and \(O(\del+h^2)\) \(L^2\) error bounds
for every degree \(k\ge1\).
For \(m=2\), \(H^3\) regularity and matching traces of the function
values and first normal derivatives give
\(O(\del^{3/2}+h)\) energy and
\(O(\del^2+\del^{1/2}h+h^2)\) \(L^2\) bounds for linear elements.
For degree \(k\ge2\), the same assumptions imply
\(O(\del^{3/2}+h^2)\) energy and \(O(\del^2+\del^{1/2}h^2+h^3)\) \(L^2\) bounds. These cases give energy and \(L^2\) error estimates in both $\del$ and $h$ that recover the results of \cite{du2024error} under a
weaker Sobolev regularity assumption. Our numerical results in \Cref{sec:Numerical} support the sharpness of the convergence rates
for the continuous problem and the error estimates in the nonlocal
energy norm for the Galerkin approximations. They also indicate that the \(L^2\) error estimates for the Galerkin
approximations can be improved when the mixed term
\(\del^{1/2}h^{\min\{m,k\}}\) cannot be absorbed into the remaining terms through Young's inequality. Obtaining sharp \(L^2\) estimates in these cases under the same
regularity and compatibility assumptions remains an open question.

For heterogeneous localization with an interaction radius comparable to
\(\del\dist(\xb,\partial\Om)\), we consider Dirichlet, Neumann, and Robin
boundary conditions on domains \(\Om\subset\mathbb R^d\).
Under the stated localization and boundary regularity assumptions, we obtain $O(\delta^{\min\{2,m\}})$ convergence in the
nonlocal energy norm for a local solution with $H^{m+1}(\Omega)$
regularity, where $m\geq1$ is an integer.
The $L^2(\Omega)$ error is $O(\delta^2)$ already under $H^2(\Omega)$ regularity.
The approximate normal derivative converges to the normal derivative
of the local solution in \(H^{-1/2}(\partial\Om)\)
at the same rate as the energy error.
For conforming Galerkin approximations with continuous piecewise
polynomials of degree \(k\ge1\), the same assumptions, together with
the required regularity of a lifting of the boundary data, give
error bounds of order
\[
    O\bigl(\del^{\min\{2,m\}}+h^{\min\{k,m\}}\bigr)
    \quad\text{in the nonlocal energy norm},
\]
and
\[
    O\bigl(\del^{2}+h^{\min\{k,m\}+1}\bigr)
    \quad\text{in }L^2(\Om).
\]
In particular, $H^2$ regularity gives $O(\delta+h)$ energy and
$O(\delta^2+h^2)$ $L^2$ error bounds for every degree $k\geq1$.
Under $H^3$ regularity, linear elements give the improved energy
bound $O(\delta^2+h)$.
For a one-dimensional homogeneous Dirichlet problem with an interaction
radius vanishing quadratically at the boundary,
\cite{chang2026finite} obtained \(O(\del^2+h)\) bounds in
both the nonlocal energy and $L^2$ norms under $H^4$ regularity.
We obtain the same energy order
under $H^3$ regularity and the improved $L^2$ bound $O(\delta^2+h^2)$
under $H^2$ regularity. Our analysis also applies in higher dimensions and covers Neumann
and Robin boundary conditions and higher polynomial degrees.
Our numerical results in \Cref{sec:Numerical} support the sharpness of these estimates
in both the nonlocal energy and \(L^2\) norms.
The analysis uses cancellation in a consistency identity
as well as the boundary-localized convolution operators of Scott and Du
\cite{scott2024nonlocal,scott2025nonlocal}.
The $L^2$ estimate also uses regularity of the local dual problem
in an Aubin--Nitsche argument.

The paper is organized as follows.
\Cref{sec:framework} develops an abstract framework for error estimates
in the energy and \(L^2\) norms.
\Cref{sec:NonlocalLaplaceI,sec:NonlocalLaplaceII} apply the framework
to the constant-horizon models of types~(I) and~(II), respectively,
and \Cref{sec:NonlocalLaplaceHeteroLoc} treats heterogeneous localization.
\Cref{sec:Numerical} presents numerical experiments illustrating
the convergence rates and examining the sharpness of the estimates.
The appendices provide the analytic properties of the half-space
nonlocal operators and the boundary-layer estimates used in the proofs.

\subsection*{AI Disclosure} The results \Cref{prop:nonlocal-boundary-layer} and \Cref{lem:hetero-H1-test-consistency} were found by modifying outputs from a series of prompts placed into ChatGPT-5.6 Sol; all other results of the manuscript were found without the use of LLMs. The entire manuscript was written by the authors, with language and grammar edited using ChatGPT-6.1 Sol. The authors take full responsibility for the correctness of all the results and proofs of this paper.

\section{Abstract framework}\label{sec:framework}  
We consider two families of Hilbert spaces \(\{(\cU_\del,\|\cdot\|_{\cU_\del}):\del\in [0, \delta_{\max}]\}\) and \(\{(\cV_\del,\|\cdot\|_{\cV_\del}):\del\in [0, \delta_{\max}]\}\)  and a family of bilinear forms \(a_\del:\cU_\del\times \cV_\del\to \R\). Let \(\cV_\del^*\) be the dual space of \(\cV_\del\). Assume that there exists a Hilbert space \(\cU_*\) such that \(\cU_\del\subset \cU_*\) for all \(\del\in [0, \delta_{\max}]\). In applications we often choose \(\cU_*\)  to be an \(L^2\) space.

Consider the following variational problems: Given \(\cF_\del\in \cV_\del^*\), find \(u_\del\in \cU_\del\) such that 
\begin{equation}\label{eq:abstract_nonlocal_eq}
        a_\del(u_\del,v)= \cF_\del(v),\quad \forall v\in \cV_\del
\end{equation}
for any \(\del\in [0, \delta_{\max}]\). Let $\cU_{\del,h}\subset \cU_\del$ and $\cV_{\del,h}\subset \cV_\del$ be the finite dimensional subspaces and $u_{\del,h}\in \cU_{\del,h}$ be the Galerkin solution given by \eqref{eq:abstract_discrete_eq} below.
Our main results are the following:
\begin{enumerate}
    \item Convergence rates in $\del$ of $u_\del$ to $u_0$ in both the $\cU_\del$ and $L^2$ norms as $\del\to 0$.
    \item Convergence rates in $\del$ and $h$ of $u_{\del,h}$ to $u_0$ in both the $\cU_\del$ and $L^2$ norms as $\del\to 0$ and $h\to 0$.
\end{enumerate}

The next three subsections provide the abstract framework for such results, which will be applied to various types of nonlocal diffusion problems with their discretizations.

\subsection{Abstract energy error estimates on the continuous level}
We start with the convergence rate of \(u_\del\) to \(u_0\) in the energy norm. 

\begin{assumption}\label{assu:ellip} The spaces \(\{\cU_\del\}_\del\), \(\{\cV_\del\}_\del\) and the bilinear forms \(\{a_\del(\cdot, \cdot)\}_\del\) satisfy the following properties:
    \begin{enumerate}[label=(\roman*)]
        \item \label{item:uniform_embedding} Uniform embedding property: there exist constants \(M_*>0\) and \(M_0>0\) independent of \(\del\in [0, \delta_{\max}]\) such that  
        \begin{equation*}
            M_*\|u\|_{\cU_*}\le \|u\|_{\cU_\del},\quad \forall u\in \cU_\del,
        \end{equation*}
        and
        \begin{equation*}
            \|u\|_{\cU_\del}\le M_0\|u\|_{\cU_0},\quad \forall u\in \cU_0.
        \end{equation*} 
        \item \(\{a_\del(\cdot, \cdot)\}_\del\) is uniformly bounded: there exists a constant \(M_a>0\) independent of \(\del\) such that 
        \begin{equation*} 
                \abs{a_\del(u,v)}\le M_a\|u\|_{\cU_\del}\|v\|_{\cV_\del},\quad\forall u\in \cU_\del,\ v\in \cV_\del. 
        \end{equation*} 
        \item \label{item:inf_sup} \(\{a_\del(\cdot, \cdot)\}_\del\) satisfies the inf-sup condition: 
        \begin{enumerate}
            \item there exists  a constant \(\al>0\) independent of \(\del\) such that 
            \begin{equation*} 
                    \inf_{u\in \cU_\del}\sup_{v\in \cV_\del}
                    \frac{\abs{a_\del(u, v)}}{\|u\|_{\cU_\del} \|v\|_{\cV_\del}} \ge \al. 
            \end{equation*}
            \item If \(v\in \cV_\del\) such that \(a_\del(u, v) = 0\) for all \(u\in \cU_\del\), then \(v = 0\).
        \end{enumerate}   
    \end{enumerate}
\end{assumption}
\begin{remark} 
    Notice that we do not assume asymptotic compact embedding of \(\{\cU_\del\}_{\del}\) into \(\cU_*\) as in \cite{tian2014asymptotically}. However, the fact that \(M_*\) is independent of \(\del\) usually requires such compactness. Also the first embedding in \Cref{assu:ellip}\ref{item:uniform_embedding} implies that once the error is bounded in \(\cU_\del\)-norm, then it is automatically bounded in \(\cU_*\)-norm, so convergence in \(\cU_\del\)-norm implies convergence in \(\cU_*\)-norm, which is one of the main results shown in \cite{tian2014asymptotically} via a compactness argument.
\end{remark}
By \Cref{assu:ellip}, given \(\del\in [0, \delta_{\max}]\), the Banach-Ne\v cas-Babu\v ska theorem \cite[Theorem~25.9]{ern2021finite} implies that  there exists a unique solution \(u_\del\) for \eqref{eq:abstract_nonlocal_eq} and
\[\norm{u_\del}_{\cU_\del}\le \frac{1}{\al}\norm{\cF_\del}_{\cV_\del^*}.\]

For the continuous problems, our goal is to analyze the error
\begin{equation*}
    e_\del:=u_\del-\uext.
\end{equation*}
By \Cref{assu:ellip}(i), \(e_\del\in \cU_\del\). It also satisfies the equation 
\begin{equation}\label{eq:error-equation}
    a_\del(e_\del,v)=\cR_\del(v),
    \qquad \forall v\in \cV_\del,
\end{equation}
where $\cR_\del$ is the residual functional given by
\begin{equation*}
    \cR_\del(v):=\cF_\del(v)-a_\del(\uext,v),
    \qquad \forall v\in \cV_\del.
\end{equation*}

\begin{lemma}\label{lem:abstract-energy-first}
Under \Cref{assu:ellip}, 
\begin{equation*}
    \norm{u_\del-\uext}_{\cU_\del}
    \le \frac{1}{\al}\norm{\cR_\del}_{\cV_\del^*}.
\end{equation*}
\end{lemma}
\begin{proof}
    Since \(e_\del\in \cU_\del\), by \Cref{assu:ellip}\ref{item:inf_sup} and \eqref{eq:error-equation}, 
    \[
    \al\norm{e_\del}_{\cU_\del}\le \sup_{v\in\cV_\del}\frac{a_\del(e_\del, v)}{\norm{v}_{\cV_\del}}=\sup_{v\in\cV_\del}\frac{\cR_\del(v)}{\norm{v}_{\cV_\del}}=\norm{\cR_\del}_{\cV_\del^*}.
    \]
The result follows after division by \(\al\).
\end{proof}
 
\begin{assumption}\label{assu:abstract-residual-estimate}
    There exist constants $C_{E}>0$ and $\theta>0$, independent of \(\del\), such that $\norm{\cR_\del}_{\cV_\del^*}\leq C_{E}\del^\theta$, namely,
\begin{equation*}
    \abs{\cR_\del(v)}
    \le C_{E}\del^\theta \norm{v}_{\cV_\del},
    \qquad \forall v\in \cV_\del.
\end{equation*}
\end{assumption}

\begin{theorem}\label{thm:abstract-energy-rate}
Suppose that  \Cref{assu:ellip} and \Cref{assu:abstract-residual-estimate} hold.  
Then
\begin{equation*}
    \norm{u_\del-\uext}_{\cU_\del}
    \le
    \frac{1}{\al} C_{E}\del^\theta.
\end{equation*}
\end{theorem}
\begin{proof} 
The result follows from \Cref{assu:abstract-residual-estimate} and \Cref{lem:abstract-energy-first}.
\end{proof}

\subsection{Abstract energy error estimate with discretization}

Now we consider discrete problems. Let \(\{\cU_{\del,h}:\del\in [0, \delta_{\max}], h\in (0, h_0]\}\) be a family of finite-dimensional subspaces of \(\cU_\del\) such that \(\cU_{\del,h}\subset \cU_{\del}\) for all \(\del\in [0, \delta_{\max}]\). Let \(\{\cV_{\del,h}:\del\in [0, \delta_{\max}], h\in (0, h_0]\}\) be a family of finite-dimensional subspaces of \(\cV_\del\) such that \(\cV_{\del,h}\subset \cV_{\del}\) for all \(\del\in [0, \delta_{\max}]\). We assume throughout that $\dim \cU_{\del,h}=\dim \cV_{\del,h}$. 
The Galerkin method seeks \(u_{\del,h}\in \cU_{\del,h}\) such that
\begin{equation}\label{eq:abstract_discrete_eq}
    a_\del(u_{\del,h},v_h ) =  \cF_\del (v_h),\quad \forall v_h \in \cV_{\del,h}.
\end{equation}
\begin{assumption}\label{assu:ellip_discrete}
    The discrete spaces \(\{\cU_{\del,h}:\del\in [0, \delta_{\max}], h\in (0, h_0]\}\) and \(\{\cV_{\del,h}:\del\in [0, \delta_{\max}], h\in (0, h_0]\}\) satisfy the discrete inf-sup condition: there exists a constant \(\tilde{\al}>0\) independent of \(\del\) and \(h\) such that
        \begin{equation*}
            \inf_{u_h\in \cU_{\del,h}}
            \sup_{v_h\in \cV_{\del,h}}
            \frac{\abs{a_\del(u_h,v_h)}}
            {\|u_h\|_{\cU_\del}\|v_h\|_{\cV_\del}}
            \ge \tilde{\al}.
        \end{equation*}
\end{assumption}
Again \cite[Theorem~25.9]{ern2021finite} shows that \(u_{\del,h}\) is the unique solution of
\eqref{eq:abstract_discrete_eq}.
\begin{remark}\label{rem:coercive}
    When \(\cU_\del=\cV_\del\) for all \(\del\in [0, \delta_{\max}]\) and \(\{a_\del(\cdot, \cdot)\}_\del\) is uniformly coercive, i.e., there exists a constant \(\al_0>0\) independent of \(\del\) such that
        \begin{equation*}
                a_\del(u,u)\ge \al_0\|u\|_{\cV_\del}^2,\quad\forall u\in \cV_\del, 
        \end{equation*}
    it is easy to see that \Cref{assu:ellip}\ref{item:inf_sup} and \Cref{assu:ellip_discrete} hold with \(\al=\tilde{\al}=\al_0\).
\end{remark}

We now present the discrete error in the energy norm.
\begin{theorem}\label{thm:abstract-energy-discrete}
    Suppose that \Cref{assu:ellip} and \Cref{assu:ellip_discrete} hold.  Then 
    \begin{equation*}
            \|u_{\del,h} - \uext\|_{\cU_\del}  \le \frac{1}{\tilde{\al}} \|\Rdel\|_{\cV_\del^*} + \left(1+\frac{M_a}{\tilde{\al}}\right)\inf_{u_h\in \cU_{\del,h}}\|\uext - u_h\|_{\cU_\del}.
    \end{equation*}  
\end{theorem}
\begin{proof}
    Let \(u_h\in \cU_{\del,h}\). By \Cref{assu:ellip_discrete}, we have
    \begin{equation*} 
            \tilde{\al}\|u_{\del,h} - u_h\|_{\cU_\del} 
            \le   \sup_{v_h\in \cV_{\del,h}} \frac{a_\del(u_{\del,h} - u_h, v_h)}{\|v_h\|_{\cV_\del}}  
            =  \sup_{v_h\in \cV_{\del,h}} \frac{\Rdel(v_h) + a_\del(\uext - u_h,v_h)}{\|v_h\|_{\cV_\del}}  
            \le   \|\Rdel\|_{\cV_\del^*} + M_a \|\uext  - u_h\|_{\cU_\del}. 
    \end{equation*}  
    Using the triangle inequality \(\|\uext -u_{\del,h}\|_{\cU_\del}\le \|\uext-u_h\|_{\cU_\del}+\|u_h-u_{\del,h}\|_{\cU_\del}\) and taking the infimum over \(u_h\in \cU_{\del,h}\) yields the desired estimate. 
\end{proof}

\subsection{From energy rates to \texorpdfstring{\(L^2\)}{L2} rates.}
From the first embedding in \Cref{assu:ellip}\ref{item:uniform_embedding}, it is clear that the convergence rate in $L^2$ norm is at least of the same order as the energy norm. Similar to the classical finite element theory, we apply a duality trick analogous to the Aubin-Nitsche lemma to improve the $L^2$ rate.

We assume \(\cU_*\) is an \(L^2\) space. For \(\phi\in L^2\), define the dual solution \(z_\del[\phi]\in \cV_\del\) by
\begin{equation}\label{eq:dual-problem-abstract}
    a_\del(u,z_\del[\phi])= (u, \phi)_{L^2},
    \qquad \forall u\in \cU_\del.
\end{equation} 
By \cite[Corollary~C.52, Lemma C.53]{ern2021finite}, for every \(\phi\in L^2\) there exists a unique solution \(z_\del[\phi]\) to \eqref{eq:dual-problem-abstract} and, thanks to the first embedding in \Cref{assu:ellip}\ref{item:uniform_embedding},
\[
\norm{z_\del[\phi]}_{\cV_\del}\le \frac{1}{\al}\sup_{u\in \cU_\del} \frac{(u, \phi)_{L^2}}{\norm{u}_{\cU_\del}}
\le \frac{1}{\al M_*}\norm{\phi}_{L^2}.
\]
Combining this bound with \Cref{assu:abstract-residual-estimate} we obtain that \(\abs{\cR_\del(z_\del[\phi])}\le \tfrac{C_E}{\al M_*} \del^\theta \norm{\phi}_{L^2}\).  This estimate motivates the following assumption which is essential for an improved $L^2$ convergence rate.
\begin{assumption}\label{assu:abstract-dualtest-residual}
There exist constants $C_{L}>0$ and $\mu\ge\theta$  such that
\begin{equation*}
    \abs{\cR_\del(z_\del[\phi])}
    \le C_{L}\del^\mu \norm{\phi}_{L^2}, \quad \forall \phi\in L^2,
\end{equation*}
where \(\theta\) is given in \Cref{assu:abstract-residual-estimate}.
\end{assumption}

\begin{theorem}\label{thm:abstract-L2-rate}
Suppose that \Cref{assu:ellip,assu:abstract-dualtest-residual} hold. Then
\begin{equation*}
    \norm{u_\del-\uext}_{L^2}
    \le C_L \del^\mu.
\end{equation*}
\end{theorem}
\begin{proof}
Let \(e_\del=u_\del-\uext\in \cU_\del\).  Choose \(\phi=e_\del\in L^2 \) and $u=e_\del$ in the dual problem \eqref{eq:dual-problem-abstract}.  Then by \eqref{eq:error-equation}
\[
    \norm{e_\del}_{L^2}^2
    =(e_\del,\phi)_{L^2}
    =a_\del(e_\del,z_\del[e_\del])
    =\cR_\del(z_\del[e_\del]).
\]
Using \Cref{assu:abstract-dualtest-residual} with \(\phi=e_\del\),
\[
    \norm{e_\del}_{L^2}^2
    \le C_L \del^\mu \norm{e_\del}_{L^2}.
\]
Dividing by \(\norm{e_\del}_{L^2}\) proves the theorem.
\end{proof}

For the discrete problem \eqref{eq:abstract_discrete_eq}, it is useful to introduce the discrete dual solution \(z_{\del,h}[\phi]\in \cV_{\del,h}\), defined for \(\phi\in L^2\) by
\begin{equation}\label{eq:dual-problem-abstract-discrete}
    a_\del(u_h,z_{\del,h}[\phi])= (u_h, \phi)_{L^2},
    \qquad \forall u_h\in \cU_{\del,h}.
\end{equation}
For the same reason, \eqref{eq:dual-problem-abstract-discrete} has a unique
solution.
However, as
\(u_{\del,h}-\uext\notin \cU_{\del,h}\), the proof of
\Cref{thm:abstract-L2-rate} cannot be applied directly to derive a similar estimate. We propose the following
assumption and prove the discrete \(L^2\) convergence rates.

\begin{assumption} 
\label{assu:abstract-dual-discrete}
Suppose that there exist nonnegative functions \(\beta(\del,h)\) and
\(\rho(\del,h)\) such that, for every \(\phi\in L^2\),
\begin{equation*}
    \abs{\Rdel(z_{\del,h}[\phi])}
    \le
    \beta(\del,h)\norm{\phi}_{L^2},
\end{equation*}
and
\begin{equation*}
    \inf_{v_h\in \cV_{\del,h}}
    \norm{z_\del[\phi]-v_h}_{\cV_\del}
    \le
    \rho(\del,h)\norm{\phi}_{L^2}.
\end{equation*}
\end{assumption}

The functions \(\beta(\del,h)\) and \(\rho(\del,h)\) will characterize the rates in terms of \(\del\) and \(h\): \(\beta(\del,h)\) measures consistency against the discrete dual solution,
and \(\rho(\del,h)\) measures approximation of the continuous dual solution.

\begin{theorem}
\label{thm:abstract-L2-rate-discrete}
Suppose that \Cref{assu:ellip,assu:abstract-dual-discrete,assu:ellip_discrete} hold.
Then
\begin{equation*}
    \norm{u_{\del,h}-\uext}_{L^2}
    \le 
    \beta(\del,h) +
    M_a\left(1+\frac{M_a}{\tilde{\al}}\right)\rho(\del,h)
    \inf_{u_h\in \cU_{\del,h}}
    \norm{\uext-u_h}_{\cU_\del}. 
\end{equation*}
\end{theorem}

\begin{proof}
We first introduce the Galerkin projection operator
\(\Pi_{\del,h}:\cU_\del\to \cU_{\del,h}\), defined by
\begin{equation*}
    a_\del(\Pi_{\del,h}u,v_h)
    =
    a_\del(u,v_h),
    \qquad \forall v_h\in \cV_{\del,h},
\end{equation*}
for \(u\in\cU_\del\). Indeed, the existence and boundedness of the operator \(\Pi_{\del,h}\) are guaranteed by the discrete inf-sup condition and uniform boundedness of the bilinear form; see \cite[Theorem~2.22]{ern2004theory}.

We next record the approximation property of this projection. For every \(u\in \cU_\del\) and
\(u_h\in \cU_{\del,h}\), the discrete inf-sup condition, the projection
identity, and continuity give 
\[\tilde{\al}\norm{\Pi_{\del,h}u-u_h}_{\cU_\del}
    \le
    \sup_{v_h\in \cV_{\del,h}}
    \frac{\abs{a_\del(\Pi_{\del,h}u-u_h,v_h)}}
    {\norm{v_h}_{\cV_\del}}
    =
    \sup_{v_h\in \cV_{\del,h}}
    \frac{\abs{a_\del(u-u_h,v_h)}}
    {\norm{v_h}_{\cV_\del}}\le
    M_a\norm{u-u_h}_{\cU_\del}.\]
The triangle inequality therefore gives
\begin{equation*}
    \norm{u-\Pi_{\del,h}u}_{\cU_\del}
    \le
    \left(1+\frac{M_a}{\tilde{\al}}\right)
    \inf_{u_h\in \cU_{\del,h}}\norm{u-u_h}_{\cU_\del},\qquad \forall u\in \cU_\del.
\end{equation*}

Notice that by \eqref{eq:dual-problem-abstract-discrete},
\[
\| u_{\del, h} - u_0 \|_{L^2}^2 = (u_{\del, h} - u_0, u_{\del, h} - u_0)_{L^2} = a_\del ( u_{\del, h}, z_{\del, h}[u_{\del, h}-u_0]) - (u_0, u_{\del, h}-u_0)_{L^2}. 
\]
For the first term,   using \eqref{eq:abstract_nonlocal_eq}, \eqref{eq:abstract_discrete_eq}, \eqref{eq:error-equation} and  the definition of \(\Pi_{\del,h}\), for any $v_h \in \cV_{\del, h}$,
\[
a_\del(u_{\del, h}, v_h) = a_\del (u_\del - u_0, v_h) + a_\del (u_0 , v_h)= \cR_\del(v_h) + a_\del (\Pi_{\del,h}u_0 , v_h).
\]
Therefore, by \eqref{eq:dual-problem-abstract-discrete},
\[
\begin{split}
\| u_{\del, h} - u_0 \|_{L^2}^2 & = \cR_\del(z_{\del,h} [u_{\del,h} - u_0]) + a_\del(\Pi_{\del, h} u_0,z_{\del,h} [u_{\del,h} - u_0] ) - (u_0, u_{\del, h}-u_0)_{L^2}  \\
&= \cR_\del(z_{\del,h} [u_{\del,h} - u_0]) + (\Pi_{\del, h} u_0 - u_0 , u_{\del, h}-u_0)_{L^2}. 
\end{split}
\]
Using the first bound in \Cref{assu:abstract-dual-discrete} with \(\phi=u_{\del,h} - u_0\) and Cauchy-Schwarz for $(\Pi_{\del, h} u_0 - u_0 , u_{\del, h}-u_0)_{L^2}$, and then canceling \(\| u_{\del, h} - u_0 \|_{L^2}\) from both sides, we obtain that
\begin{equation}\label{eq:tmp0}
\norm{u_{\del,h}-\uext}_{L^2}
    \le 
    \beta(\del,h) + \norm{u_0 - \Pi_{\del, h} u_0}_{L^2}.
\end{equation}

It remains to estimate the last term in the previous inequality. Let
\(z_\del[\uext-\Pi_{\del,h}\uext]\) solve \eqref{eq:dual-problem-abstract}. For every \(v_h\in \cV_{\del,h}\), the definition of \(\Pi_{\del,h}\) yields
\[
\begin{aligned}
    \norm{\uext-\Pi_{\del,h}\uext}_{L^2}^2
    &=
    a_\del(\uext-\Pi_{\del,h}\uext,z_\del[\uext-\Pi_{\del,h}\uext])\\
    &=
    a_\del(\uext-\Pi_{\del,h}\uext,z_\del[\uext-\Pi_{\del,h}\uext]-v_h)\\
    &\le
    M_a\norm{\uext-\Pi_{\del,h}\uext}_{\cU_\del}
    \norm{z_\del[\uext-\Pi_{\del,h}\uext]-v_h}_{\cV_\del}.
\end{aligned}
\]
Taking the infimum over \(v_h\), using
the second bound in \Cref{assu:abstract-dual-discrete}, and then applying
the approximation bound for \(\Pi_{\del,h}\) to $u_0$, we obtain
\begin{equation*}
    \norm{\uext-\Pi_{\del,h}\uext}_{L^2}
    \le
    M_a\left(1+\frac{M_a}{\tilde{\al}}\right)\rho(\del,h)
    \inf_{u_h\in \cU_{\del,h}}
    \norm{\uext-u_h}_{\cU_\del}.
\end{equation*}
Combining this estimate with \eqref{eq:tmp0} completes the proof.
\end{proof}

\section{Applications to nonlocal Laplacians with constant horizon of type (I)}\label{sec:NonlocalLaplaceI}

Throughout this paper, we consider the Poisson equation as the model problem. Let \(\Om\subset \R^d\) be a bounded Lipschitz domain. Given \(f\in L^2(\Om)\) and \(g\in H^{1/2}(\p\Om)\),
\begin{equation}\label{eq:local_Poisson_strong}
    \left\{\begin{aligned}
        -\Delta U_0 &= f \ \text{in}\ \Om,\\
        U_0 &= g \ \text{on} \ \partial\Om.
    \end{aligned}\right.
\end{equation} 

Our goal is to construct corresponding nonlocal problems with appropriate nonlocal boundary conditions and study the convergence rate of nonlocal solutions to the local one. In this section, we study the constant horizon Laplacian of type (I) which is given by a single integral involving a kernel function with finite second moment. This type of nonlocal Laplacian has been studied extensively. In particular, \cite{tian2014asymptotically} shows that the nonlocal solution \(U_\del\) converges to the local solution \(U_0\) in \(L^2\) norm as \(\del\to 0\), and the discretized nonlocal solution \(U_{\del,h}\) converges to the local solution \(U_0\) in \(L^2\) norm as \(\del\to 0\) and \(h\to 0\) provided that  the finite element spaces  contain continuous piecewise linear functions. Their proof is based on a compactness embedding of energy space into \(L^2\) space as nonlocality vanishes, hence there is no explicit convergence rate. In this section, we establish convergence rates in both nonlocal and \(L^2\) norms with constants depending on the regularity of the local solution.

For the convergence-rate results, we assume in addition that \(\Om\) is
\(C^{1,1}\). We fix \(\del_0>0\) as specified in
\eqref{eq:def-delta-0} below.
In the notation of \Cref{sec:framework}, we take
\(\del_{\max}:=\del_0\). Throughout this section,
\(0<\del\le\del_0\).

For every \(\varepsilon>0\), define
\begin{equation}\label{eq:notation-collars}
\begin{aligned}
    \Omega^+_\varepsilon
    &:=\{\xb\in\mathbb R^d\setminus\Omega:\dist(\xb,\Omega)<\varepsilon\},
    &\widehat\Omega_\varepsilon&:=\Omega\cup\Omega^+_\varepsilon,\\
    \Omega^\circ_\varepsilon
    &:=\{\xb\in\Omega:\dist(\xb,\partial\Omega)>\varepsilon\},
    &\Omega^-_\varepsilon&:=\{\xb\in\Omega:\dist(\xb,\partial\Omega)<\varepsilon\}.
\end{aligned}
\end{equation}

\begin{figure}[htbp]
\centering
\begin{tikzpicture}[>=Stealth,line join=round]
  \def\bl#1#2{2.0+0.21*cos(2*(#1))+0.12*sin(3*(#1))+(#2)}
  \def\ep{0.36}

  %------------------------------- panel (a): the four sets
  \fill[black!12]  plot[variable=\t,domain=0:360,samples=200,smooth cycle] (\t:{\bl{\t}{\ep}});
  \fill[orange!35] plot[variable=\t,domain=0:360,samples=200,smooth cycle] (\t:{\bl{\t}{0}});
  \fill[blue!14]   plot[variable=\t,domain=0:360,samples=200,smooth cycle] (\t:{\bl{\t}{-\ep}});

  \draw[dashed,black!55] plot[variable=\t,domain=0:360,samples=200,smooth cycle] (\t:{\bl{\t}{\ep}});
  \draw[very thick]      plot[variable=\t,domain=0:360,samples=200,smooth cycle] (\t:{\bl{\t}{0}});
  \draw[dashed,blue!65]  plot[variable=\t,domain=0:360,samples=200,smooth cycle] (\t:{\bl{\t}{-\ep}});

  \node at (0,0) {$\Omega^{\circ}_{\varepsilon}$};

  \coordinate (Pin)  at (132:{\bl{132}{-0.18}});
  \coordinate (Pout) at (62:{\bl{62}{0.18}});
  \coordinate (Pbd)  at (8:{\bl{8}{0}});
  \coordinate (Phat) at (-118:{\bl{-118}{0.18}});

  \node (Lin)  at (-3.6,1.9) {$\Omega^{-}_{\varepsilon}$};
  \node (Lout) at (1.3,3.1)  {$\Omega^{+}_{\varepsilon}$};
  \node (Lbd)  at (3.5,0.6)  {$\partial\Omega$};
  \node (Lhat) at (-1.6,-3.2){$\widehat{\Omega}_{\varepsilon}=\Omega\cup\Omega^{+}_{\varepsilon}$};

  \draw[->,thin] (Lin)  -- (Pin);
  \draw[->,thin] (Lout) -- (Pout);
  \draw[->,thin] (Lbd)  -- (Pbd);
  \draw[->,thin] (Lhat) -- (Phat);

  \draw[<->,thin,black!70] (-38:{\bl{-38}{0}}) --
        node[fill=white,inner sep=1pt,font=\footnotesize] {$\varepsilon$} (-38:{\bl{-38}{\ep}});

  \node at (0,-3.9) {(a)};

  %------------------------------- panel (b): the two-sided collar
  \begin{scope}[xshift=7.4cm]
    \def\ee{0.19}
    \fill[teal!30] plot[variable=\t,domain=0:360,samples=200,smooth cycle] (\t:{\bl{\t}{\ee}});
    \fill[white]   plot[variable=\t,domain=0:360,samples=200,smooth cycle] (\t:{\bl{\t}{-3*\ee}});

    \draw[dashed,black!55] plot[variable=\t,domain=0:360,samples=200,smooth cycle] (\t:{\bl{\t}{\ee}});
    \draw[very thick]      plot[variable=\t,domain=0:360,samples=200,smooth cycle] (\t:{\bl{\t}{0}});
    \draw[dashed,teal!70]  plot[variable=\t,domain=0:360,samples=200,smooth cycle] (\t:{\bl{\t}{-3*\ee}});

    \node at (0,0) {$\Omega^{\circ}_{a\varepsilon}$};

    \coordinate (Q1) at (132:{\bl{132}{-1.5*\ee}});
    \coordinate (Q2) at (62:{\bl{62}{0.5*\ee}});
    \node (LQ1) at (-3.6,1.9) {$\Omega^{-}_{a\varepsilon}$};
    \node (LQ2) at (1.3,3.1)  {$\Omega^{+}_{b\varepsilon}$};
    \draw[->,thin] (LQ1) -- (Q1);
    \draw[->,thin] (LQ2) -- (Q2);

    \node[align=center] at (0,-3.2)
      {$(\partial\Om)^{a,b}_{\varepsilon}=\Omega^{-}_{a\varepsilon}\cup\Omega^{+}_{b\varepsilon}$ (shaded),\\
       drawn for $a=3$, $b=1$};
    \node at (0,-3.9) {(b)};
  \end{scope}
\end{tikzpicture}
\caption{The sets introduced in \eqref{eq:notation-collars} and \eqref{eq:near-boundary-layer}.}
\label{fig:domains}
\end{figure}
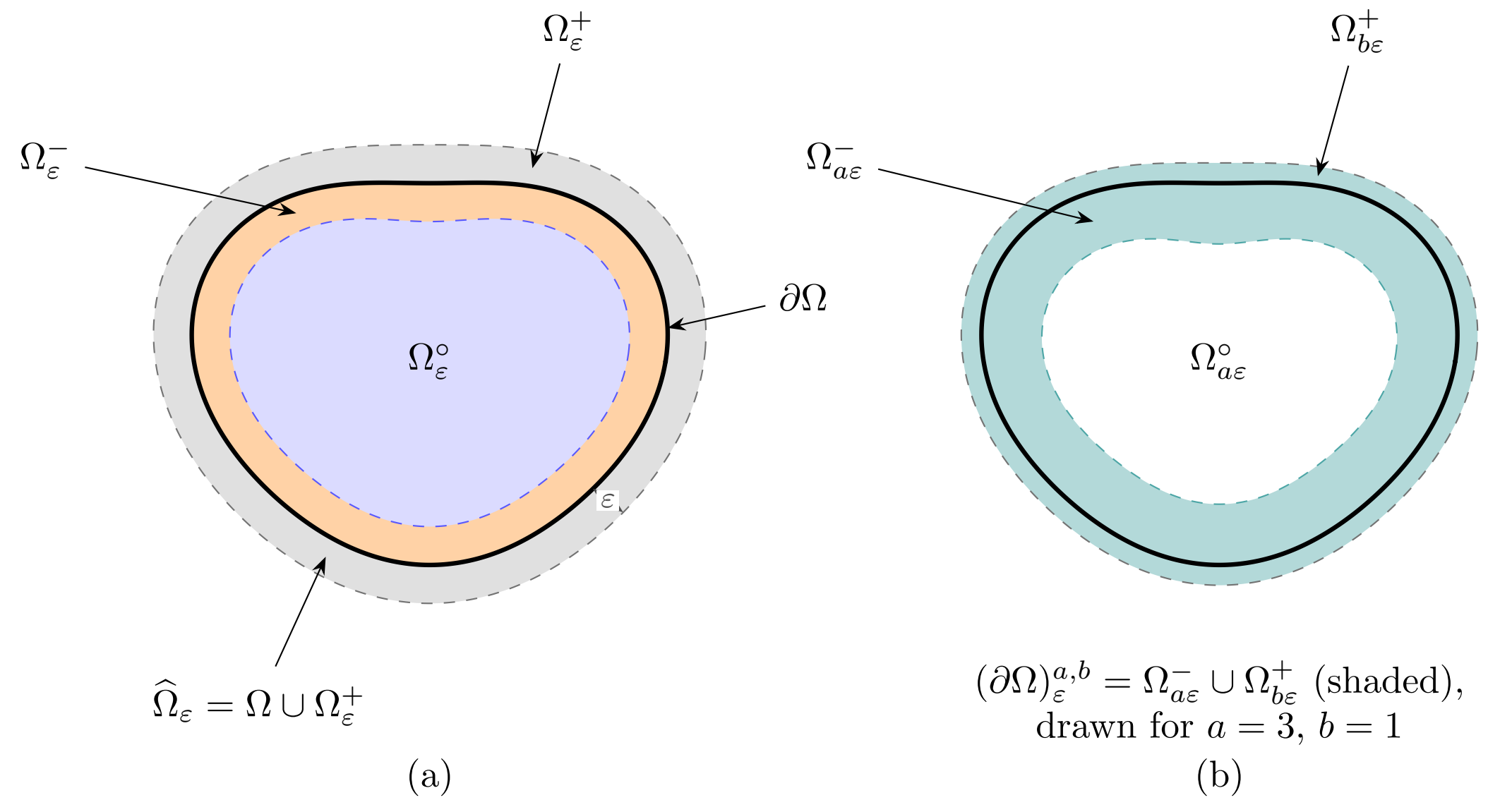

For \(a,b>0\), we also set 
\begin{equation}\label{eq:near-boundary-layer}
    (\partial\Om)^{a,b}_\varepsilon
    :=\Omega^-_{a\varepsilon}\cup\Omega^+_{b\varepsilon}.
\end{equation}
These sets are shown in \Cref{fig:domains}. In the following these sets will be used for \(\vep=\del\in (0,\del_0]\).

% =====================================================================
%  FIGURE 1 --- domains and collars.   
% =====================================================================

Throughout the manuscript, \(\Tr\) denotes the boundary trace from the
domain on which its argument is defined. For a function defined on a
neighborhood of \(\partial\Om\), it denotes the trace of the restriction
to \(\partial\Om\). 

For the local Dirichlet boundary datum \(g\in H^{1/2}(\p\Om)\), we fix an extension \(G\in H^1(\hat{\Om}_{\del_0})\) such that
\begin{equation*}
    \norm{G}_{H^1(\hat{\Om}_{\del_0})}\le C \norm{g}_{H^{1/2}(\p\Om)},\quad \Tr G = g\ \text{ on }\ \p\Om.
\end{equation*}
Higher regularity of \(G\) will be assumed later for the convergence rate results.

The nonlocal Poisson problem in strong form reads as follows
\cite{tian2020asymptotically}:
\begin{equation}\label{eq:fullSpace:BVP:strong}
    \left\{\begin{aligned} 
        -\cL_{\del} U_\del &= f\ \text{in }\Om,\\
        U_\del &= G\ \text{in }\Omp,
    \end{aligned}\right.
\end{equation}
where the nonlocal Laplacian \(\cL_\del\) is defined as 
\begin{equation*}
    \cL_{\del} u (\xb) = \operatorname{P.V.}\;2\int_{\R^d} \gamma_\del(\yb - \xb) (u(\yb) - u(\xb))\,d\yb,\quad \xb\in\Om.
\end{equation*}

The nonlocal kernel \(\gamma_\del\) is defined by  
\begin{equation}\label{eq:gamma-scaling}
    \gamma_\del(\zb) := \frac{1}{\del^{d+2}}\gamma\left(\frac{\zb}{\del}\right),
\end{equation}
where \(\ga\) satisfies the following assumption.

\begin{assumption}
\label{assu:kernel}
The kernel \(\gamma\) is radial, nonnegative, compactly supported in \(B_1(\bm 0)\), and its second moment equals \(d\), i.e.,
\begin{equation}\label{eq:gamma:2ndMoment}
    \int_{B_1(\bm 0)} |\zb|^2 \gamma(\zb)\,d\zb = d.
\end{equation}
In addition, we assume the following nondegeneracy condition: there exist constants \(c_0>0\) and \(c_\ga\in(0,1]\) such that
\begin{equation}\label{eq:gamma-lower-bound}
    \gamma(\zb)\ge c_0
    \quad\text{for a.e. }\zb\in B_{c_\ga}(\bm 0).
\end{equation}
\end{assumption}

We first define 
\begin{equation*} 
    \abs{u}_{\cS_\del(\mathbb R^d)}^2
    :=
    \iint_{\mathbb R^d\times\mathbb R^d}
    \ga_\del(\xb-\yb)
    \bigl(u(\xb)-u(\yb)\bigr)^2\,d\yb\,d\xb,
\end{equation*}
and set
\begin{equation*} 
    \cS_\del(\mathbb R^d)
    :=
    \left\{u\in L^2(\mathbb R^d):
    \abs{u}_{\cS_\del(\mathbb R^d)}<\infty\right\}.
\end{equation*}
This space is equipped with the full norm
\[
    \norm{u}_{\cS_\del(\mathbb R^d)}
    :=
    \left(
    \norm{u}_{L^2(\mathbb R^d)}^2+
    \abs{u}_{\cS_\del(\mathbb R^d)}^2
    \right)^{1/2}.
\]

On the bounded domain \(\Om\), define 
\[
\cS_{\del}(\Om):=\left\{u\in L^2(\Omh):
    \abs{u}_{\cS_\del(\Om)}<\infty \right\},
\]
where 
\begin{equation*}
    \abs{u}_{\cS_\del(\Om)}^2
    :=
    \iint_{\Omh\times\Omh}
    \ga_\del(\xb-\yb)
    \bigl(u(\xb)-u(\yb)\bigr)^2\,d\yb\,d\xb.
\end{equation*}
Equip this space with the norm
\[
    \norm{u}_{\cS_\del(\Om)}
    :=\left(
        \norm{u}_{L^2(\widehat\Om_\del)}^2
        +\abs{u}_{\cS_\del(\Om)}^2
    \right)^{1/2}.
\]
The homogeneous volume constraint is imposed by the closed subspace
\cite{mengesha2013analysis,tian2014asymptotically}
\begin{equation*}
    \cS_{\del,0}(\Om)
    :=
    \left\{u\in \cS_\del(\Om):
        u=0\ \text{a.e. in }\Omp\right\}.
\end{equation*} 
Both spaces are Hilbert spaces  
\cite{mengesha2013analysis}.
We use the following uniform Poincar\'e inequality
\cite[Lemma~3.2]{tian2014asymptotically}: there exist
constants \(\del_{P}=\del_{P}(\Om,d,\ga)>0\) and \(C_P=C_P(\del_{P})>0\)  such that
for \(\del\in (0, \del_P)\),
\[
    \norm{u}_{L^2(\Om)}
    \le C_P\abs{u}_{\cS_\del(\Om)},
    \qquad \forall u\in\cS_{\del,0}(\Om).
\]
Consequently, the seminorm is equivalent to the full norm on \(\cS_{\del,0}(\Om)\):
\[\norm{u}_{\cS_\del(\Om)} \le (1+C_P^2)^{\frac{1}{2}}\abs{u}_{\cS_\del(\Om)},\qquad \forall u\in \cS_{\del,0}(\Om).\]  
Also, by \cite[Lemma~3.1]{tian2014asymptotically} there exists a constant \(C_0=C(\Om, d) > 0\) such that 
\begin{equation}\label{eq:H1embedS_fullSpace}
    \abs{u}_{\cS_\del(\Om)}\le C_0 \abs{u}_{H^1(\Om)},\quad \forall u\in H^1_0(\Om).
\end{equation}
Here  we identify functions in \(H^1_0(\Om)\) with their zero extension to \(\mathbb R^d\) so that they can be naturally compared with functions in the nonlocal space \(\cS_{\del,0}(\Om)\).

The local limit of \eqref{eq:fullSpace:BVP:strong} as \(\del\to 0\) is the local Poisson problem \eqref{eq:local_Poisson_strong} with the unique solution  \(U_0\).

To apply the abstract framework in \Cref{sec:framework}, we define
\[
    \cU_\del=\cV_\del:=\cS_{\del,0}(\Om),
    \qquad \del\in (0,\del_0]
\]
equipped with the full \(\cS_\del(\Om)\) norm.  Set
\(\cU_0=\cV_0:=H^1_0(\Om)\), equipped with the full norm
\[
    \norm{u}_{\cU_0}
    :=\left(
        \norm{u}_{L^2(\Om)}^2
        +\norm{\nabla u}_{L^2(\Om,\mathbb R^d)}^2
    \right)^{1/2}.
\]
Without loss of generality, we assume \(C_P\) is also the Poincar\'e constant for the local space \(H^1_0(\Om)\), i.e., \(\norm{u}_{L^2(\Om)}\le C_P\norm{\nabla u}_{L^2(\Om,\mathbb R^d)}\) for \(u\in H^1_0(\Om)\).
For \(\del>0\), define \(a_\del\) on
\(\cS_\del(\Om)\times\cS_\del(\Om)\), and define \(a_0\) on
\(H^1(\Om)\times H^1(\Om)\), by
\begin{equation}\label{eq:fullSpace:bilinearForm}
    a_\del(u,v) = \left\{\begin{aligned}
        &\displaystyle
        \iint_{\Omh\times\Omh}
        \ga_\del(\xb-\yb)
        \bigl(u(\xb)-u(\yb)\bigr)
        \bigl(v(\xb)-v(\yb)\bigr)\,d\yb\,d\xb,
        \quad (\del>0),\\
        &\int_{\Om} \nabla u(\xb)\cdot \nabla v(\xb)\,d\xb,\quad (\del = 0).
    \end{aligned}\right.
\end{equation} 

It is clear that 
\[
\abs{a_\del(u,v)}
    \le
    \abs{u}_{\cU_\del}\abs{v}_{\cV_\del}
    \le \norm{u}_{\cU_\del}\norm{v}_{\cV_\del},\quad u,v\in \cS_\del(\Om),
\]
and for $u\in \cS_{\del,0}(\Om)$,
\[
a_\del(u,u)=\abs{u}_{\cU_\del}^2\ge (1+C_P^2)^{-1}\norm{u}_{\cU_\del}^2.
\]
Thus \Cref{assu:ellip} holds with \(M_*=1\), \(M_0=(1+C_0^2)^{1/2}\), \(M_a=1\) and \(\al = (1+C_P^2)^{-1}\).
Set \(u_\del:=U_\del-G|_{\Omh}\in \cU_\del\) for \(\del\in[0,\del_0]\). Then the weak formulations of
\eqref{eq:fullSpace:BVP:strong} and \eqref{eq:local_Poisson_strong} can be
written as follows:

\begin{equation}\label{eq:fullSpace:Poisson:weak}
    \left\{\begin{aligned}
        &\text{Find }u_\del \in \cU_\del \text{ such that }\\
        &a_\del(u_\del ,v)=(f,v)_{L^2(\Om)}-a_\del(G,v),
        \qquad \forall v\in \cV_\del.
    \end{aligned}\right.
\end{equation} 

Note that \(u_0=U_0-G|_{\Om}\) is the solution to the local Poisson problem with homogeneous boundary conditions. Let \(E_0u_0\) denote the zero extension of \(u_0\), and set \(\Uext:=E_0u_0+G\). Thus
\begin{equation*}
    \Uext(\xb)=\begin{cases}
        U_0(\xb), & \xb\in\Om,\\
        G(\xb), & \xb\in\Omega^+_{\del_0}.
    \end{cases}
\end{equation*}
For each \(\del\in (0,\del_0]\), the same notation \(\Uext\) denotes the
restriction of this fixed glued extension to \(\widehat\Om_\del\).
Apply the abstract framework to the homogeneous problem \eqref{eq:fullSpace:Poisson:weak}. Then, by
the definitions of \(\edel\) and \(\Rdel\),
\(\edel=u_\del -E_0u_0 =U_\del-\Uext\), and the corresponding
residual is
\begin{equation*}
\begin{aligned}
    \cR_\del(v)
    &=(f,v)_{L^2(\Om)}-a_\del(G,v)-a_\del(E_0u_0,v)\\
    &=(f,v)_{L^2(\Om)}-a_\del(\Uext,v).
\end{aligned}
\end{equation*}

Assume that \(U_0\in H^2(\Om)\). Then \(-\Delta U_0 =f\) by \eqref{eq:local_Poisson_strong}. Thus 
\begin{equation}\label{eq:Rdel_formula}
    \Rdel(v) = \int_{\Omega} -\Delta U_0(\xb)v(\xb)\,d\xb - \iint_{\Omh\times\Omh} \ga_\del(\xb-\yb)(\Uext(\xb)-\Uext(\yb))(v(\xb) - v(\yb))\,d\yb\,d\xb.
\end{equation}

We will verify the assumptions in \Cref{sec:framework} and apply
\Cref{thm:abstract-energy-rate,thm:abstract-L2-rate} to obtain the convergence rates.

\subsection{Continuous convergence rates}\label{sec:main-results}
We first state the continuous convergence rates in the energy and \(L^2\)
norms.
Analogous to discretizations for local problems, higher orders of convergence are guaranteed by sufficient regularity of the solution, data, and domain. To this end, we assume the following for some $m \in \mathbb{N}$:

\begin{assumption}\label{assu:compatibility}
    $\Omega$ is a $C^{m,1}$ bounded domain, \(U_0\in H^{m+1}(\Om)\), and \(g \in H^{m+1/2}(\p \Om)\) with extension $G \in H^{m+1}(\widehat\Om_{\del_0})$ and
    \begin{equation}\label{eq:inhom-trace-compatibility}
    \Tr(\normal^j U_0)=\Tr(\normal^j G)
    \quad \text{on }\partial\Om,
    \qquad j=0,1,\dots,m-1.
    \end{equation}
\end{assumption}
The regularity of $g$ assumed guarantees the existence of a Sobolev extension $G$ with normal traces $\Tr(\normal^j G) \in H^{m+1/2-j}(\p \Om)$ for $j$ up to order $m$; see \Cref{thm:normal-traces}. In addition, there exists a constant \(C=C(m,\Om,\del_0)>0\) such that
\begin{equation*}
    \norm{G}_{H^{m+1}(\widehat\Om_{\del_0})}
    \le C\norm{g}_{H^{m+1/2}(\p\Om)}.
\end{equation*}
The boundary compatibility \eqref{eq:inhom-trace-compatibility} is not guaranteed to be satisfied by the solution $U_0$, hence its inclusion in the assumption. See the appendix for the precise definitions of the boundary trace operator and normal derivatives.

Note that \eqref{eq:inhom-trace-compatibility} always holds with \(m=1\) due to the boundary condition \(\Tr U_0=g=\Tr G\) on \(\p\Om\).

\begin{remark}\label{rem:LapI-sufficient-data-regularity}
Let \(m\ge1\) be an integer and suppose that \(\Om\) is a bounded
\(C^{m,1}\) domain.  If
\[
    f\in H^{m-1}(\Om),
    \qquad
    g\in H^{m+1/2}(\p\Om),
\]
then the regularity of the solution to \eqref{eq:local_Poisson_strong} (see \cite[Section~2.5.1]{grisvard2011elliptic}) gives
\begin{equation*}
    \norm{U_0}_{H^{m+1}(\Om)}
    \le C_{\rm reg}\left(
        \norm{f}_{H^{m-1}(\Om)}
        +\norm{g}_{H^{m+1/2}(\p\Om)}
    \right)
\end{equation*}
for some constant \(C_{\rm reg}>0\).
Consequently, with \(u_0=U_0-G|_\Om\), using the lifting bound above,
\begin{equation*}
    \norm{u_0}_{H^{m+1}(\Om)}
    +\norm{G}_{H^{m+1}(\widehat\Om_{\del_0})}
    \le C\left(
        \norm{f}_{H^{m-1}(\Om)}
        +\norm{g}_{H^{m+1/2}(\p\Om)}
    \right).
\end{equation*}
Therefore, in all the convergence rate theorems below one may replace the bound \(\norm{u_0}_{H^{m+1}(\Om)}
    +\norm{G}_{H^{m+1}(\widehat\Om_{\del_0})}\) by the right-hand side data \(\norm{f}_{H^{m-1}(\Om)}
        +\norm{g}_{H^{m+1/2}(\p\Om)}\).
\end{remark}

Once the compatibility order \(m\) has been fixed, choose whole-space Sobolev
extensions \(\bar U_0\) of \(U_0\) and \(\bar G\) of \(G\):
\begin{equation*}
\begin{aligned}
    \norm{\bar U_0}_{H^{m+1}(\mathbb R^d)}
    &\le C\norm{U_0}_{H^{m+1}(\Om)},\\
    \norm{\bar G}_{H^{m+1}(\mathbb R^d)}
    &\le C\norm{G}_{H^{m+1}(\widehat\Om_{\del_0})},
\end{aligned}
\end{equation*}
where \(C=C(m,\Om,\del_0)\). 

The fixed choice of \(\del_0\) is
\begin{equation}\label{eq:def-delta-0}
    \del_0:=\frac{1}{2}\min\{\rho_0,\rho_1,\del_P\},
\end{equation}
where 
\(\rho_0\) and \(\rho_1\) are given in \Cref{lem:boundary-layer-est} and \eqref{eq:def_rho_1}, and \(\del_P\) is given by the uniform Poincar\'e inequality above.
The bound \(\del_0<\rho_1\) ensures that the fixed outer domain
\(\widehat\Omega_{\del_0}\) is Lipschitz \cite[Corollary~2]{lewicka2020which} and admits bounded Sobolev
extension operators.

\begin{theorem}\label{thm:LapI-rates-continuous}
Let \(u_\del\in \cU_\del\) solve \eqref{eq:fullSpace:Poisson:weak}, and let \(E_0u_0\) be the zero extension of \(u_0\). Suppose that \Cref{assu:compatibility} holds for some \(m\in\N_+\).
Then there exists a constant \(C=C(m,d,\Om,\ga,\del_0)>0\) such that for every \(\del\in (0,\del_0]\),
\begin{equation}\label{eq:LapI-energy-rate}
    \norm{u_\del-E_0u_0}_{\cU_\del}
    \le
    C\del^{\min\{2,m-\frac{1}{2}\}}
    \left(
        \norm{u_0}_{H^{m+1}(\Om)}
        +\norm{G}_{H^{m+1}(\widehat\Omega_{\del_0})}
    \right),
\end{equation}
and
\begin{equation*}
    \norm{u_\del-u_0}_{L^2(\Om)} 
    \le
    C\del^{\min\{2,m\}}
    \left(
        \norm{u_0}_{H^{m+1}(\Om)}
        +\norm{G}_{H^{m+1}(\widehat\Omega_{\del_0})}
    \right).
\end{equation*}
\end{theorem}

\begin{remark}\label{rem:LapI-energy-rate-high-regularity}
    Note that the function $\Uext$ defined above belongs to $H^{m}(\widehat\Om_{\del_0})$ thanks to \Cref{assu:compatibility}.
    If we assume that \eqref{eq:inhom-trace-compatibility} holds additionally with \(j=m\), then \(\Uext\in H^{m+1}(\widehat\Om_{\del_0})\) and using the same argument, we can improve the convergence rate in \eqref{eq:LapI-energy-rate} to \(\del^{\min\{2, m\}}\), i.e.,
    \[
    \norm{u_\del-E_0u_0}_{\cU_\del}
    \le
    C\del^{\min\{2,m\}}\| \Uext \|_{H^{m+1}(\widehat\Om_{\del_0})}.
    \]
\end{remark}
\begin{remark}
The compatibility condition in \Cref{assu:compatibility} describes
the agreement between the exterior data and the local solution
through their boundary traces.
When \(m=1\), only the boundary values match, and
\Cref{thm:LapI-rates-continuous} gives \(O(\del)\) convergence
in \(L^2(\Om)\) under \(H^2\) regularity.
This is consistent with the first-order convergence in the maximum
norm obtained for constant extensions of Dirichlet values, which
generally fail to match the first normal derivative
\cite{lee2022second}.
When \(m=2\), the first normal derivatives also match, and the theorem
gives \(O(\del^2)\) convergence in \(L^2(\Om)\) under \(H^3\)
regularity.
This boundary agreement is achieved by linear extrapolation using
the boundary derivatives of the local solution, which yields
second-order convergence in the maximum norm for smooth
one-dimensional problems \cite{du2019uniform}.
Second-order boundary treatments have also been constructed using
nonlocal gradients \cite{lee2022second} and auxiliary local problems
\cite{delia2020physically}.
The present theorem quantifies the dependence on boundary compatibility
under Sobolev regularity and distinguishes the convergence rates
in the nonlocal energy and \(L^2\) norms.
\end{remark}

\subsection{Discrete Convergence rates}\label{subsec:LapI-discrete-results} 
Fix a polynomial degree \(k\ge1\). For sufficiently small
\(h_0\in (0,1)\), let \(\{\mathcal T_h\}_{0<h\le h_0}\) be a
shape-regular family of exact conforming triangulations
of \(\Om\), obtained by the
construction in \cite[Section~6.2]{bernardi1989optimal}, with \(\theta=m\). Although the construction is stated for \(C^{m+1}\) domains,
the interpolation and derivative estimates in its proof require
only \(W^{m+1,\infty}\) bounds on the boundary charts.
These bounds are available for \(C^{m,1}\) boundaries.
The resulting family is therefore regular of order \(m\) in the
sense of \cite[Definition~3.2]{bernardi1989optimal}, with
\(C^{m,1}\) element maps replacing \(C^{m+1}\) maps. 

For each \(K\in\mathcal T_h\), let
\(F_K:\widehat K\to K\) be the associated \(C^{m,1}\)
diffeomorphism from the fixed reference simplex
\(\widehat K\). Define
\[
 V_{h,0}^{k}
 :=
 \left\{
 v_h\in H^1_0(\Om):
 v_h|_K\circ F_K\in\P_k(\widehat K)
 \text{ for every }K\in\mathcal T_h
 \right\}.
\]
Let \(\cU_{\del,h}=\cV_{\del,h}\subset\cU_\del\) be any
conforming finite element space containing \(V_{h,0}^{k}\),
where functions on \(\Om\) are identified with their zero
extensions to \(\widehat\Omega_\del\). 

By
\cite[Theorem~5.1 and Corollary~5.2]{bernardi1989optimal}, with a slight generalization to \(C^{m,1}\) maps \(F_K\), there exist an operator
\(I_h:H^1_0(\Om)\to V_{h,0}^{k}\) and a constant \(C>0\) independent of \(h\) such that
\begin{equation}\label{eq:FE-H1-approximation}
 \norm{v-I_hv}_{H^1(\Om)}
 \le C h^s\norm{v}_{H^{s+1}(\Om)},
 \qquad
 v\in H^{s+1}(\Om)\cap H^1_0(\Om),
 \quad s=1,\ldots,\min\{m,k\}.
\end{equation}  

Let \(u_{\del,h}\) be the unique solution of the discrete nonlocal
Poisson problem:
\begin{equation}\label{eq:fullSpace:Poisson:weak:discrete}
    \left\{\begin{aligned}
        &\text{Find }u_{\del,h}\in \cU_{\del,h}\ \text{such that}\\
        &a_\del(u_{\del,h},v_h)=(f,v_h)_{L^2(\Om)}-a_\del(G,v_h),\ \forall v_h\in \cV_{\del,h}.
    \end{aligned}\right.
\end{equation}
The solution with boundary condition is \(U_{\del,h}:=u_{\del,h}+G\).
By \Cref{rem:coercive} and coercivity of \(a_\del\), \Cref{assu:ellip_discrete} holds with \(\tilde{\al} = (1+C_P^2)^{-1}\).

\begin{theorem}\label{thm:LapI-rates-discrete}
Let \(u_{\del,h}\) solve
\eqref{eq:fullSpace:Poisson:weak:discrete}, let
\(U_{\del,h}=u_{\del,h}+G\), and let
\(\Uext=E_0u_0+G\). Suppose that  \Cref{assu:compatibility} holds for some \(m\in\N_+\).  
 Then there exists a constant \(C=C(m,k,d,\Om,\ga,\del_0)>0\) such that for every \(\del\in (0,\del_0]\) and
\(0<h\le h_0\),
\begin{equation}\label{eq:LapI-energy-rate-discrete}
    \norm{u_{\del,h}-E_0u_0}_{\cU_\del}
    \le C\left(
        \del^{\min\{2,m-\frac{1}{2}\}}
        +h^{\min\{m,k\}}
    \right)
    \left(\norm{u_0}_{H^{m+1}(\Om)}+\norm{G}_{H^{m+1}(\widehat\Om_{\del_0})}\right),
\end{equation}
and
\begin{equation}\label{eq:LapI-L2-rate-discrete}
    \norm{u_{\del,h}-u_0}_{L^2(\Om)} 
    \le C\left(
        \del^{\min\{2,m\}}
        +(\del^{\frac{1}{2}}+h)h^{\min\{m,k\}}
    \right)
    \left(\norm{u_0}_{H^{m+1}(\Om)}+\norm{G}_{H^{m+1}(\widehat\Om_{\del_0})}\right).
\end{equation}
\end{theorem}
\begin{remark}
\label{rem:LapI-L2-separated}
    If \(m=1\) or \(m\ge 2\) and \(\min\{m,k\}\ge 3\), the \(L^2\) convergence rate in \eqref{eq:LapI-L2-rate-discrete} would be \(O(\del^{\min\{2,m\}}
        +h^{\min\{m,k\}+1})\), as the mixed term \(\del^{\frac{1}{2}}h^{\min\{m,k\}}\) can be absorbed into other terms using Young's inequality. Indeed, if \(m=1\), then \(\del^{\frac{1}{2}}h\le \tfrac12 (\del+h^2)\le \tfrac12 (\del+h^{\min\{m,k\}+1})\). Similarly, if \(m\ge 2\) and \(\min\{m,k\}\ge 3\), then \(\frac{4}{3}\min\{m,k\}\ge \min\{m,k\}+1\) and Young's inequality implies that
        \[
        \del^{\frac{1}{2}}h^{\min\{m,k\}}
        \le \frac{1}{4} \del^2 + \frac{3}{4} h^{\frac{4}{3}\min\{m,k\}}\le \frac{1}{4} \del^2 + \frac{3}{4} h^{\min\{m,k\}+1}.
        \]
\end{remark}

Let \(\mathcal T_h\) be as above and denote the piecewise constant finite element space by
\[
 V_h^0
 :=
 \left\{
 v_h\in L^2(\Om):
 v_h|_K\in\P_0(K)\ \text{for every }K\in\mathcal T_h
 \right\},
\]
with every \(v_h\in V_h^0\) extended by zero to
\(\widehat\Omega_\del\).  
\begin{proposition}
\label{prop:LapI-piecewise-constant}
Suppose that, for some
\(0\le\alpha <\frac{1}{2}\),
\begin{equation*}
 \int_{B_1(\bm0)}\abs{\zb}^{2\alpha}\ga(\zb)\,d\zb<\infty.
\end{equation*}
Then \(V_h^0\subset\cU_\del\). Let
\(u_{\del,h}\in V_h^0\) be the Galerkin solution
\eqref{eq:fullSpace:Poisson:weak:discrete} with
\(\cU_{\del,h}=\cV_{\del,h}=V_h^0\).  Suppose that  \Cref{assu:compatibility} holds for some \(m\in\N_+\).  Then there exists \(C=C(m,k,d,\Om,\ga,\al,\del_0)>0\) independent of \(\del\) and \(h\),
such that
\begin{equation*}
\begin{split}
 \norm{u_{\del,h}-E_0u_0}_{\cU_\del}
 \le C\left(
  \del^{\min\{2,m-\frac{1}{2}\}}
  +\left(\frac{h}{\del}\right)^{1-\alpha}
 \right)
 \left(
  \norm{u_0}_{H^{m+1}(\Om)}
  +\norm{G}_{H^{m+1}(\widehat\Om_{\del_0})}
 \right).
\end{split}
\end{equation*}
Moreover,
\begin{equation*}
 \norm{u_{\del,h}-u_0}_{L^2(\Om)}
 \le C\Bigg(
   \del^{\min\{2,m\}}
  +\del^{\frac{1}{2}}
   \left(\frac{h}{\del}\right)^{1-\alpha}
  +\left(\frac{h}{\del}\right)^{2-2\alpha}
 \Bigg) \left(
  \norm{u_0}_{H^{m+1}(\Om)}
  +\norm{G}_{H^{m+1}(\widehat\Om_{\del_0})}
 \right). 
\end{equation*}
\end{proposition}

\subsection{Proofs of the main results}\label{sec:proofs-main-results}
In this subsection we prove the main estimates. We first prove the energy rates. Then, applying the energy rates to the homogeneous problem and using a duality argument, we prove the \(L^2\) rates.

Throughout the proof, \(C\) denotes a constant independent of the horizon \(\del\), of the mesh parameter \(h\), and of the test functions. Its value may change from line to line and may depend on fixed data such as \(m,d,\Omega,\gamma\), \(\del_0\), and the mesh shape-regularity constant.

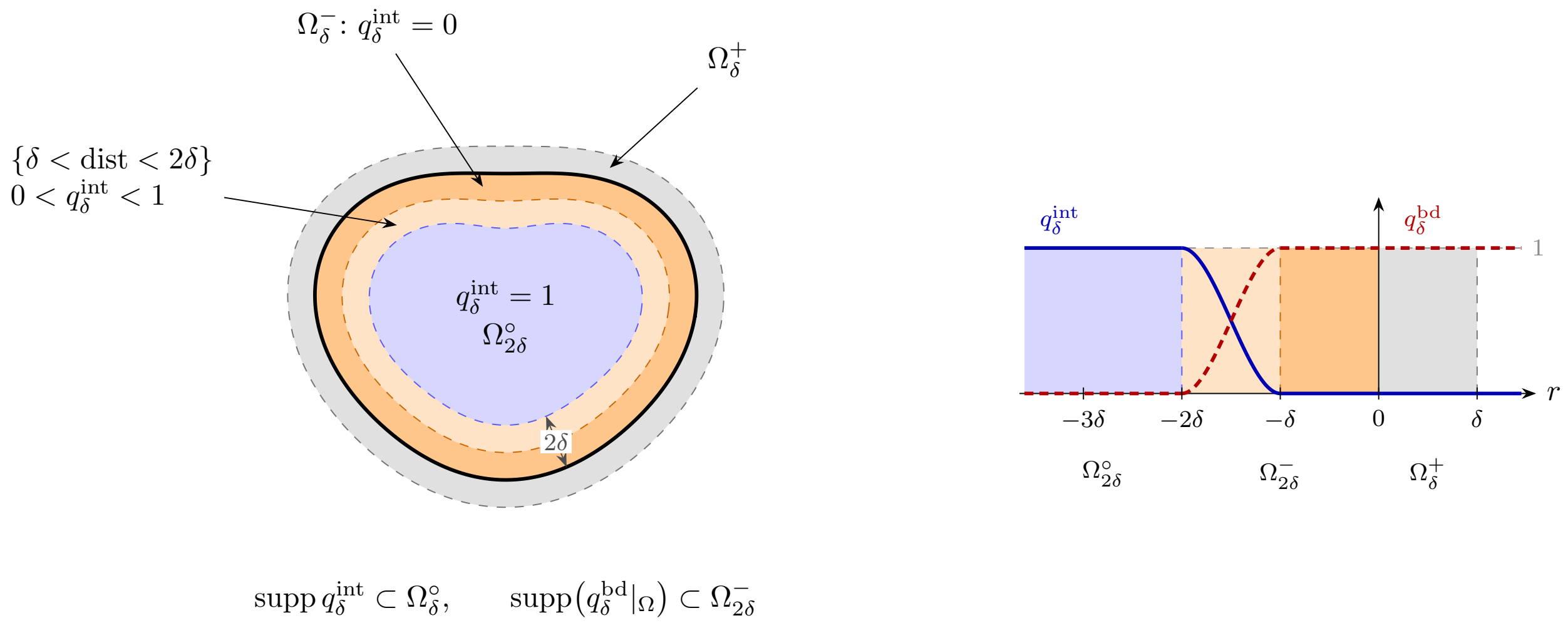
\begin{figure}[htbp]
\centering
\resizebox{\textwidth}{!}{%
\begin{tikzpicture}[>=Stealth,line join=round]
  \def\bl#1#2{2.0+0.21*cos(2*(#1))+0.12*sin(3*(#1))+(#2)}
  \def\dl{0.32}

  %------------------------------- panel (a): geometry
  \fill[black!12]  plot[variable=\t,domain=0:360,samples=200,smooth cycle] (\t:{\bl{\t}{\dl}});
  \fill[orange!45] plot[variable=\t,domain=0:360,samples=200,smooth cycle] (\t:{\bl{\t}{0}});
  \fill[orange!22] plot[variable=\t,domain=0:360,samples=200,smooth cycle] (\t:{\bl{\t}{-\dl}});
  \fill[blue!16]   plot[variable=\t,domain=0:360,samples=200,smooth cycle] (\t:{\bl{\t}{-2*\dl}});

  \draw[dashed,black!55]        plot[variable=\t,domain=0:360,samples=200,smooth cycle] (\t:{\bl{\t}{\dl}});
  \draw[very thick]             plot[variable=\t,domain=0:360,samples=200,smooth cycle] (\t:{\bl{\t}{0}});
  \draw[dashed,orange!80!black] plot[variable=\t,domain=0:360,samples=200,smooth cycle] (\t:{\bl{\t}{-\dl}});
  \draw[dashed,blue!65]         plot[variable=\t,domain=0:360,samples=200,smooth cycle] (\t:{\bl{\t}{-2*\dl}});

  \node[align=center] at (0,0) {$ q^{\mathrm{int}}_{\delta}=1$\\[1pt]$\Omega^{\circ}_{2\delta}$};

  \coordinate (A1) at (140:{\bl{140}{-1.5*\dl}});
  \coordinate (A2) at (100:{\bl{100}{-0.5*\dl}});
  \coordinate (A3) at (55:{\bl{55}{0.5*\dl}});
  \node[align=left] (LA1) at (-4.6,1.6) {$\{\delta<\operatorname{dist}<2\delta\}$\\ $0< q^{\mathrm{int}}_{\delta}<1$};
  \node[align=left] (LA2) at (-1.5,3.4) {$\Omega^{-}_{\delta}$:\ $ q^{\mathrm{int}}_{\delta}=0$};
  \node (LA3) at (2.6,3.0) {$\Omega^{+}_{\delta}$};
  \draw[->,thin] (LA1) -- (A1);
  \draw[->,thin] (LA2) -- (A2);
  \draw[->,thin] (LA3) -- (A3);

  \draw[<->,thin,black!70] (-68:{\bl{-68}{-2*\dl}}) --
        node[fill=white,inner sep=1pt,font=\footnotesize] {$2\delta$} (-68:{\bl{-68}{0}});

  \node[align=center] at (0,-3.3)
     {$\operatorname{supp} q^{\mathrm{int}}_{\delta}\subset\Omega^{\circ}_{\delta}$,\qquad
      $\operatorname{supp}\bigl( q^{\mathrm{bd}}_{\delta}|_{\Omega}\bigr)\subset\Omega^{-}_{2\delta}$};
  \node at (0,-3.9) {(a)};

 %------------------------------- panel (b): signed-distance profiles
  \begin{scope}[xshift=10.2cm,yshift=-0.9cm]
    \def\ux{1.15}
    \def\uy{1.7}

    \fill[blue!16]   (-3.6*\ux,0) rectangle (-2*\ux,\uy);
    \fill[orange!22] (-2*\ux,0)   rectangle (-1*\ux,\uy);
    \fill[orange!45] (-1*\ux,0)   rectangle (0,\uy);
    \fill[black!12]  (0,0)        rectangle (1*\ux,\uy);

    \draw[->] (-3.65*\ux,0) -- (1.6*\ux,0) node[right] {$r$};
    \draw[->] (0,0) -- (0,1.35*\uy);
    \draw[dashed,black!45] (-3.6*\ux,\uy) -- (1.45*\ux,\uy)
      node[right,font=\footnotesize] {$1$};
    \draw[dashed,blue!65] (-2*\ux,0) -- (-2*\ux,\uy);
    \draw[dashed,orange!80!black] (-1*\ux,0) -- (-1*\ux,\uy);
    \draw[dashed,black!55] (1*\ux,0) -- (1*\ux,\uy);

    % Interior cutoff: solid blue.
    \draw[very thick,blue!70!black]
      (-3.6*\ux,\uy) -- (-2*\ux,\uy);
    \draw[very thick,blue!70!black]
      plot[variable=\x,domain=-2:-1,samples=80]
      ({\x*\ux},{(1-3*(\x+2)^2+2*(\x+2)^3)*\uy});
    \draw[very thick,blue!70!black]
      (-1*\ux,0) -- (1.45*\ux,0);
    \node[blue!70!black,font=\footnotesize,anchor=south west]
      at (-3.55*\ux,1.02*\uy) {$q_\delta^{\mathrm{int}}$};

    % Boundary cutoff: dashed red.
    \draw[very thick,red!70!black,densely dashed]
      (-3.6*\ux,0) -- (-2*\ux,0);
    \draw[very thick,red!70!black,densely dashed]
      plot[variable=\x,domain=-2:-1,samples=80]
      ({\x*\ux},{(3*(\x+2)^2-2*(\x+2)^3)*\uy});
    \draw[very thick,red!70!black,densely dashed]
      (-1*\ux,\uy) -- (1.45*\ux,\uy);
    \node[red!70!black,font=\footnotesize,anchor=south west]
      at (0.15*\ux,1.02*\uy) {$q_\delta^{\mathrm{bd}}$};

    \foreach \v/\l in
      {-3/{-3\delta},-2/{-2\delta},-1/{-\delta},0/{0},1/{\delta}}{
      \draw (\v*\ux,2pt) -- (\v*\ux,-2pt)
        node[below,font=\footnotesize] {$\l$};
    }

    \node[font=\footnotesize] at (-2.8*\ux,-0.55*\uy)
      {$\Omega^{\circ}_{2\delta}$};
    \node[font=\footnotesize] at (-1*\ux,-0.55*\uy)
      {$\Omega^{-}_{2\delta}$};
    \node[font=\footnotesize] at (0.5*\ux,-0.55*\uy)
      {$\Omega^{+}_{\delta}$};

    \node at (-1*\ux,-3.0) {(b)};
  \end{scope}
\end{tikzpicture}}
\caption{The type~(I) cutoffs.
(a) The interior region, transition layer, and exterior collar.
(b) Profiles of \(q_\delta^{\mathrm{int}}\) (solid) and
\(q_\delta^{\mathrm{bd}}\) (dashed) against the signed distance
\(r\) to \(\partial\Omega\), negative inside \(\Omega\)
and positive outside.}
\label{fig:cutoff-typeI}
\end{figure}

The residual \(\Rdel\) is given in \eqref{eq:Rdel_formula}, and has two different sources. Away from the boundary, after a suitable extension of the solution to the whole space, the error is essentially the consistency error of \(\cL_\del\), which is estimated by Fourier analysis in \Cref{subsec:interior}. Near the boundary, the error is analyzed by the  boundary-layer estimates in \Cref{subsec:boundary}. To separate these two effects, we introduce a smooth cutoff \( q_\del^{\mathrm{int}}\in C^\infty_c(\mathbb R^d)\) such that
\[
    0\le  q_\del^{\mathrm{int}}\le 1,\qquad
    \supp q_\del^{\mathrm{int}}\subset \Omin,\qquad
     q_\del^{\mathrm{int}}=1\ \text{on }\Ominn,\qquad
    |\nabla q_\del^{\mathrm{int}}|\le C\del^{-1}.
\]
Set \( q_\del^{\mathrm{bd}}:=1- q_\del^{\mathrm{int}}\). 
The two cutoffs and their supports are shown in \Cref{fig:cutoff-typeI}. 

% =====================================================================
%  FIGURE 2 --- type (I) cutoff.   Section 4.2, at the definition of the cutoff.
% =====================================================================

Recall that the residual $\Rdel$ is defined in \eqref{eq:Rdel_formula}. We introduce the interior and boundary residuals \(\Rint(v):=\Rdel( q_\del^{\mathrm{int}} v)\) and \(\Rbd(v):=\Rdel( q_\del^{\mathrm{bd}} v)\) for \(v\in \cV_\del\). By linearity, \(\Rdel(v)=\Rint(v)+\Rbd(v)\).
We write both terms explicitly and further decompose $\Rbd(v)$:
\begin{equation*}
    \Rint(v)=\int_{\Omega} -\Delta U_0(\xb) q_\del^{\mathrm{int}}(\xb)v(\xb)\,d\xb - \iint_{\Om\times\Om} \ga_\del(\xb-\yb)(U_0(\xb)-U_0(\yb))( q_\del^{\mathrm{int}}(\xb)v(\xb) -  q_\del^{\mathrm{int}}(\yb)v(\yb))\,d\yb\,d\xb,
\end{equation*}
and 
\begin{equation}\label{eq:Rbd}
    \Rbd(v)=\underbrace{(-\Delta u_0,  q_\del^{\mathrm{bd}} v)_{L^2(\Omega)} - a_\del(E_0 u_0,  q_\del^{\mathrm{bd}} v)}_{\RbdI(v)} + \underbrace{(-\Delta G,  q_\del^{\mathrm{bd}} v)_{L^2(\Omega)} - a_\del( G,  q_\del^{\mathrm{bd}} v)}_{\RbdII(v)}. 
\end{equation}
We will estimate \(\Rint\) and \(\RbdII\) using the whole-space consistency error via Fourier analysis, and estimate \(\RbdI\) using boundary-layer estimates.

\subsubsection{Interior residual estimates}\label{subsec:interior}
We now estimate the interior residual by a Fourier symbol
argument. Define
\begin{equation}\label{eq:full-space-symbol}
    \lambda_\del(\xib)
    :=2\int_{\R^d}\ga_\del(\zb)
    \bigl(1-\cos(2\pi\xib\cdot\zb)\bigr)\,d\zb.
\end{equation}
The Fourier representation of the energy and the symbol estimates
used below are established in
\Cref{lem:fourier-bond-energy,lem:bond-symbol-estimates}
in Appendix~\ref{sec:halfspace-analytic-foundations}.

\begin{lemma}\label{lem:interior-consistency}
Let \(j\in\N_+\). There exists \(C=C(j,d,\ga)>0\)  such that for every \(q\in H^{1+j}(\mathbb R^d)\) and every \(\psi\in \cS_\del(\mathbb R^d)\),
\begin{equation*}
\begin{aligned}
    &\left|\int_{\R^d}\Delta q(\xb)\psi(\xb)\,d\xb
    +\iint_{\mathbb R^d\times\R^d}
       \ga_\del(\xb-\yb)(q(\xb)-q(\yb))
       (\psi(\xb)-\psi(\yb))\,d\xb\,d\yb\right|\\
    &\qquad\le
    C\del^{\min\{2,j\}}\norm{q}_{H^{1+j}(\mathbb R^d)}
    \abs{\psi}_{\cS_\del(\mathbb R^d)}.
\end{aligned}
\end{equation*}
\end{lemma}

\begin{proof} 
    Since \(\widehat{\Delta q}(\xib) = -(2\pi|\xib|)^2\hat{q}(\xib)\), by Plancherel's theorem and \Cref{lem:fourier-bond-energy} we have
    \begin{align*}
        &\quad\int_{\R^d} \Delta q(\xb)\psi(\xb)\,d\xb+\iint_{\mathbb R^d\times \R^d}\ga_\del(\xb-\yb)(q(\xb)-q(\yb))(\psi(\xb)-\psi(\yb))\,d\xb\,d\yb \\
        &= \int_{\R^d} \left(\lambda_\del(\xib) - (2\pi|\xib|)^2\right)\hat{q}(\xib)\overline{\hat{\psi}(\xib)}\,d\xib.
    \end{align*} 
    Note that 
    \[\norm{q}_{H^{1+j}(\mathbb R^d)}^2 \asymp \int_{\R^d} (1+(2\pi|\xib|)^2)^{1+j} |\hat{q}(\xib)|^2 \, d\xib,\qquad \abs{\psi}_{\cS_\del(\mathbb R^d)}^2 = \int_{\R^d} \lambda_\del(\xib) |\hat{\psi}(\xib)|^2 \, d\xib.\]
    Here \(A\asymp B\) means \(cB\le A\le CB\) for some constants \(c,C>0\).
    To prove \Cref{lem:interior-consistency}, it suffices to prove the multiplier estimate
    \begin{equation}\label{eq:multiplier-estimate}
        \abs{\lambda_\del(\xib) - (2\pi|\xib|)^2} \le C\del^{\min\{2,j\}}(1+(2\pi|\xib|)^2)^{\frac{1+j}{2}}\lambda_\del(\xib)^{1/2}.
    \end{equation}
    Indeed, if \eqref{eq:multiplier-estimate} holds, then by the Cauchy--Schwarz inequality,
    \begin{align*}
        &\quad\abs{\int_{\R^d} \left(\lambda_\del(\xib) - (2\pi|\xib|)^2\right)\hat{q}(\xib)\overline{\hat{\psi}(\xib)}\,d\xib} \\
        &\le C\del^{\min\{2,j\}}\left(\int_{\R^d} (1+(2\pi|\xib|)^2)^{1+j} |\hat{q}(\xib)|^2 \, d\xib\right)^{1/2} \left(\int_{\R^d} \lambda_\del(\xib) |\hat{\psi}(\xib)|^2 \, d\xib\right)^{1/2} \\
        &\le C\del^{\min\{2,j\}}\norm{q}_{H^{1+j}(\mathbb R^d)}\abs{\psi}_{\cS_\del(\mathbb R^d)}.
    \end{align*}
    It remains to prove \eqref{eq:multiplier-estimate}.
By \eqref{eq:bond-symbol-scaling}, it suffices to show that
\begin{equation}\label{eq:multiplier-estimate-scaled}
\begin{aligned}
    \bigl|\lambda_\ga(\xib)-(2\pi)^2|\xib|^2\bigr|
    &\le C\del^{1+\min\{2,j\}}
    \bigl(1+(2\pi)^2\del^{-2}|\xib|^2\bigr)^{(1+j)/2}
    \lambda_\ga(\xib)^{1/2}.
\end{aligned}
\end{equation}

For $|\xib|\le1$, \Cref{lem:bond-symbol-estimates} gives
\[
    \bigl|\lambda_\ga(\xib)-(2\pi)^2|\xib|^2\bigr|
    \le C|\xib|^4,
    \qquad \lambda_\ga(\xib)^{1/2}\ge c|\xib|.
\]
If $j=1$, the expression on the right of
\eqref{eq:multiplier-estimate-scaled}, without its constant, satisfies
\[
    \del^2\bigl(1+(2\pi)^2\del^{-2}|\xib|^2\bigr)
    \lambda_\ga(\xib)^{1/2}
    \ge c|\xib|^3\ge c|\xib|^4.
\]
If $j\ge2$, use
\begin{equation*}
    x^3\le(1+x^2)^{(j+1)/2},\qquad x\ge0,\quad j\ge2,
\end{equation*}
with $x=2\pi\del^{-1}|\xib|$ to obtain
\[
    \del^3\bigl(1+(2\pi)^2\del^{-2}|\xib|^2\bigr)^{(1+j)/2}
    \lambda_\ga(\xib)^{1/2}
    \ge c|\xib|^4.
\]
Thus \eqref{eq:multiplier-estimate-scaled} holds for $|\xib|\le1$.

For $|\xib|>1$, the bounds in the same lemma, with $\del=1$, give
\[
    \bigl|\lambda_\ga(\xib)-(2\pi)^2|\xib|^2\bigr|
    \le C|\xib|^2,
    \qquad \lambda_\ga(\xib)^{1/2}\ge c.
\]
Since $0<\del\le\del_0\le1$, we have
\[
\begin{aligned}
    &\del^{1+\min\{2,j\}}
    \bigl(1+(2\pi)^2\del^{-2}|\xib|^2\bigr)^{(1+j)/2}
    \lambda_\ga(\xib)^{1/2}\\
    &\qquad\ge c\del^{\min\{2,j\}-j}|\xib|^{1+j}
    \ge c|\xib|^2.
\end{aligned}
\]
This proves \eqref{eq:multiplier-estimate-scaled}, and hence
\eqref{eq:multiplier-estimate}.
\end{proof}

\begin{lemma}\label{lem:Rint-bound}
Let \(U_0\in H^{m+1}(\Om)\) and \(\Rint\) be the interior residual defined above. Then
\begin{equation*}
    |\Rint(v)|
    \le
    C\del^{\min\{2,m\}}\norm{U_0}_{H^{m+1}(\Om)}
    \norm{v}_{\cV_\del},
    \qquad \forall v\in \cV_\del.
\end{equation*}
\end{lemma}
\begin{proof} 
We use Sobolev extension to pass from the domain to the whole space and then apply \Cref{lem:interior-consistency}. Take a Sobolev extension \(\bar U_0\in H^{1+m}(\mathbb R^d)\) of \(U_0\). For \(v\in\cV_\del\), set \(\psi= q_\del^{\mathrm{int}}E_0v\). On \(\supp q_\del^{\mathrm{int}}\), the ball \(B_\del(\xb)\) is contained in \(\Om\), so
\[
    \abs{\Rint(v)}
    =\abs{\int_{\R^d} \Delta \bar U_0(\xb)\psi(\xb)\,d\xb+\iint_{\mathbb R^d\times \R^d}\ga_\del(\xb-\yb)(\bar U_0(\xb)-\bar U_0(\yb))(\psi(\xb)-\psi(\yb))\,d\xb\,d\yb}.
\]
Applying \Cref{lem:interior-consistency} with \(j=m\), the extension estimate \(\norm{\bar U_0}_{H^{m+1}(\R^d)}\le C\norm{U_0}_{H^{m+1}(\Om)}\), and the first estimate in \Cref{lem:thin-layer-cutoff-stability} with \( q_\del= q_\del^{\mathrm{int}}\) proves \Cref{lem:Rint-bound}.
\end{proof} 

\begin{proof}[Proof of \Cref{rem:LapI-energy-rate-high-regularity}]
    The same idea of using Sobolev extension for \(\Uext\) from \(\widehat\Om_{\del_0}\) to \(\R^d\) as in \Cref{lem:Rint-bound} and applying \Cref{lem:interior-consistency} yields 
\[\abs{\Rdel(v)}
    \le
    C\del^{\min\{2,m\}}\norm{\Uext}_{H^{m+1}(\widehat\Om_{\del_0})}
    \norm{v}_{\cV_\del},
    \qquad \forall v\in \cV_\del.
\]
Then \Cref{rem:LapI-energy-rate-high-regularity} holds by \Cref{thm:abstract-energy-rate}.
\end{proof}

Note that the second term \(\RbdII(v)\) in \eqref{eq:Rbd} can be estimated in a similar way as in the proof of \Cref{lem:Rint-bound}. We record the result here for later use.
\begin{equation}\label{eq:RbdII-bound}
    \abs{\RbdII(v)}\le C\del^{\min\{2,m\}}\norm{G}_{H^{m+1}(\widehat\Om_{\del_0})}
    \norm{v}_{\cV_\del},
    \qquad \forall v\in \cV_\del.
\end{equation}
Indeed, taking \(\psi_\del =  q_\del^{\mathrm{bd}} E_0 v\) and applying \Cref{lem:interior-consistency} with \(j=m\) yields the bound \eqref{eq:RbdII-bound}.

\subsubsection{Boundary residual and energy estimates}\label{subsec:boundary}

The boundary residual will be tested against functions in \(\cV_\del=\cS_{\del,0}(\Om)\), which vanish in the exterior collar. We prove two auxiliary results in the following. The first estimate is the nonlocal Poincar\'e inequality in a thin layer; it shows that the $L^2$ norm on a boundary layer with thickness proportional to $\del$ can be controlled by $\del$ multiplied by the nonlocal energy seminorm.

Let \(\Om\) be a \(C^{1,1}\) domain. By \cite[Corollary~2]{lewicka2020which}, there exists \(\rho_1>0\) depending on \(\Om\) such that 
\begin{equation}\label{eq:def_rho_1}
    B_{\rho_1}\bigl(\pb+\rho_1\nb(\pb)\bigr)
    \subset \mathbb R^d\setminus\Om,
    \qquad \forall\pb\in\partial\Om.
\end{equation} 
\begin{proposition}\label{prop:nonlocal-boundary-layer}
Let \(\Om\subset \mathbb R^d\) be a bounded \(C^{1,1}\) domain and \(\rho_1\) 
be as in \eqref{eq:def_rho_1}, and let
\(\ga\) satisfy 
\Cref{assu:kernel}. For every fixed  \(a\ge1\), there
exists \(C_a=C(a,d,\Om,\ga)>0\)  such that
for
every \(\del\in (0,4\rho_1)\), 
\[
   \norm{v}_{L^2(\Omega^-_{a\del})}
    \le C_a\del\abs{v}_{\cS_\del(\Om)},
    \qquad \forall v\in\cS_{\del,0}(\Om).
\]
\end{proposition}
% =====================================================================
%  FIGURE 4 --- covering and exterior-ball geometry (Prop. 4.6).
% =====================================================================
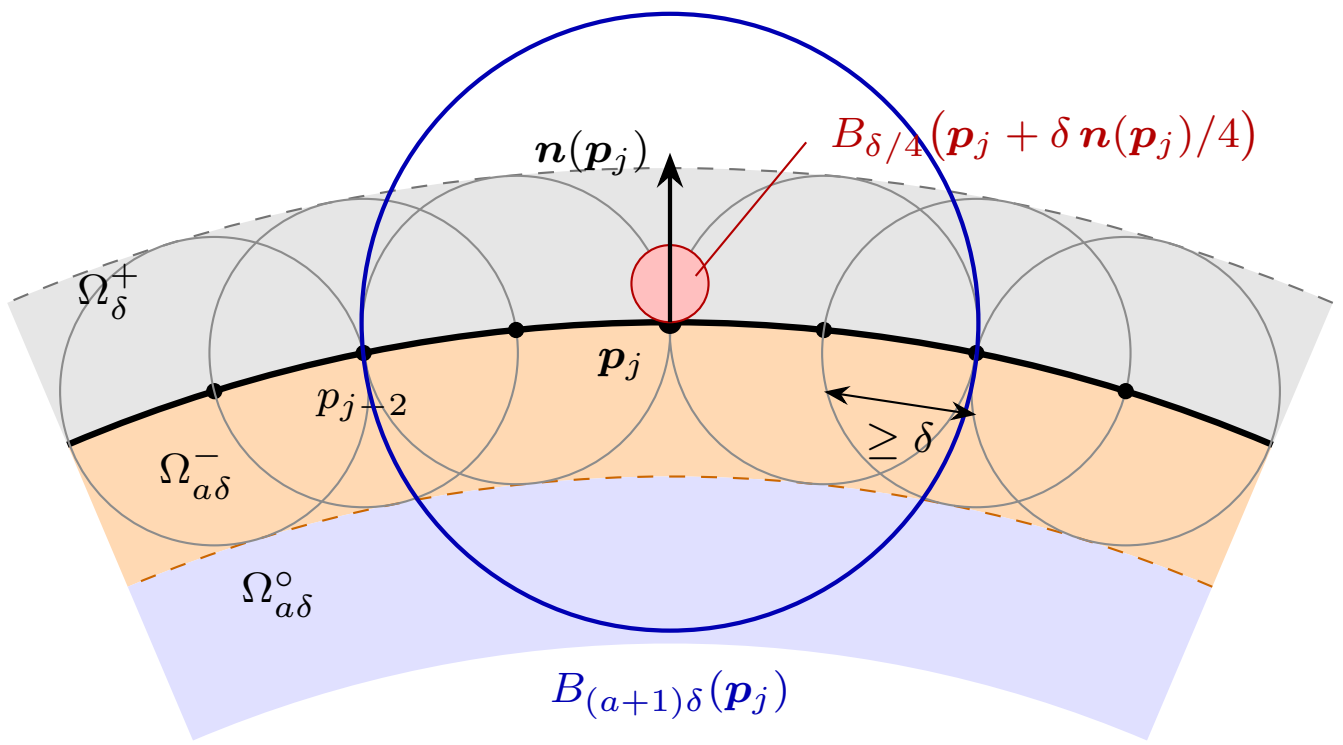
\begin{figure}[htbp]
\centering
\begin{tikzpicture}[>=Stealth,line join=round,scale=1.7]
  \def\dl{0.60}                        % \delta in picture units
  % \p\Om is the arc of the circle of radius 6 centred at (0,-6); the
  % angular window 67..113 degrees gives |x| <~ 2.34.  Bands are drawn as
  % closed arc pairs so that the picture stays bounded.
  \def\band#1#2#3{\fill[#3] ({#1*cos(67)},{-6+#1*sin(67)}) arc (67:113:#1)
      -- ({#2*cos(113)},{-6+#2*sin(113)}) arc (113:67:#2) -- cycle;}

  \band{6.6}{6}{black!10}
  \band{6}{5.4}{orange!30}
  \band{5.4}{4.75}{blue!12}

  \draw[dashed,black!55]        ({6.6*cos(67)},{-6+6.6*sin(67)}) arc (67:113:6.6);
  \draw[dashed,orange!80!black] ({5.4*cos(67)},{-6+5.4*sin(67)}) arc (67:113:5.4);
  \draw[very thick]             ({6*cos(67)},{-6+6*sin(67)})     arc (67:113:6);

  % maximal \del-separated set on \p\Om  (arc spacing \del <-> 5.73 degrees)
  \foreach \k in {-3,-2,-1,1,2,3}{
     \draw[black!45,thin] ({6*sin(5.73*\k)},{-6+6*cos(5.73*\k)}) circle (\dl);
     \fill[black] ({6*sin(5.73*\k)},{-6+6*cos(5.73*\k)}) circle (0.9pt);}
  \fill[black] (0,0) circle (1.3pt);
  \node[font=\footnotesize,anchor=north east] at (-0.02,-0.03) {$\pb_j$};
  \node[font=\footnotesize,anchor=north] at ({6*sin(-11.46)},{-6+6*cos(-11.46)-0.07}) {$p_{j-2}$};

  % B_{R\del}(p_j),  R = a+1 = 2
  \draw[thick,blue!70!black] (0,0) circle ({2*\dl});
  \node[blue!70!black,font=\footnotesize,anchor=north] at (0,-1.28) {$B_{(a+1)\delta}(\pb_j)$};

  % tangent exterior ball,  \rho = 1/4
  \fill[red!25] (0,{0.25*\dl}) circle ({0.25*\dl});
  \draw[red!70!black] (0,{0.25*\dl}) circle ({0.25*\dl});
  \node[red!70!black,font=\footnotesize,anchor=west] at (0.55,0.72)
       {$B_{\delta/4}\bigl(\pb_j+\delta\,\boldsymbol{n}(\pb_j)/4\bigr)$};
  \draw[red!70!black,thin] (0.53,0.70) -- (0.10,{0.30*\dl});

  \draw[->,thick] (0,0) -- (0,0.66) node[left,font=\footnotesize] {$\boldsymbol{n}(\pb_j)$};

  \draw[<->,thin] ({6*sin(5.73)},{-6+6*cos(5.73)-0.24}) --
        node[below,font=\footnotesize] {$\ge\delta$} ({6*sin(11.46)},{-6+6*cos(11.46)-0.24});

  \node[font=\footnotesize,anchor=east] at (-1.98,0.12)  {$\Omega^{+}_{\delta}$};
  \node[font=\footnotesize,anchor=east] at (-1.62,-0.58) {$\Omega^{-}_{a\delta}$};
  \node[font=\footnotesize,anchor=east] at (-1.30,-1.06) {$\Omega^{\circ}_{a\delta}$};
\end{tikzpicture}
\caption{Geometry of the proof of
Proposition~\ref{prop:nonlocal-boundary-layer}, drawn for $a=1$. }
\label{fig:boundary-layer-covering}
\end{figure}

\begin{proof}

The geometry is shown in \Cref{fig:boundary-layer-covering}. We first
record the only fixed-scale estimate that will be used. Let
\[
    B:=B_{a+1}(\boldsymbol 0),
    \qquad
    B':=B_{1/4}(\mathbf e_d/4)\subset B.
\]
By \eqref{eq:gamma-lower-bound} and the fixed-horizon nonlocal Poincar\'e
inequality \cite[Corollary~1]{mengesha2013analysis}, there exists
\(C_*=C(a,d,\ga)>0\) such that
\begin{equation}\label{eq:fixed-ball-nonlocal-poincare}
    \int_B |w(\boldsymbol z)|^2\,d\boldsymbol z
    \le C_*
    \iint_{B\times B}
    \ga(\boldsymbol z-\boldsymbol\zeta)
    \bigl(w(\boldsymbol z)-w(\boldsymbol\zeta)\bigr)^2
    \,d\boldsymbol\zeta\,d\boldsymbol z
\end{equation}
whenever \(w=0\) a.e. on \(B'\).  

By \eqref{eq:def_rho_1}, for \(\del\in (0,4\rho_1)\),
\begin{equation}\label{eq:uniform-exterior-small-ball}
    B_{\del/4}\bigl(\pb+\del\nb(\pb)/4\bigr)
    \subset \Om^c\cap B_{(a+1)\del}(\pb),
    \qquad \forall\pb\in\partial\Om.
\end{equation}

Choose a finite set
\(\{\pb_j\}_{j=1}^{N_\del}\subset\partial\Om\) such that
\[
    \partial\Om\subset\bigcup_{j=1}^{N_\del}B_\del(\pb_j),
    \qquad
    |\pb_i-\pb_j|\ge\del\quad(i\ne j).
\]
Write
\[
    B_{j,\del}:=B_{(a+1)\del}(\pb_j),
    \qquad
    D_{j,\del}:=\Omh\cap B_{j,\del}.
\]
The balls \(B_{j,\del}\)
have overlap bounded by a number \(M=M(a,d)\),
independently of \(N_\del\) and \(\del\), i.e.,
\begin{equation}\label{eq:bounded-overlap}
    \sum_{j=1}^{N_\del}\mathbf 1_{B_{(a+1)\del}(\pb_j)}(\xb)
    \le M(a,d),
    \qquad \xb\in\R^d.
\end{equation}
Indeed, for fixed \(\xb\in\R^d\),
if \( J(\xb):=\{j:\xb\in B_{(a+1)\del}(\pb_j)\}\), 
then \(\pb_j\in B_{(a+1)\del}(\xb)\) for \(j\in J(\xb)\) and  the balls
\(\{B_{\del/2}(\pb_j):j\in J(\xb)\}\) are pairwise disjoint and are contained in
\(B_{(a+\frac{3}{2})\del}(\xb)\).  Comparing volumes yields \eqref{eq:bounded-overlap}.
Let \(E_0v\) denote the zero extension of  \(v\in\cS_{\del,0}(\Om)\). For each \(j\), choose an orthogonal matrix \(Q_j\) such
that \(Q_j\mathbf e_d=\nb(\pb_j)\), and set
\[
    w_j(\boldsymbol z):=(E_0v)(\pb_j+\del Q_j\boldsymbol z),
    \qquad \boldsymbol z\in B.
\]
By \eqref{eq:uniform-exterior-small-ball}, \(w_j=0\) on \(B'\). Applying
\eqref{eq:fixed-ball-nonlocal-poincare} and using the radiality of \(\ga\),
together with
\[
    \ga_\del\bigl(\del Q_j(\boldsymbol z-\boldsymbol\zeta)\bigr)
    =\del^{-d-2}\ga(\boldsymbol z-\boldsymbol\zeta),
\]
we obtain after change of variables
\begin{equation}\label{eq:local-scaled-hardy}
\begin{split}
    \int_{\Om\cap B_{j,\del}} |v(\xb)|^2\,d\xb
    &\le C_*\del^2
    \iint_{B_{j,\del}\times B_{j,\del}}
    \ga_\del(\xb-\yb)
    \bigl((E_0v)(\xb)-(E_0v)(\yb)\bigr)^2\,d\yb\,d\xb \\
    &= C_*\del^2
    \iint_{D_{j,\del}\times D_{j,\del}}
    \ga_\del(\xb-\yb)
    \bigl(v(\xb)-v(\yb)\bigr)^2\,d\yb\,d\xb,
\end{split}
\end{equation}
The last equality follows from \(\supp\ga_\del\subset \overline B_\del(0)\):
after zero extension, the integrand vanishes whenever one variable lies
outside \(\Omh\).

It remains to sum \eqref{eq:local-scaled-hardy}. 
Note that
\[
    \Omega^-_{a\del}
    \subset\bigcup_j\bigl(\Om\cap B_{j,\del}\bigr).
\]
Summing \eqref{eq:local-scaled-hardy} and using \eqref{eq:bounded-overlap} yields 
\[
    \norm{v}_{L^2(\Omega^-_{a\del})}^2
    \le C_*M\del^2\abs{v}_{\cS_\del(\Om)}^2.
\]
Taking square roots completes the proof.
\end{proof}
The second lemma will be used to localize the errors in the interior or near the boundary of $\Om$.
\begin{lemma}\label{lem:thin-layer-cutoff-stability}
Let \(E_0v\) denote the zero extension of \(v\in\cV_\del\) to \(\mathbb R^d\). Let \( q_\del\) denote either \( q_\del^{\mathrm{int}}\) or \( q_\del^{\mathrm{bd}}\). For \(\del\in (0,\del_0]\), there exists \(C=C(d,\Om,\ga,\del_0)>0\) such that
\begin{equation*}
\begin{aligned}
    \abs{ q_\del E_0v}_{\cS_\del(\mathbb R^d)}
    &\le C\abs{v}_{\cV_\del}
        +C\del^{-1}\norm{v}_{L^2(\Omega^-_{3\del})} \\
    &\le C\abs{v}_{\cV_\del}.
\end{aligned}
\end{equation*}
Moreover, if \(q\in H^1_0(\Om)\), then
\begin{equation*}
    \abs{ q_\del^{\mathrm{bd}}E_0q}_{\cS_\del(\mathbb R^d)}
    \le
    C\norm{\nabla q}_{L^2(\Omega^-_{3\del},\mathbb R^d)}
    +C\del^{-1}\norm{q}_{L^2(\Omega^-_{3\del})}.
\end{equation*}
\end{lemma}

\begin{proof}
We first note that for $v\in \cV_\del$,
\begin{equation}\label{eq:zero-extension-energy-identity}
    \abs{E_0v}_{\cS_\del(\mathbb R^d)}
    =\abs{v}_{\cV_\del}.
\end{equation}
Indeed, the kernel is supported in \(B_\del(\bm0)\), and all additional whole-space interactions have both values equal to zero.

We next estimate the cost of the cutoff.  Using
\[
     q_\del(\xb)E_0v(\xb)- q_\del(\yb)E_0v(\yb)
    =
     q_\del(\xb)(E_0v(\xb)-E_0v(\yb))
    +
    E_0v(\yb)( q_\del(\xb)- q_\del(\yb)),
\]
we obtain
\[ 
    \abs{ q_\del E_0v}_{\cS_\del(\mathbb R^d)}^2
    \le
    2\abs{E_0v}_{\cS_\del(\mathbb R^d)}^2 
    +
    2\iint_{\mathbb R^d\times\mathbb R^d}
    \ga_\del(\xb-\yb)|E_0v(\yb)|^2
    | q_\del(\xb)- q_\del(\yb)|^2\,d\xb\,d\yb .
\]
The first term is handled by \eqref{eq:zero-extension-energy-identity}. For the second term, since \(\supp \nabla q_\del \subset \overline{\Om_{2\del}^-\cap \Om_\del^\circ}\) and \(\supp \,\ga_\del \subset \overline{B_\del(\bm 0)}\), the domain of integration is restricted to the set \(\{(\xb,\yb)\in \Om\times \Om\backslash \left((\Om_{2\del}^\circ \times \Om_{2\del}^\circ )\cup (\Om_\del^- \times \Om_\del^-)\right): \abs{\xb-\yb}<\del\}\). Then the Lipschitz bound \(\abs{\nabla q_\del}\le C\del^{-1}\) and the finite second moment of \(\ga_\del\)  give  
\[
\begin{aligned}
    &\iint_{\mathbb R^d\times\mathbb R^d}
    \ga_\del(\xb-\yb)|E_0v(\yb)|^2
    | q_\del(\xb)- q_\del(\yb)|^2\,d\xb\,d\yb \\
    &\qquad\le
    C\del^{-2}\norm{v}_{L^2(\Omega^-_{3\del})}^2
    \int_{\mathbb R^d}\ga_\del(\zb)|\zb|^2\,d\zb
    \le
    C\del^{-2}\norm{v}_{L^2(\Omega^-_{3\del})}^2 .
\end{aligned}
\]
This proves the first inequality in \Cref{lem:thin-layer-cutoff-stability}; \Cref{prop:nonlocal-boundary-layer}, with \(a=3\), proves its second inequality.

It remains to show the second estimate of the lemma. We repeat the proof above using the two-term decomposition with a special treatment of the first term. If \(q\in H^1_0(\Om)\), then \(E_0q\in H^1(\mathbb R^d)\). Since \( q_\del^{\mathrm{bd}}(\xb) = 0\) for \(\xb\in \Om_{2\del}^\circ\), the first term can be estimated as follows using \Cref{lem:shifted-H1-bound}:
\[
\begin{split}
    \int_{\R^d}\int_{\R^d} \ga_\del(\xb-\yb)
    |E_0q(\xb)-E_0q(\yb)|^2
    \mathbf 1_{\{ q_\del^{\mathrm{bd}}(\xb)\ne0\}}\,d\yb\,d\xb
    &\le \int_{B_\del(\bm 0)}\ga_\del(\zb)\int_{(\partial\Om)^{2,1}_\del}\abs{E_0 q(\xb + \zb) - E_0 q(\xb)}^2\,d\xb\,d\zb\\
    &\le C \norm{\nabla E_0 q}_{L^2((\partial\Om)^{3,2}_\del,\mathbb R^d)}^2.
\end{split}\]
The second term is bounded by \(C\del^{-2}\norm{q}_{L^2(\Omega^-_{3\del})}^2\), as above. Taking square roots proves the second estimate.
\end{proof}

We are now ready to prove the boundary residual estimate.

\begin{lemma}\label{lem:Rbd-bound}
Let \(u_0 =U_0-G|_{\Om}\in H^{m+1}(\Om)\) satisfy \Cref{assu:compatibility}, i.e.,
\[
    \Tr(\normal^j u_0)=0,
    \qquad j=0,\ldots,m-1.
\]
Let \(\RbdI\) be given as in \eqref{eq:Rbd}.
Then
\begin{equation*}
    |\RbdI(v)|
    \le
    C\del^{m-\frac{1}{2}}
    \norm{u_0}_{H^{m+1}(\Om)}
    \abs{ q_\del^{\mathrm{bd}}E_0v}_{\cS_\del(\mathbb R^d)}
    \le
    C\del^{m-\frac{1}{2}}
    \norm{u_0}_{H^{m+1}(\Om)}
    \norm{v}_{\cV_\del},
    \qquad \forall v\in\cV_\del.
\end{equation*}
\end{lemma}

\begin{proof}
Write \(\phi_\del= q_\del^{\mathrm{bd}}E_0v\).   Decompose \(\RbdI\) into local and nonlocal parts:
\[
    \RbdI(v)=I(v)+II(v),
\]
where
\[
    I(v):=\int_\Om -\Delta u_0(\xb) q_\del^{\mathrm{bd}}(\xb)v(\xb)\,d\xb=\int_{\Omega^-_{2\del}} -\Delta u_0(\xb)\phi_\del(\xb)\,d\xb
\]
and
\[
II(v):=-a_\del(E_0 u_0,  q_\del^{\mathrm{bd}} v)=-\iint_{\Omh\times \Omh} \ga_\del(\xb-\yb)(E_0 u_0(\xb)-E_0 u_0(\yb))(\phi_\del(\xb)-\phi_\del(\yb))\,d\xb\,d\yb.
\]
We first estimate the nonlocal part \(II(v)\). Using the support of \( q_\del^{\mathrm{bd}}\) we can restrict the domain of integration to \((\partial\Om)^{3,1}_\del\times (\partial\Om)^{3,1}_\del\) in the double integral where \((\partial\Om)^{3,1}_\del\) is defined in \eqref{eq:near-boundary-layer}.
Since \((\partial\Om)^{3,1}_\del + B_\del(\bm 0)\subset (\partial\Om)^{4,2}_\del\), by the finite second moment of \(\ga_\del\) and \Cref{lem:shifted-H1-bound}, we have
\[
\begin{split}
    \int_{(\partial\Om)^{3,1}_\del}\int_{(\partial\Om)^{3,1}_\del} \ga_\del(\xb-\yb)\abs{E_0 u_0(\xb)-E_0 u_0(\yb)}^2 d\xb d\yb&\le \int_{B_\del(\bm 0)}\ga_\del(\zb)\int_{(\partial\Om)^{3,1}_\del}\abs{E_0 u_0(\xb+\zb)-E_0 u_0(\xb)}^2 d\xb d\zb\\
    &\le \left(\int_{B_\del(\bm 0)}\abs{\zb}^2\ga_\del(\zb)d\zb\right)\norm{\nabla E_0 u_0}_{L^2((\partial\Om)^{4,2}_\del,\R^d)}^2.
\end{split}
\]
Then Cauchy-Schwarz gives
\[
\begin{aligned}
\abs{II(v)}
&\le \left(\int_{(\partial\Om)^{3,1}_\del}\int_{(\partial\Om)^{3,1}_\del}
 \ga_\del(\xb-\yb)\abs{E_0 u_0(\xb)-E_0 u_0(\yb)}^2
 \,d\xb\,d\yb\right)^{1/2}
 \abs{\phi_\del}_{\cS_\del(\mathbb R^d)}\\
&\le \sqrt{d}\norm{\nabla E_0 u_0}_{L^2((\partial\Om)^{4,2}_\del,\R^d)}
 \abs{\phi_\del}_{\cS_\del(\mathbb R^d)}.
\end{aligned}
\]
Since \(E_0 u_0\) vanishes in the exterior collar, its weak gradient also
vanishes there. Hence
\[
\norm{\nabla E_0u_0}_{L^2((\partial\Om)^{4,2}_\del,\mathbb R^d)}
=
\norm{\nabla E_0u_0}_{L^2((\partial\Om)^{4,1}_\del,\mathbb R^d)}.
\]
Since \(E_0u_0\) has zero normal traces across \(\partial\Om\) up to
order \(m-1\), \Cref{cor:local-boundary-double-layer-est}, applied with
\(a=4\) and \(b=1\), gives
\[
\norm{\nabla E_0 u_0}_{L^2((\partial\Om)^{4,2}_\del,\mathbb R^d)}
\le C\del^{m-\frac{1}{2}}\norm{u_0}_{H^{m+1}(\Om)}.
\]
Thus
\[
    |II(v)|
    \le 
    C\del^{m-\frac{1}{2}}\norm{u_0}_{H^{m+1}(\Om)}
    \abs{\phi_\del}_{\cS_\del(\mathbb R^d)}.
\]
We next estimate \(I(v)\). Since \(\phi_\del\) is supported in the boundary collar, \Cref{prop:nonlocal-boundary-layer}, applied to its restriction to \(\widehat\Omega_\del\), gives
\[
    \norm{\phi_\del}_{L^2(\Omega^-_{2\del})}
    \le C\del\abs{\phi_\del}_{\cS_\del(\mathbb R^d)}.
\]
If \(m=1\), then \(\norm{\Delta u_0}_{L^2(\Omega^-_{2\del})}\le C\norm{u_0}_{H^2(\Om)}\), and therefore
\[
    |I(v)|
    \le C\del\norm{u_0}_{H^2(\Om)}
    \abs{\phi_\del}_{\cS_\del(\mathbb R^d)}
    \le C\del^{1/2}\norm{u_0}_{H^2(\Om)}
    \abs{\phi_\del}_{\cS_\del(\mathbb R^d)}.
\]
If \(m\ge2\), the zero traces also make \(\nabla^2u_0\) small in the collar: \Cref{cor:local-boundary-double-layer-est}, with \(a=2\), \(b=1\),
\(S^-=u_0\), and \(S^+=0\), gives
\(\norm{\nabla^2u_0}_{L^2(\Omega^-_{2\del},\mathbb R^{d\times d})}\le C\del^{m-\frac{3}{2}}\norm{u_0}_{H^{m+1}(\Om)}\). Combining this with the preceding estimate gives
\[
    |I(v)|
    \le
    C\del^{m-\frac{3}{2}}\norm{u_0}_{H^{m+1}(\Om)}\cdot
    \del\abs{\phi_\del}_{\cS_\del(\mathbb R^d)}
    =
    C\del^{m-\frac{1}{2}}\norm{u_0}_{H^{m+1}(\Om)}
    \abs{\phi_\del}_{\cS_\del(\mathbb R^d)}.
\]
Since \(\del_0<1\), these bounds prove the first inequality in \Cref{lem:Rbd-bound}. The second inequality in \Cref{lem:Rbd-bound} follows from \Cref{lem:thin-layer-cutoff-stability}.
\end{proof}

\begin{proof}[Energy estimate in \Cref{thm:LapI-rates-continuous}]
Combining \Cref{lem:Rbd-bound} and \eqref{eq:RbdII-bound}, we obtain 
\begin{equation*} 
    |\Rbd(v)|
    \le
    C\del^{\min\left\{2, m-\frac{1}{2}\right\}}
    \left(\norm{u_0}_{H^{m+1}(\Om)}+\norm{G}_{H^{m+1}(\widehat\Om_{\del_0})}\right)
    \norm{v}_{\cV_\del},
    \qquad \forall v\in\cV_\del.
\end{equation*}
Then together with \Cref{lem:Rint-bound},
\begin{equation}\label{eq:residual-bound}
    \norm{\Rdel}_{\cV_\del^*}
    \le
    C\del^{\min\{2,m-\frac{1}{2}\}}\left(\norm{u_0}_{H^{m+1}(\Om)}+\norm{G}_{H^{m+1}(\widehat\Om_{\del_0})}\right).
\end{equation}

Thus \Cref{thm:abstract-energy-rate} gives \eqref{eq:LapI-energy-rate}.
\end{proof}

\subsubsection{ \(L^2\) estimate}

We now prove the \(L^2\) estimate using a duality argument. This is a nonlocal analog of the gain in the Aubin--Nitsche argument in finite element analysis. 

For \(\phi\in L^2(\Om)\), let \(z_\del[\phi]\) solve the dual problem \eqref{eq:dual-problem-abstract} with \(a_\del(\cdot,\cdot)\) given by \eqref{eq:fullSpace:bilinearForm}. Since \(a_\del(\cdot,\cdot)\) is symmetric, \(z_\del[\phi]\) also solves the primal problem \eqref{eq:fullSpace:Poisson:weak} with \(f=\phi\) and \(G=0\). Let \(z\) solve the homogeneous Poisson problem
\begin{equation}\label{eq:local-dual}
    -\Delta z=\phi\quad\text{in }\Om,
    \qquad
    z=0\quad\text{on }\partial\Om.
\end{equation}
Since \(\Om\in C^{1,1}\),
\begin{equation}\label{eq:local-dual-reg}
    \norm{z}_{H^2(\Om)}\le C_{\rm reg}\norm{\phi}_{L^2(\Om)}.
\end{equation}
Applying the already proved energy estimate \eqref{eq:LapI-energy-rate} with \(m=1\) and \(G=0\) to the solutions \(z_\del[\phi]\) and \(z\) gives
\begin{equation}\label{eq:dual-energy-approx-phi}
    \norm{z_\del[\phi]-E_0z}_{\cU_\del}
    \le
    C\del^{1/2}\norm{\phi}_{L^2(\Om)}.
\end{equation}

\begin{proposition}\label{prop:dual-boundary-layer}
Extend \(z_\del[\phi]\) by zero to \(\mathbb R^d\). Then
\begin{equation*}
    \abs{ q_\del^{\mathrm{bd}}E_0z_\del[\phi]}_{\cS_\del(\mathbb R^d)}
    \le
    C\del^{1/2}\norm{\phi}_{L^2(\Om)}.
\end{equation*}
\end{proposition}

\begin{proof}
We first separate the local dual solution from the nonlocal error:
\[
    z_\del[\phi]=E_0z+(z_\del[\phi]-E_0z).
\]

Extend the two terms in the preceding decomposition by zero to \(\mathbb R^d\). By the first estimate in \Cref{lem:thin-layer-cutoff-stability} and \eqref{eq:dual-energy-approx-phi},
\[
    \abs{ q_\del^{\mathrm{bd}}E_0(z_\del[\phi]-E_0z)}_{\cS_\del(\mathbb R^d)}
    \le C\norm{z_\del[\phi]-E_0z}_{\cU_\del}
    \le C\del^{1/2}\norm{\phi}_{L^2(\Om)}.
\]
For the local term, the second estimate in \Cref{lem:thin-layer-cutoff-stability},
\Cref{cor:local-boundary-double-layer-est} applied with
\(m=1\), \(a=3\), \(b=1\), \(S^-=z\), and \(S^+=0\), and
\eqref{eq:local-dual-reg} give
\[
\begin{aligned}
    \abs{ q_\del^{\mathrm{bd}}E_0z}_{\cS_\del(\mathbb R^d)}
    &\le
    C\norm{\nabla z}_{L^2(\Omega^-_{3\del},\mathbb R^d)}
    +C\del^{-1}\norm{z}_{L^2(\Omega^-_{3\del})} \\
    &\le C\del^{1/2}\norm{z}_{H^2(\Om)}
    \le C\del^{1/2}\norm{\phi}_{L^2(\Om)}.
\end{aligned}
\]
Combining the two terms proves \Cref{prop:dual-boundary-layer}.
\end{proof}

\begin{proof}[\(L^2\) estimate in \Cref{thm:LapI-rates-continuous}]
We first apply the \(\RbdI(v)\) estimate, i.e., the first inequality in
\Cref{lem:Rbd-bound}, to \(v=z_\del[\phi]\in \cV_\del\). The additional \(\del^{\frac{1}{2}}\) factor
in \Cref{prop:dual-boundary-layer} gives
\begin{equation}\label{eq:RbdI-bound-dual}
\begin{aligned}
    |\RbdI(z_\del[\phi])|
    &\le
    C\del^{m-\frac{1}{2}}\norm{u_0}_{H^{m+1}(\Om)}
    \abs{ q_\del^{\mathrm{bd}}E_0z_\del[\phi]}_{\cS_\del(\mathbb R^d)} \\
    &\le
    C\del^m\norm{u_0}_{H^{m+1}(\Om)}
    \norm{\phi}_{L^2(\Om)}.
\end{aligned}
\end{equation}
This is precisely the gain \(\del^{m-\frac{1}{2}}\del^{\frac{1}{2}}=\del^m\).

We now use the stability of the dual problem,
\(\norm{z_\del[\phi]}_{\cV_\del}\le C\norm{\phi}_{L^2(\Om)}\), in \Cref{lem:Rint-bound} and \eqref{eq:RbdII-bound}. Together with \eqref{eq:RbdI-bound-dual}, this yields
\[
    |\Rdel(z_\del[\phi])|
    \le
    C\del^{\min\{2,m\}}\left(\norm{u_0}_{H^{m+1}(\Om)}+\norm{G}_{H^{m+1}(\widehat\Om_{\del_0})}\right)\norm{\phi}_{L^2(\Om)}.
\]
Applying \Cref{thm:abstract-L2-rate} gives the \(L^2\) estimate in \Cref{thm:LapI-rates-continuous}.
\end{proof} 

\subsubsection{Discrete estimates}
We now turn to the discrete problem and prove corresponding energy and \(L^2\) estimates. The energy estimate is an application of \Cref{thm:abstract-energy-discrete}, while the \(L^2\) estimate is an application of \Cref{thm:abstract-L2-rate-discrete}. The main task is to verify the assumptions of these abstract results.
\begin{proof}[Proof of \Cref{thm:LapI-rates-discrete}]
We first prove the energy estimate \eqref{eq:LapI-energy-rate-discrete}. We identify finite element functions with their zero extensions when they are measured in \(\cU_\del\). Since \(u_0\in H^1_0(\Om)\), the function \(I_hu_0\) belongs to
\(V_{h,0}^{k}\subset\cU_{\del,h}\).  The embedding
\eqref{eq:H1embedS_fullSpace} and the approximation property
\eqref{eq:FE-H1-approximation} therefore give
\begin{equation}\label{eq:LapI-energy-approximation}
    \inf_{u_h\in \cU_{\del,h}}\norm{E_0 u_0-u_h}_{\cU_\del}
    \le
    \norm{E_0 u_0-I_hu_0}_{\cU_\del}
    \le
    C\norm{u_0-I_hu_0}_{H^1(\Om)}
    \le
    C h^{\min\{m,k\}}\norm{u_0}_{H^{m+1}(\Om)}.
\end{equation}

Together with the residual estimate \eqref{eq:residual-bound} and \Cref{thm:abstract-energy-discrete}, this proves \eqref{eq:LapI-energy-rate-discrete}.

We next verify the two assumptions of
\Cref{thm:abstract-L2-rate-discrete}. Let \(\phi\in L^2(\Om)\), let
\(z_{\del,h}[\phi]\)  solve the discrete dual problem \eqref{eq:dual-problem-abstract-discrete} with \(a_\del(\cdot,\cdot)\) given by \eqref{eq:fullSpace:bilinearForm}, and let
\(z_\del[\phi]\) and \(z\) be given as before, i.e., they solve \eqref{eq:dual-problem-abstract} with the same \(a_\del(\cdot,\cdot)\) and
\eqref{eq:local-dual}, respectively. We first approximate the continuous dual
solution. Since \(I_hz\in V_{h,0}^{k}\subset\cV_{\del,h}\), by
\eqref{eq:dual-energy-approx-phi},
\eqref{eq:H1embedS_fullSpace}, \eqref{eq:FE-H1-approximation}, and local
\(H^2\)-regularity,
\begin{equation}\label{eq:dual-approximation}
\begin{aligned}
    \inf_{v_h\in \cV_{\del,h}}
    \norm{z_\del[\phi]-v_h}_{\cV_\del}
    &\le
    \norm{z_\del[\phi]-I_hz}_{\cV_\del} \\
    &\le
    \norm{z_\del[\phi]-E_0z}_{\cV_\del}
    +C\norm{z-I_hz}_{H^1(\Om)} \\
    &\le
    C(\del^{\frac{1}{2}}+h)\norm{\phi}_{L^2(\Om)}.
\end{aligned}
\end{equation}

Thus the dual approximation rate is
\(\rho(\del,h)=C(\del^{1/2}+h)\) in \Cref{assu:abstract-dual-discrete}.

In view of \eqref{eq:dual-approximation}, Galerkin orthogonality and Céa's estimate applied to the dual problems give
\begin{equation}\label{eq:discrete-dual-approximation}
    \norm{z_\del[\phi]-z_{\del,h}[\phi]}_{\cV_\del}
    \le C\inf_{v_h\in \cV_{\del,h}}
    \norm{z_\del[\phi]-v_h}_{\cV_\del}\le
    C(\del^{\frac{1}{2}}+h)\norm{\phi}_{L^2(\Om)}.
\end{equation}
Coercivity of the discrete dual problem also gives
\begin{equation}\label{eq:discrete-dual-stability}
    \norm{z_{\del,h}[\phi]}_{\cV_\del}
    \le C\norm{\phi}_{L^2(\Om)}.
\end{equation}

It remains to estimate the complete residual against the discrete dual
solution. The estimate in \Cref{lem:Rint-bound} and
\eqref{eq:discrete-dual-stability} give
\[
    |\Rint(z_{\del,h}[\phi])|
    \le
    C\del^{\min\{2,m\}}\norm{U_0}_{H^{m+1}(\Om)}
    \norm{\phi}_{L^2(\Om)}.
\]
Similarly, \eqref{eq:RbdII-bound} and \eqref{eq:discrete-dual-stability} give
\[
    |\RbdII(z_{\del,h}[\phi])|
    \le
    C\del^{\min\{2,m\}}\norm{G}_{H^{m+1}(\widehat\Om_{\del_0})}
    \norm{\phi}_{L^2(\Om)}.
\]
For the boundary part, we compare the discrete and continuous dual solutions.
By \eqref{eq:RbdI-bound-dual}, the second inequality in \Cref{lem:Rbd-bound}, and
\eqref{eq:discrete-dual-approximation},
\[
\begin{aligned}
    |\RbdI(z_{\del,h}[\phi])|
    &\le
    |\RbdI(z_\del[\phi])|
    +|\RbdI(z_\del[\phi]-z_{\del,h}[\phi])| \\
    &\le
    C\left(\del^m+\del^{m-\frac{1}{2}}(\del^{\frac{1}{2}}+h)\right)
    \norm{u_0}_{H^{m+1}(\Om)}\norm{\phi}_{L^2(\Om)} \\
    &\le
    C\left(\del^{\min\{2,m\}}+\del^{m-\frac{1}{2}}h\right)
    \norm{u_0}_{H^{m+1}(\Om)}\norm{\phi}_{L^2(\Om)}.
\end{aligned}
\]
Here we used \(\min\{2,m\}\le m\) and \(\del_0\le1\). Combining the
interior and boundary estimates gives
\begin{equation}\label{eq:discrete-dual-consistency-homo}
    \abs{\Rdel(z_{\del,h}[\phi])}
    \le
    C\left(\del^{\min\{2,m\}}+\del^{m-\frac{1}{2}}h\right)
    \left(\norm{u_0}_{H^{m+1}(\Om)}+\norm{G}_{H^{m+1}(\widehat\Om_{\del_0})}\right)
    \norm{\phi}_{L^2(\Om)}.
\end{equation}
This verifies the first bound in \Cref{assu:abstract-dual-discrete} with
\[
    \beta(\del,h)
    =
    C\left(\del^{\min\{2,m\}}+\del^{m-\frac{1}{2}}h\right)
    \left(\norm{u_0}_{H^{m+1}(\Om)}+\norm{G}_{H^{m+1}(\widehat\Om_{\del_0})}\right).
\]
Since \(M_a=1\) and \(\tilde{\al} = (1+C_P^2)^{-1}\) in the present application,
\Cref{thm:abstract-L2-rate-discrete}, together with the approximation
estimate \eqref{eq:LapI-energy-approximation} and \eqref{eq:dual-approximation}, gives
\[
\begin{aligned}
    \norm{u_{\del,h}-u_0}_{L^2(\Om)}
    \le C\Bigl(&
        \del^{\min\{2,m\}}+\del^{m-\frac{1}{2}}h\\
        &+(\del^{\frac{1}{2}}+h)h^{\min\{m,k\}}
    \Bigr)\left(\norm{u_0}_{H^{m+1}(\Om)}+\norm{G}_{H^{m+1}(\widehat\Om_{\del_0})}\right).
\end{aligned}
\]
Note that the second term in \(\del\) and \(h\) can be absorbed into the first and the last. Indeed, Young's
inequality and \(0<\del\le1\) give
\begin{equation*}
\begin{aligned}
    \del^{m-\frac{1}{2}}h
    &\le
    C\del^{(m-\frac{1}{2})
    \frac{\min\{m,k\}+1}{\min\{m,k\}}}
    +C h^{\min\{m,k\}+1}\\
    &\le
    C\left(
        \del^{\min\{2,m\}}+h^{\min\{m,k\}+1}
    \right).
\end{aligned}
\end{equation*}
This proves \eqref{eq:LapI-L2-rate-discrete}. 
\end{proof}

\begin{remark}\label{rem:improved-rates-integrable-kernel}
If \(\ga\in L^1(\R^d)\), the convergence rates in    \eqref{eq:LapI-energy-rate-discrete} and \eqref{eq:LapI-L2-rate-discrete} can be improved by
\[
 \del^{\min\{2,m-1/2\}}+h^{\min\{m,k\}} \min\{1,\del^{-1}h\}
\]
and
\[
 \del^{\min\{2,m\}}+
 \bigl(\del^{1/2}+ h \min\{1,\del^{-1}h\}\bigr) h^{\min\{m,k\}} 
 \min\{1,\del^{-1}h\},
\]
respectively. This gives better estimates when $h < \delta$. 
Indeed, for \(v\in L^2(\Om)\), zero extension and
\eqref{eq:gamma-scaling} give
\[
\begin{aligned}
 \norm{E_0v}_{\cU_\del}^2
 &=\norm{v}_{L^2(\Om)}^2+
   \int_{\R^d}\ga_\del(\zb)
   \norm{E_0v(\,\cdot+\zb)-E_0v}_{L^2(\R^d)}^2\,d\zb\\
 &\le\bigl(1+4\del^{-2}\norm{\ga}_{L^1(\R^d)}\bigr)
    \norm{v}_{L^2(\Om)}^2.
\end{aligned}
\]
Therefore, the term in \eqref{eq:LapI-energy-approximation} can also be bounded by 
\[
\inf_{u_h\in \cU_{\del,h}}\norm{E_0 u_0-u_h}_{\cU_\del}
    \le 
    \norm{E_0 u_0-I_hu_0}_{\cU_\del}  \le C \delta^{-1}\norm{u_0 -I_h u_0 }_{L^2(\Omega)} \leq C \delta^{-1} h^{\min\{m,k\}+1} \| u_0\|_{H^{m+1}(\Omega)} 
\]
using the \(L^2\) approximation estimate
\[
 \norm{v-I_hv}_{L^2(\Om)}
 \le Ch^{s+1}\norm{v}_{H^{s+1}(\Om)},
 \qquad s=1,\ldots,\min\{m,k\}
\]
for \(v\in H^{s+1}(\Om)\cap H^1_0(\Om)\) in place of \eqref{eq:FE-H1-approximation}. This gives the improved rate in the energy norm. The constant \(C\) may additionally depend on \(\norm{\ga}_{L^1}\), but is
independent of \(\del\) and \(h\). 
 Applying the same argument to the
local \(H^2\) dual solution replaces the dual approximation factor in \eqref{eq:dual-approximation} by
\(\del^{1/2}+\min\{h,\delta^{-1}h^2\}\). Therefore, the same argument as before and \Cref{thm:abstract-L2-rate-discrete} give the improved $L^2$ estimate. 
\end{remark}

\begin{proof}[Proof of \Cref{prop:LapI-piecewise-constant}]
Define 
\(Q_h:L^2(\Om)\to V_h^0\)  by
\begin{equation}\label{eq:Qh-definition}
    (Q_hu)|_K=\frac{1}{|K|}\int_K u(\xb)\,d\xb.
\end{equation}
We first establish the approximation estimate that replaces
\eqref{eq:FE-H1-approximation}, namely,
\begin{equation}\label{eq:piecewise-constant-approximation}
     \norm{u-Q_hu}_{\cU_\del}
 \le
 C\left(\frac{h}{\del}\right)^{1-\alpha}
 \norm{u}_{H^1(\Om)},\qquad \forall u\in H^1_0(\Om).
\end{equation}
First,
\cite[Lemma~3.1]{tian2014asymptotically}, \eqref{eq:gamma-scaling} and the moment assumption in \Cref{prop:LapI-piecewise-constant} give
\begin{equation*}
 \norm{u}_{\cU_\del}
 \le C\left(\int_{B_1(\bm 0)} \abs{\zb}^{2\alpha}\ga_\del(\zb)d\zb\right)^{\frac{1}{2}}\norm{u}_{H^{\alpha}(\Om)}\le
 C\del^{\alpha-1}
 \norm{u}_{H^{\alpha}(\Om)},
 \qquad \forall u\in H^{\alpha}_0(\Om).
\end{equation*}
This estimate also proves that every function in \(V_h^0\)
belongs to \(\cS_\del(\Om)\), since piecewise constants belong to
\(H^\alpha(\Om)\) for every \(0\le \alpha<1/2\).

By \cite[(1.3)]{belgacem2001some}, with suitable adaptation to the present curved meshes, we have
\begin{equation}\label{eq:fractional-approximation}
     \norm{u-Q_hu}_{H^{\alpha}(\Om)}
 \le
 Ch^{1-\alpha}\norm{u}_{H^1(\Om)},
 \qquad \forall u\in H^1(\Om).
\end{equation}
Applying the energy bound above to
\(u-Q_hu\) yields \eqref{eq:piecewise-constant-approximation}.

We now apply the abstract estimates.  Taking \(u=u_0\) in
\eqref{eq:piecewise-constant-approximation}, and combining the
result with \eqref{eq:residual-bound} and
\Cref{thm:abstract-energy-discrete}, proves
the energy estimate in \Cref{prop:LapI-piecewise-constant}.

For the \(L^2\) estimate, let \(z_\del[\phi]\) and \(z\) be the nonlocal and
local dual solutions used in the proof of
\Cref{thm:LapI-rates-discrete}.  The dual energy estimate
\eqref{eq:dual-energy-approx-phi},
\eqref{eq:piecewise-constant-approximation}, and
\eqref{eq:local-dual-reg} give
\[
 \inf_{v_h\in V_h^0}
 \norm{z_\del[\phi]-v_h}_{\cV_\del}
 \le
 \norm{z_\del[\phi]-E_0Q_hz}_{\cV_\del}
 \le
 C\left(
  \del^{1/2}
  +\left(\frac{h}{\del}\right)^{1-\alpha}
 \right)\norm{\phi}_{L^2(\Om)}.
\]
Galerkin orthogonality and coercivity now give
\[
\begin{split}
 \norm{z_\del[\phi]-z_{\del,h}[\phi]}_{\cV_\del}
 &\le C\left(
  \del^{1/2}
  +\left(\frac{h}{\del}\right)^{1-\alpha}
 \right)\norm{\phi}_{L^2(\Om)},\\
 \norm{z_{\del,h}[\phi]}_{\cV_\del}
 &\le C\norm{\phi}_{L^2(\Om)}.
\end{split}
\]

The residual calculation in
\eqref{eq:discrete-dual-consistency-homo} uses the finite element space only
through these two dual estimates.  Replacing
\eqref{eq:discrete-dual-approximation} by the preceding bound gives
\[
\begin{split}
 \abs{\Rdel(z_{\del,h}[\phi])}
 \le C\left(
  \del^{\min\{2,m\}}
  +\del^{m-\frac{1}{2}}
   \left(\frac{h}{\del}\right)^{1-\alpha}
 \right)
 \left(
  \norm{u_0}_{H^{m+1}(\Om)}
  +\norm{G}_{H^{m+1}(\widehat\Om_{\del_0})}
 \right)
 \norm{\phi}_{L^2(\Om)}.
\end{split}
\]
Applying \Cref{thm:abstract-L2-rate-discrete}, using
\eqref{eq:piecewise-constant-approximation} for the primal
approximation, and observing that
\(\del^{m-\frac{1}{2}}\le\del^{1/2}\) for \(m\ge1\), proves
the \(L^2\) estimate in \Cref{prop:LapI-piecewise-constant}.  
\end{proof}

\section{Applications to nonlocal Laplacians with constant horizon of type (II)}\label{sec:NonlocalLaplaceII}

This section develops the convergence-rate theory for the second type of nonlocal Laplacian operators, constructed using half-space
first-order operators. More specifically, we define half-space gradient and divergence operators and consider their composition as a nonlocal Laplacian. These types of operators have been studied in \cite{han2023nonlocal,han2024compactnessa} and originate from the peridynamic correspondence models \cite{silling2007peridynamic,silling2017stability}.
We prove the same convergence rates as in \Cref{sec:NonlocalLaplaceI} for the first type of nonlocal Laplacian operators.

Technical facts about half-space operators are presented in
\Cref{sec:halfspace-analytic-foundations}.

To begin with, we recall that our goal is to prove that the solutions to nonlocal Poisson problems converge to that of the local problem \eqref{eq:local_Poisson_strong}. In \Cref{sec:NonlocalLaplaceI}, the nonlocal operator \(\cL_\del\) is defined by an integral involving a kernel \(\ga_\del\) that is radially symmetric and has finite second moment. In this section, we consider a different type of nonlocal Laplacian, which is defined as the composition of a half-space gradient and a half-space divergence. The half-space gradient operator is defined using a kernel \(w_\del\) that is radial and has finite first moment. Another key difference is that we impose nonlocal volumetric boundary conditions on a collar of width \(2\del\), see \Cref{fig:cutoff-typeII}(a), instead of the collar of width \(\del\) used in \Cref{sec:NonlocalLaplaceI}. This is because both gradient and divergence operators need an interaction region of width \(\del\) to be well-defined, and the composition of these two operators requires a collar of width \(2\del\).

 Let \(\Om\subset\R^d\) be a bounded \(C^{1,1}\) domain and fix
\(\nub\in\mathbb S^{d-1}\). In this section, \(\del_0\) is chosen by \eqref{eq:halfspace-admissible-horizon}.
We take \(\del_{\max}:=\del_0/2\) and use
\(0<\del\le\del_{\max}\) throughout this section. Define
\[
    \chi_{\nub}(\zb):=\mathbf 1_{\{\zb\cdot\nub\ge0\}}.
\]

We assume that the kernel \(w\) satisfies the following properties.
\begin{assumption}
\label{assu:halfspace-kernel}
The kernel function \(w:\R^d\to[0,\infty]\) is measurable, radial, supported in
\(B_1(\bm0)\), and
\begin{equation*}
    \int_{\R^d}|\zb|w(\zb)\,d\zb=2d.
\end{equation*}
If \(d=1\), we additionally require that $w$ is nonincreasing on $(0,\infty)$.
For \(0<\del\le\del_0/2\), define
\begin{equation*}
    w_\del(\zb):=\del^{-d-1}w(\zb/\del).
\end{equation*}
\end{assumption}

By \Cref{assu:halfspace-kernel}, the kernel \(w_\del\) is supported in \(B_\del(\bm0)\) and satisfies
\begin{equation*}
    \int_{\R^d}|\zb|w_\del(\zb)\,d\zb=2d.
\end{equation*}

Moreover,
\begin{equation*}
    \int_{\R^d}
    \chi_{\nub}(\zb)w(\zb)
    \frac{\zb\otimes\zb}{|\zb|}\,d\zb
    =I_d.
\end{equation*}
This identity is the key for the consistency of nonlocal operators.

For smooth functions on the relevant domains, the pointwise half-space gradient \(\cG^{\nub}_{w_\del}u\) is defined on \(\Omh\) from the values of \(u\) on \(\Omhh\), while \(\cD^{-\nub}_{w_\del}\vb\) is defined on \(\Om\) from the values of \(\vb\) on \(\Omh\). Thus their domains decrease by one interaction layer at each application. We use throughout the distributional gradient defined in \Cref{def:distributional-gradient,sec:halfspace-analytic-foundations}, denoted by \(\fG^{\nub}_{w_\del,\Omhh}\), and the distributional divergence \(\fD^{-\nub}_{w_\del,\Omh}:L^2(\Omh;\mathbb R^d)\to\mathcal D'(\Om)\). For the remainder of this section, we abbreviate \(\fG^{\nub}_{w_\del,\Omhh}\) by \(\fG^{\nub}_{w_\del}\), and \(\fD^{-\nub}_{w_\del,\Omh}\) by \(\fD^{-\nub}_{w_\del}\). 
The nonlocal \(H^1\) space is
\begin{equation*}
    \cS_{w_\del}^{\nub}(\Omhh):=\{u\in L^2(\Omhh):\fG^{\nub}_{w_\del} u\in L^2(\Omh,\R^d)\}, \quad \norm{u}_{\cS_{w_\del}^{\nub}(\Omhh)}^2
:=\norm{u}_{L^2(\Omhh)}^2 + \abs{u}_{\cS_{w_\del}^{\nub}(\Omhh)}^2.
\end{equation*}
Here the seminorm is
defined as
\begin{equation*}
    \abs{u}_{\cS_{w_\del}^{\nub}(\Omhh)} 
    :=\norm{\fG^{\nub}_{w_\del} u}_{L^2(\Omh,\R^d)}.
\end{equation*}
Then the energy space is   defined as
\begin{equation*}
\cS_{w_\del,0}^{\nub}(\Om)
:=\{u\in\cS_{w_\del}^{\nub}(\Omhh):
        u=0\ \text{a.e. in }\Omega_{2\del}^+\}.
\end{equation*}
For the type~(II) application, the  trial and test spaces are given
by
\begin{equation*}
\cU_\del=\cV_\del:=\cS_{w_\del,0}^{\nub}(\Om),
\end{equation*}
equipped with the full graph norm
\[
\norm{v}_{\cU_\del}=\norm{v}_{\cV_\del}
    :=\norm{v}_{\cS_{w_\del}^{\nub}(\Omhh)},
    \qquad
    \abs{v}_{\cU_\del}=\abs{v}_{\cV_\del}
    :=\abs{v}_{\cS_{w_\del}^{\nub}(\Omhh)}.
\]
  By \cite[Theorem~3.14]{han2024compactnessa}
, there exist
\(\del_P=\del_P(\Om,d,w)>0\) and \(C_P>0\) such that
\begin{equation}\label{eq:halfspace-Poincare}
    \norm{u}_{L^2(\Om)}\le C_P\abs{u}_{\cU_\del},
    \qquad u\in\cU_\del,\quad 0<\del<\del_P.
\end{equation}
  The constant  \(C_P\) is independent of \(\del\).   Here \(\del_P\) and
\(C_P\) refer to the present operator. The choice
\eqref{eq:halfspace-admissible-horizon} gives
\(\del_0/2<\del_P\), so \eqref{eq:halfspace-Poincare} holds
throughout the admissible horizon range. Enlarging \(C_P\) if necessary,
we also use it for the local Poincar\'e inequality at \(\del=0\).

We also need the nonlocal space on the whole domain \(\R^d\), which plays a role similar to \(H^1(\R^d)\), defined as
\begin{equation*}
    \cS_{w_\del}^{\nub}(\R^d):=\{u\in L^2(\R^d):\fG^{\nub}_{w_\del,\R^d} u\in L^2(\R^d,\R^d)\}, \quad
    \norm{u}_{\cS_{w_\del}^{\nub}(\R^d)}^2
    :=\norm{u}_{L^2(\R^d)}^2 + \abs{u}_{\cS_{w_\del}^{\nub}(\R^d)}^2.
\end{equation*}
Here the semi-norm is defined as
\begin{equation*}
    \abs{u}_{\cS_{w_\del}^{\nub}(\R^d)} 
    :=\norm{\fG^{\nub}_{w_\del,\R^d} u}_{L^2(\R^d,\R^d)}.
\end{equation*}

We are ready to state the nonlocal Poisson problem. The nonlocal Laplacian is
the composition
\begin{equation*}
    \cL^{\nub}_{w_\del}u:=\fD^{-\nub}_{w_\del}\fG^{\nub}_{w_\del}u.
\end{equation*}
Fix \(G\in H^1(\widehat\Om_{\del_0})\) as a lifting of the boundary data \(g\) such that \(\Tr G=g\) on \(\partial\Om\). Given \(f\in L^2(\Om)\) and \(0<\del\le\del_0/2\), the strong form of the nonlocal Poisson problem is
\begin{equation*}
    \left\{
        \begin{aligned}
            -\cL^{\nub}_{w_\del}U_\del&=f\quad\text{in }\mathcal D'(\Om),\\
            U_\del&=G\quad\text{in }\Om_{2\del}^+.
        \end{aligned}
    \right.
\end{equation*}
We turn to the homogeneous boundary value problem with solution \(u_\del= U_\del - G\). The weak form of the nonlocal Poisson problem is to find \(u_\del\in \cU_\del\) such that  
\begin{equation}\label{eq:halfspace-Poisson-weak}
    a_\del(u_\del,v)
    =\int_\Om f(\xb)v(\xb)\,d\xb - a_\del(G,v),
    \quad\forall v\in \cV_\del,
\end{equation}
where the bilinear form is defined as
\begin{equation*}
    a_\del(u,v)
    :=\int_{\widehat\Omega_\del}
       \fG_{w_\del}^{\nub}u(\xb)\cdot
       \fG_{w_\del}^{\nub}v(\xb)\,d\xb.
\end{equation*}
With \(U_\del:=u_\del+G\), the definition of the distributional
divergence and the density of \(C_c^\infty(\Om)\) in \(\cV_\del\)
show that this weak formulation is equivalent to the strong
problem for \(U_\del\).
By Cauchy--Schwarz and the uniform Poincaré inequality
\eqref{eq:halfspace-Poincare},
\[
    \abs{a_\del(u,v)}
    \le \norm{u}_{\cU_\del}\norm{v}_{\cV_\del},
    \qquad
    a_\del(u,u)
    \ge(1+C_P^2)^{-1}
       \norm{u}_{\cU_\del}^2.
\]
Thus the continuity constant is \(M_a=1\), while the coercivity constant is
\((1+C_P^2)^{-1}\). The Lax--Milgram theorem gives a unique solution
\(u_\del\in \cU_\del\) to
\eqref{eq:halfspace-Poisson-weak}. Moreover,
\begin{equation*}
    \norm{u_\del}_{\cU_\del}
    \le C\left(
       \norm{f}_{L^2(\Om)}
       +\norm{\fG_{w_\del}^{\nub}G}
          _{L^2(\widehat\Omega_\del,\R^d)}
    \right),
\end{equation*}
where \(C\) is independent of \(\del\).

We record the uniform embedding \(H^1_0(\Om)\hookrightarrow \cU_\del\): there exists a constant \(C>0\) such that
\begin{equation}\label{eq:halfspace-uniform-H1-bound}
    \norm{E_0 v}_{\cU_\del}\le C\norm{v}_{H^1(\Om)},\qquad \forall v\in H^1_0(\Om).
\end{equation}

To apply the abstract framework in \Cref{sec:framework}, take
\[
    \cU_*=\cV_*:=L^2(\Om),\qquad
    \cU_0=\cV_0:=H^1_0(\Om).
\]
Equip these spaces with their full norms, retain the previously defined \(\cU_\del\) and \(\cV_\del\) for
  \(0<\del\le\del_0/2\), and set
\(a_0(u,v):=\int_\Om\nabla u(\xb)\cdot\nabla v(\xb)\,d\xb\).
Then \eqref{eq:halfspace-uniform-H1-bound} and the preceding estimates give
\[
    M_*=1,\qquad M_0=C,\qquad M_a=1,\qquad
    \alpha=(1+C_P^2)^{-1}.
\]
Let \(U_0\) solve \eqref{eq:local_Poisson_strong}, set
\(u_0:=U_0-G\), and define
\begin{equation*}
    \Uext:=E_0u_0+G.
\end{equation*}
The residual functional is given by
\begin{equation*}
    \Rdel(v)
    :=\int_\Om f(\xb)v(\xb)\,d\xb
      -a_\del(\Uext,v),
    \qquad \forall v\in\cV_\del.
\end{equation*}
Assume that \(U_0\in H^2(\Om)\). Then \(-\Delta U_0 = f\) by \eqref{eq:local_Poisson_strong}. Thus, 
\begin{equation}\label{eq:halfspace-Rdel}
    \Rdel(v) = \int_\Om -\Delta U_0(\xb)v(\xb)d\xb - \int_{\widehat\Omega_\del} \fG_{w_\del}^{\nub}\Uext(\xb)\cdot \fG_{w_\del}^{\nub}v(\xb)d\xb.
\end{equation}
We will use this definition of \(\Rdel\) throughout the rest of this section. The reader should distinguish this from the residual functional \(\Rdel\) defined in \Cref{sec:NonlocalLaplaceI}.

\subsection{Continuous convergence rates}
\label{subsec:halfspace-main-results}
Throughout this subsection, we assume that \Cref{assu:halfspace-kernel,assu:compatibility} hold. The following theorem states the convergence rates for the continuous half-space problem.
As in the type~(I) application, the solution  norms on the right-hand side can be
bounded in terms of sufficiently regular problem data under the hypotheses
stated in \Cref{rem:LapI-sufficient-data-regularity}.

The fixed collar width is chosen independently of its value in
\Cref{sec:NonlocalLaplaceI}:
\begin{equation}\label{eq:halfspace-admissible-horizon}
    \del_0:=\frac{1}{2}\min\{\rho_0,\rho_1,\del_P\}.
\end{equation}
Here \(\rho_0\) and \(\rho_1\) are given in
\Cref{lem:boundary-layer-est} and \eqref{eq:def_rho_1}, while
\(\del_P\) is given in \eqref{eq:halfspace-Poincare} depending on \(w\).
For \(0<\del\le\del_0/2\), this gives
\[
    2\del\le\del_0<\min\{\rho_0,\rho_1\},
    \qquad \del<\del_P,
    \qquad \widehat\Omega_{2\del}\subset\widehat\Omega_{\del_0}.
\]
The compatibility and Sobolev extension assumptions are understood on the
fixed domain \(\widehat\Omega_{\del_0}\), which is Lipschitz since
\(\del_0<\rho_1\); see \cite[Corollary~2]{lewicka2020which}.

\begin{theorem}
\label{thm:LapII-rates-continuous}
Let \(u_\del\) solve \eqref{eq:halfspace-Poisson-weak}, and let \(E_0 u_0\) be the extension of \(u_0\). Suppose that
\Cref{assu:compatibility} holds. Then there exists a constant \(C=C(m, d,\Om, w, \del_0)>0\), such that for any \(0<\del\le\del_0/2\),
\begin{equation}\label{eq:LapII-energy-rate}
    \norm{u_\del - E_0 u_0}_{\cU_\del} \le C \del^{\min\{2,m-\frac{1}{2}\}} \left(\norm{u_0}_{H^{m+1}(\Om)} + \norm{G}_{H^{m+1}(\widehat\Om_{\del_0})}\right),
\end{equation}
and
\begin{equation*}
    \norm{u_\del - u_0}_{L^2(\Om)} \le C \del^{\min\{2,m\}} \left(\norm{u_0}_{H^{m+1}(\Om)} + \norm{G}_{H^{m+1}(\widehat\Om_{\del_0})}\right).
\end{equation*}
\end{theorem}
\begin{remark}\label{rem:LapII-energy-rate-high-regularity}
    Similar to \Cref{rem:LapI-energy-rate-high-regularity}, if \Cref{assu:compatibility} holds with \(j=0,\cdots, m\), then \(\Uext\in H^{m+1}(\widehat\Om_{\del_0})\) and  convergence rates in \eqref{eq:LapII-energy-rate} can be improved to \(\del^{\min\{2, m\}}\), i.e.,
    \[
    \norm{u_\del-E_0u_0}_{\cU_\del}
    \le
    C\del^{\min\{2,m\}}
    \left(
        \norm{u_0}_{H^{m+1}(\Om)}
        +\norm{G}_{H^{m+1}(\widehat\Om_{\del_0})}
    \right).
    \]
\end{remark}

\subsection{Discrete convergence rates}
\label{subsec:halfspace-discrete-results}
We use the meshes \(\{\cT_h\}_{0<h\le h_0}\) and the space \(V_{h,0}^{k}\) from
\Cref{subsec:LapI-discrete-results} and recall \eqref{eq:FE-H1-approximation}. 
For each admissible \(\del\), let
\(\cU_{\del,h}=\cV_{\del,h}\) be any finite-dimensional
conforming space such that
\[
 V_{h,0}^{k}
 \subset\cU_{\del,h}
 \subset\cU_\del,
 \qquad k\ge1,
\]
where functions on \(\Om\) are identified with their zero extensions to
\(\Omhh\). 
  The discrete half-space problem is to find
\(u_{\del,h}\in\cU_{\del,h}\) such that
\begin{equation*}
    a_\del(u_{\del,h},v_h)
    =\int_\Om f(\xb)v_h(\xb)\,d\xb
      -a_\del(G,v_h),
    \qquad \forall v_h\in \cV_{\del,h},
\end{equation*}
with \(U_{\del,h}:=u_{\del,h}+G\).

\begin{theorem}
\label{thm:LapII-rates-discrete}
Under the hypotheses of
\Cref{thm:LapII-rates-continuous}, the discrete problem has a unique solution and there exists a constant \(C=C(m,k,d,\Om, w, \del_0)>0\), such that for any \(0<\del\le\del_0/2\) and \(0<h\le h_0\),
\begin{equation}\label{eq:LapII-energy-rate-discrete}
    \norm{u_{\del,h}-E_0 u_0}_{\cU_\del}
    \le C\left(
       \del^{\min\{2,m-\frac{1}{2}\}}
       +h^{\min\{m,k\}}
    \right)
    \left(\norm{u_0}_{H^{m+1}(\Om)} + \norm{G}_{H^{m+1}(\widehat\Om_{\del_0})}\right),
\end{equation}
and
\begin{equation}\label{eq:LapII-L2-rate-discrete}
    \norm{u_{\del,h}-u_0}_{L^2(\Om)}
    \le C\left(
       \del^{\min\{2,m\}}
       +(\del^{\frac{1}{2}}+h)h^{\min\{m,k\}}
    \right)
    \left(\norm{u_0}_{H^{m+1}(\Om)} + \norm{G}_{H^{m+1}(\widehat\Om_{\del_0})}\right).
\end{equation} 
\end{theorem}

\begin{proposition}
\label{prop:LapII-piecewise-constant}
Let \(V_h^0\) be the piecewise constant space in
\Cref{prop:LapI-piecewise-constant}, now identified with its zero extension
to \(\Omhh\). Suppose that, for some \(0\le\alpha<\frac{1}{2}\),
\begin{equation*}
 \int_{B_1(\bm0)}\abs{\zb}^{\alpha}w(\zb)\,d\zb<\infty.
\end{equation*}
Then \(V_h^0\subset\cU_\del\). Let
\(u_{\del,h}\in V_h^0\) solve
the discrete variational problem with
\(\cU_{\del,h}=\cV_{\del,h}=V_h^0\).  Suppose that
\Cref{assu:compatibility} holds.  Then there exists \(C=C(m,k,d,\Om,w,\al,\del_0)>0\) independent of \(\del\) and \(h\),   such that
\begin{equation*}
\begin{split}
 \norm{u_{\del,h}-E_0u_0}_{\cU_\del}
 \le C\left(
  \del^{\min\{2,m-\frac{1}{2}\}}
  +\left(\frac{h}{\del}\right)^{1-\alpha}
 \right)
 \left(
  \norm{u_0}_{H^{m+1}(\Om)}
  +\norm{G}_{H^{m+1}(\widehat\Om_{\del_0})}
 \right).
\end{split}
\end{equation*}
Moreover,
\begin{equation*}
 \norm{u_{\del,h}-u_0}_{L^2(\Om)}
 \le C\Bigg(
  \del^{\min\{2,m\}}
  +\del^{\frac{1}{2}}
   \left(\frac{h}{\del}\right)^{1-\alpha}
  +\left(\frac{h}{\del}\right)^{2-2\alpha}
 \Bigg)\left(
  \norm{u_0}_{H^{m+1}(\Om)}
  +\norm{G}_{H^{m+1}(\widehat\Om_{\del_0})}
 \right).
\end{equation*}

Consequently, the estimates converge as \(\del\to0\) whenever
\(h/\del\to0\). 
\end{proposition}

\subsection{Proofs of the main results}
\label{subsec:halfspace-main-results-proofs}
To prove the main results, we follow the same strategy as in \Cref{sec:NonlocalLaplaceI}. First we prove the energy estimate \(O(\del^{\min\{2,m-\frac{1}{2}\}})\) using the consistency error described by the residual functional, where the interior part \(O(\del^{\min\{2,m\}})\) is estimated by Fourier analysis and the boundary part \(O(\del^{\min\{2,m-\frac{1}{2}\}})\) is obtained by the boundary-layer estimate. Then, for the \(L^2\) estimate, we show a dual boundary layer estimate to gain an additional \(O(\del^{\frac{1}{2}})\) for the boundary part so that in combination the \(L^2\) error is \(O(\del^{\min\{2,m\}})\). For the discrete problems, we estimate the additional approximation error and apply the abstract theorems \Cref{thm:abstract-energy-discrete,thm:abstract-L2-rate-discrete}. Note that we need auxiliary boundary-layer estimates analogous to \Cref{prop:nonlocal-boundary-layer} and \Cref{lem:thin-layer-cutoff-stability} in \Cref{sec:NonlocalLaplaceI}.

Throughout the proofs, \(C\) may depend on \(m,d,\Om,w,\del_0\)
and the mesh regularity constants, but is independent of \(\del\), \(h\),
and the test functions. Its value may change from line to line.

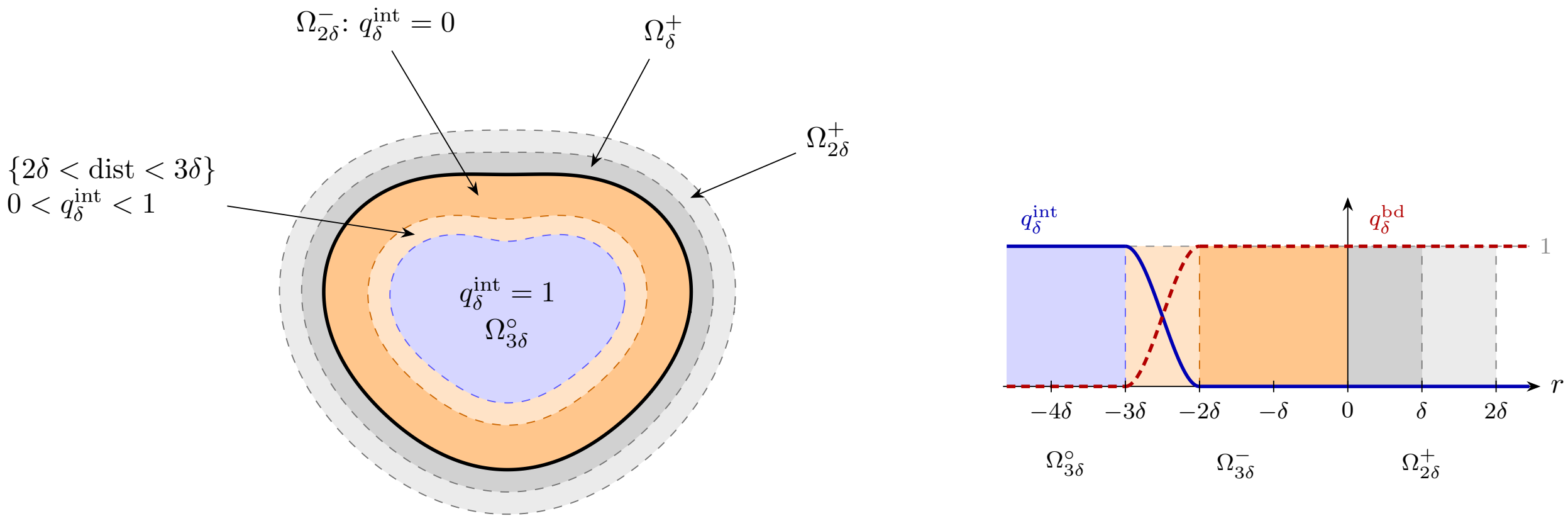
\begin{figure}[htbp]
\centering
\resizebox{\textwidth}{!}{%
\begin{tikzpicture}[>=Stealth,line join=round]
  \def\bl#1#2{2.0+0.21*cos(2*(#1))+0.12*sin(3*(#1))+(#2)}
  \def\dl{0.27}

  %------------------------------- panel (a): geometry
  \fill[black!8]
    plot[variable=\t,domain=0:360,samples=200,smooth cycle]
    (\t:{\bl{\t}{2*\dl}});
  \fill[black!18]
    plot[variable=\t,domain=0:360,samples=200,smooth cycle]
    (\t:{\bl{\t}{\dl}});
  \fill[orange!45]
    plot[variable=\t,domain=0:360,samples=200,smooth cycle]
    (\t:{\bl{\t}{0}});
  \fill[orange!22]
    plot[variable=\t,domain=0:360,samples=200,smooth cycle]
    (\t:{\bl{\t}{-2*\dl}});
  \fill[blue!16]
    plot[variable=\t,domain=0:360,samples=200,smooth cycle]
    (\t:{\bl{\t}{-3*\dl}});

  \draw[dashed,black!55]
    plot[variable=\t,domain=0:360,samples=200,smooth cycle]
    (\t:{\bl{\t}{2*\dl}});
  \draw[dashed,black!55]
    plot[variable=\t,domain=0:360,samples=200,smooth cycle]
    (\t:{\bl{\t}{\dl}});
  \draw[very thick]
    plot[variable=\t,domain=0:360,samples=200,smooth cycle]
    (\t:{\bl{\t}{0}});
  \draw[dashed,orange!80!black]
    plot[variable=\t,domain=0:360,samples=200,smooth cycle]
    (\t:{\bl{\t}{-2*\dl}});
  \draw[dashed,blue!65]
    plot[variable=\t,domain=0:360,samples=200,smooth cycle]
    (\t:{\bl{\t}{-3*\dl}});

  \node[align=center] at (0,0)
    {$q^{\mathrm{int}}_{\delta}=1$\\[1pt]
     $\Omega^{\circ}_{3\delta}$};

  \coordinate (B1) at (140:{\bl{140}{-2.5*\dl}});
  \coordinate (B2) at (104:{\bl{104}{-\dl}});
  \coordinate (B3) at (60:{\bl{60}{0.5*\dl}});
  \coordinate (B4) at (32:{\bl{32}{1.5*\dl}});

  \node[align=left] (LB1) at (-4.8,1.5)
    {$\{2\delta<\operatorname{dist}<3\delta\}$\\
     $0<q^{\mathrm{int}}_{\delta}<1$};
  \node[align=left] (LB2) at (-1.6,3.5)
    {$\Omega^{-}_{2\delta}$:\ $q^{\mathrm{int}}_{\delta}=0$};
  \node (LB3) at (1.9,3.4) {$\Omega^{+}_{\delta}$};
  \node (LB4) at (3.9,2.1) {$\Omega^{+}_{2\delta}$};

  \draw[->,thin] (LB1) -- (B1);
  \draw[->,thin] (LB2) -- (B2);
  \draw[->,thin] (LB3) -- (B3);
  \draw[->,thin] (LB4) -- (B4);

  \node at (0,-3.9) {(a)};

  %------------------------------- panel (b): signed-distance profiles
  \begin{scope}[xshift=10.2cm,yshift=-0.9cm]
    \def\ux{0.90}
    \def\uy{1.7}

    % Regions matching panel (a).
    \fill[blue!16]   (-4.6*\ux,0) rectangle (-3*\ux,\uy);
    \fill[orange!22] (-3*\ux,0)   rectangle (-2*\ux,\uy);
    \fill[orange!45] (-2*\ux,0)   rectangle (0,\uy);
    \fill[black!18]  (0,0)        rectangle (1*\ux,\uy);
    \fill[black!8]   (1*\ux,0)    rectangle (2*\ux,\uy);

    \draw[->] (-4.65*\ux,0) -- (2.6*\ux,0)
      node[right] {$r$};
    \draw[->] (0,0) -- (0,1.35*\uy);
    \draw[dashed,black!45] (-4.6*\ux,\uy) -- (2.45*\ux,\uy)
      node[right,font=\footnotesize] {$1$};

    \draw[dashed,blue!65] (-3*\ux,0) -- (-3*\ux,\uy);
    \draw[dashed,orange!80!black] (-2*\ux,0) -- (-2*\ux,\uy);
    \draw[dashed,black!55] (1*\ux,0) -- (1*\ux,\uy);
    \draw[dashed,black!55] (2*\ux,0) -- (2*\ux,\uy);

    % Interior cutoff: solid blue.
    \draw[very thick,blue!70!black]
      (-4.6*\ux,\uy) -- (-3*\ux,\uy);
    \draw[very thick,blue!70!black]
      plot[variable=\x,domain=-3:-2,samples=80]
      ({\x*\ux},{(1-3*(\x+3)^2+2*(\x+3)^3)*\uy});
    \draw[very thick,blue!70!black]
      (-2*\ux,0) -- (2.45*\ux,0);

    \node[blue!70!black,font=\footnotesize,anchor=south west]
      at (-4.55*\ux,1.02*\uy) {$q_\delta^{\mathrm{int}}$};

    % Boundary cutoff: dashed red.
    \draw[very thick,red!70!black,densely dashed]
      (-4.6*\ux,0) -- (-3*\ux,0);
    \draw[very thick,red!70!black,densely dashed]
      plot[variable=\x,domain=-3:-2,samples=80]
      ({\x*\ux},{(3*(\x+3)^2-2*(\x+3)^3)*\uy});
    \draw[very thick,red!70!black,densely dashed]
      (-2*\ux,\uy) -- (2.45*\ux,\uy);

    \node[red!70!black,font=\footnotesize,anchor=south west]
      at (0.15*\ux,1.02*\uy) {$q_\delta^{\mathrm{bd}}$};

    \foreach \v/\l in
      {-4/{-4\delta},-3/{-3\delta},-2/{-2\delta},-1/{-\delta},
       0/{0},1/{\delta},2/{2\delta}}{
      \draw (\v*\ux,2pt) -- (\v*\ux,-2pt)
        node[below,font=\footnotesize] {$\l$};
    }

    \node[font=\footnotesize] at (-3.8*\ux,-0.55*\uy)
      {$\Omega^{\circ}_{3\delta}$};
    \node[font=\footnotesize] at (-1.5*\ux,-0.55*\uy)
      {$\Omega^{-}_{3\delta}$};
    \node[font=\footnotesize] at (1*\ux,-0.55*\uy)
      {$\Omega^{+}_{2\delta}$};

    \node at (-1*\ux,-3.0) {(b)};
  \end{scope}
\end{tikzpicture}}
\caption{The type~(II) cutoffs.
(a) The interior region, transition layer, and two-layer exterior collar.
(b) Profiles of \(q_\delta^{\mathrm{int}}\) (solid) and
\(q_\delta^{\mathrm{bd}}\) (dashed) against the signed distance
\(r\) to \(\partial\Omega\), negative inside \(\Omega\)
and positive outside.}
\label{fig:cutoff-typeII}
\end{figure}

For the cutoff functions, we place the transition layer one horizon
farther inside \(\Om\) than in \Cref{sec:NonlocalLaplaceI}. We choose \( q_\del^{\mathrm{int}}\in C_c^\infty(\Om)\) such that
\[
0\le q_\del^{\mathrm{int}}\le1,
\quad
 q_\del^{\mathrm{int}}=1\text{ on }\Omega_{3\del}^o,
\quad
\supp q_\del^{\mathrm{int}}\subset\overline{\Omega_{2\del}^o},
\quad
|\nabla q_\del^{\mathrm{int}}|\le C\del^{-1},
\]
and set \( q_\del^{\mathrm{bd}}:=1- q_\del^{\mathrm{int}}\). The illustration is in \Cref{fig:cutoff-typeII} and it should be compared with \Cref{fig:cutoff-typeI}.

% =====================================================================
%  FIGURE 3 --- type (II) cutoff and the two-layer domain ladder.
% =====================================================================

\begin{lemma}
\label{lem:halfspace-localization-bridge}
Let \(v\in\cV_\del\). Its zero extension satisfies
\begin{equation*}
\begin{aligned}
E_0v&\in\cS_{w_\del}^{\nub}(\R^d),\\
\fG_{w_\del,\R^d}^{\nub}(E_0v)
&=
\begin{cases}
\fG_{w_\del}^{\nub}v&\text{a.e. in }\Omh,\\
0&\text{a.e. in }\R^d\setminus\Omh,
\end{cases}\\
\norm{E_0v}_{\cS_{w_\del}^{\nub}(\R^d)}
&=\norm{v}_{\cV_\del}.
\end{aligned}
\end{equation*}
For either \( q_\del= q_\del^{\mathrm{int}}\) or \( q_\del= q_\del^{\mathrm{bd}}\),
\begin{equation*}
\begin{gathered}
 q_\del v\in\cV_\del,
\qquad
 q_\del E_0v\in\cS_{w_\del}^{\nub}(\R^d),\\
\fG_{w_\del,\R^d}^{\nub}( q_\del E_0v)
=
\begin{cases}
\fG_{w_\del}^{\nub}( q_\del v)&\text{a.e. in }\Omh,\\
0&\text{a.e. in }\R^d\setminus\Omh,
\end{cases}\\
\abs{ q_\del E_0v}_{\cS_{w_\del}^{\nub}(\R^d)}
=\abs{ q_\del v}_{\cV_\del}.
\end{gathered}
\end{equation*}
Moreover,
\begin{equation}\label{eq:halfspace-localized-gradient-vanishing}
\begin{aligned}
\fG_{w_\del,\R^d}^{\nub}( q_\del^{\mathrm{int}} E_0v)&=0
&&\text{a.e. in }\R^d\setminus\overline{\Omin},\\
\fG_{w_\del,\R^d}^{\nub}( q_\del^{\mathrm{bd}} E_0v)&=0
&&\text{a.e. in }\R^d\setminus\overline{(\partial\Om)_\del^{4,1}}.
\end{aligned}
\end{equation}
\end{lemma}
\begin{proof}
The first three assertions follow from the zero-extension isometry in \Cref{lem:halfspace-zero-extension-isometry}. Since \( q_\del\in C^\infty(\R^d)\cap W^{1,\infty}(\R^d)\), by \Cref{lem:halfspace-product-rule}, \( q_\del E_0v\in \cS_{w_\del}^{\nub}(\R^d)\). 
Applying \Cref{lem:halfspace-zero-extension-isometry}   proves the second part of \Cref{lem:halfspace-localization-bridge}.

Finally, \( q_\del^{\mathrm{int}} E_0v\) vanishes a.e. outside \(\overline{\Omega_{2\del}^\circ}\), whereas \( q_\del^{\mathrm{bd}} E_0v\) vanishes a.e. outside \(\overline{\Omega_{3\del}^-}\). Then \Cref{lem:halfspace-finite-range-vanishing} and the inclusions
\[
    \Omega_{2\del}^\circ+B_\del(\bm0)\subset\Omega_\del^\circ,
    \qquad
    \Omega_{3\del}^-+B_\del(\bm0)\subset(\partial\Om)_\del^{4,1}
\]
give \eqref{eq:halfspace-localized-gradient-vanishing}. These inclusions can be read off \Cref{fig:cutoff-typeII}(a).
\end{proof}

Let \(\Rdel\) be given by \eqref{eq:halfspace-Rdel}. The preceding lemma shows that the following localized residuals are well-defined for every \(v\in\cV_\del\):
\begin{equation*}
    \Rint(v):=\Rdel( q_\del^{\mathrm{int}}v),
    \qquad
    \Rbd(v):=\Rdel( q_\del^{\mathrm{bd}}v).
\end{equation*}
Then \(\Rdel=\Rint+\Rbd\). By \eqref{eq:halfspace-localized-gradient-vanishing}, the nonlocal gradient of \( q_\del^{\mathrm{int}} v\) vanishes a.e. in \(\Omh\setminus\overline{\Omin}\), and hence
\begin{equation*}
    \Rint(v)
    =\int_\Om -\Delta U_0(\xb) q_\del^{\mathrm{int}}(\xb)v(\xb)\,d\xb
      -\int_{\Omin} \fG_{w_\del}^{\nub}\Uext(\xb)\cdot \fG_{w_\del}^{\nub}( q_\del^{\mathrm{int}} v)(\xb)\,d\xb,
\end{equation*}
and 
\begin{equation}\label{eq:halfspace-Rbd}
    \Rbd(v)=\underbrace{(-\Delta u_0,  q_\del^{\mathrm{bd}} v)_{L^2(\Omega)} - a_\del(E_0 u_0,  q_\del^{\mathrm{bd}} v)}_{\RbdI(v)} + \underbrace{(-\Delta G,  q_\del^{\mathrm{bd}} v)_{L^2(\Omega)} - a_\del( G,  q_\del^{\mathrm{bd}} v)}_{\RbdII(v)}. 
\end{equation}

We next record the localized representations used below. If \(\bar U_0\in H^{m+1}(\R^d)\) is a Sobolev extension of \(U_0\), then \(\bar U_0=\Uext=U_0\) on \(\Om\). For \(\xb\in\Omin\), the nonlocal gradient uses only values in \(B_\del(\xb)\subset\Om\), so \(\fG_{w_\del,\R^d}^{\nub}\bar U_0 = \fG_{w_\del}^{\nub}\Uext\) a.e. on \(\Omin\) by \Cref{lem:halfspace-locality}. Thus \Cref{lem:halfspace-localization-bridge} gives
\begin{equation*}
\Rint(v)
=-\int_{\R^d}\fG_{w_\del,\R^d}^{\nub}\bar U_0(\xb)\cdot
       \fG_{w_\del,\R^d}^{\nub}( q_\del^{\mathrm{int}} E_0v)(\xb)\,d\xb
  -\int_{\R^d}\Delta\bar U_0(\xb) q_\del^{\mathrm{int}}(\xb)E_0v(\xb)\,d\xb.
\end{equation*}
Choose a fixed Sobolev extension \(\bar G\in H^{m+1}(\R^d)\) of \(G\)
such that
\begin{equation*}
    \norm{\bar G}_{H^{m+1}(\R^d)}
    \le C\norm{G}_{H^{m+1}(\widehat\Om_{\del_0})},
    \qquad C=C(m,\Om,\del_0).
\end{equation*}
Its contribution has the analogous representation
\begin{equation*}
\RbdII(v)
=-\int_{\R^d}\fG_{w_\del,\R^d}^{\nub}\bar G(\xb)\cdot
       \fG_{w_\del,\R^d}^{\nub}( q_\del^{\mathrm{bd}} E_0v)(\xb)\,d\xb
  -\int_{\R^d}\Delta \bar G(\xb) q_\del^{\mathrm{bd}}(\xb)E_0v(\xb)\,d\xb.
\end{equation*}
Finally, since \(u_0\in H_0^1(\Om)\), the boundary part of \eqref{eq:halfspace-localized-gradient-vanishing} yields  
\begin{equation}\label{eq:halfspace-RbdI-gradient-localization}
a_\del(E_0u_0, q_\del^{\mathrm{bd}} v)
=\int_{(\partial\Om)_\del^{4,1}}
 \fG_{w_\del,\R^d}^{\nub}(E_0u_0)(\xb)\cdot
 \fG_{w_\del,\R^d}^{\nub}( q_\del^{\mathrm{bd}} E_0v)(\xb)\,d\xb.
\end{equation}

\subsubsection{Interior residual estimates} 

\begin{lemma}
\label{lem:halfspace-Fourier-consistency}
Let \(j\in \N_+\).  There exists \(C=C(j,d,w)>0\) such that for every
\(q\in H^{j+1}(\R^d)\) and \(\psi\in \cS^{\nub}_{w_\del}(\R^d)\),
\begin{equation*}
\begin{aligned}
\left|\int_{\R^d}
      \fG_{w_\del,\R^d}^{\nub}q(\xb)\cdot
      \fG_{w_\del,\R^d}^{\nub}\psi(\xb)\,d\xb
+\int_{\R^d}\Delta q(\xb) \psi(\xb)\,d\xb\right|
\le C\del^{\min\{2,j\}}
\norm{q}_{H^{j+1}(\R^d)}
\norm{\fG_{w_\del,\R^d}^{\nub}\psi}_{L^2(\R^d,\R^d)}.
\end{aligned}
\end{equation*}
\end{lemma}
\begin{proof}
By \eqref{eq:grad-Fourier} and Plancherel's theorem,
\[
\int_{\R^d}\fG_{w_\del,\R^d}^{\nub}q(\xb)\cdot
        \fG_{w_\del,\R^d}^{\nub}\psi(\xb)\,d\xb
      +\int_{\R^d}\Delta q(\xb)\psi(\xb)\,d\xb\qquad=\int_{\R^d}
        \left(\abs{\lab_{w_\del}^{\nub}(\xib)}^2-(2\pi)^2\abs{\xib}^2\right)
        \widehat q(\xib)\overline{\widehat\psi(\xib)}\,d\xib.
\]
Using Cauchy--Schwarz as in the proof of \Cref{lem:interior-consistency} reduces
the conclusion of \Cref{lem:halfspace-Fourier-consistency} to the multiplier estimate
\begin{equation*}
    \abs{\abs{\lab_{w_\del}^{\nub}(\xib)}^2 - (2\pi)^2 \abs{\xib}^2} \le C\del^{\min\{2,j\}} \left(1 + (2\pi)^2 \abs{\xib}^2\right)^{\frac{j+1}{2}} \abs{\lab_{w_\del}^{\nub}(\xib)}.
\end{equation*}
The scaling identity \eqref{eq:scaling-Fourier-symbol}, the estimate
for small frequencies in \Cref{lem:small-frequency-symbol}, and the
symbol bounds \eqref{eq:grad-Fourier-symbol-bounds}
give this estimate by the low and high frequency argument in the proof of
\Cref{lem:interior-consistency}.
\end{proof}

\begin{lemma}\label{lem:halfspace-Rint-bound}
    Let \(U_0\in H^{m+1}(\Om)\) and \(\Rint\) be the interior residual defined above. Then there exists \(C>0\) such that
    \begin{equation*}
        \abs{\Rint(v)} \le C\del^{\min\{2,m\}} \norm{U_0}_{H^{m+1}(\Om)} \norm{v}_{\cV_\del}, \quad \forall v\in \cV_\del.
    \end{equation*}
\end{lemma}

\begin{proof}
    Take a Sobolev extension \(\bar U_0\in H^{m+1}(\R^d)\) of \(U_0\). The whole-space representation of \(\Rint\) gives
    \[
        \abs{\Rint(v)}
    = \abs{\int_{\R^d} \fG_{w_\del,\R^d}^{\nub}\bar U_0(\xb)\cdot \fG_{w_\del,\R^d}^{\nub}( q_\del^{\mathrm{int}} E_0 v)(\xb)\,d\xb + \int_{\R^d} \Delta \bar U_0(\xb) q_\del^{\mathrm{int}}(\xb) E_0 v(\xb)\,d\xb}.
    \]
    Applying \Cref{lem:halfspace-Fourier-consistency} with \(\psi =  q_\del^{\mathrm{int}} E_0 v\) and \(j=m\), the extension estimate \(\norm{\bar U_0}_{H^{m+1}(\R^d)} \le C \norm{U_0}_{H^{m+1}(\Om)}\), and the first estimate in \Cref{lem:halfspace-cutoff-energy} with \( q_\del= q_\del^{\mathrm{int}}\) proves \Cref{lem:halfspace-Rint-bound}.
\end{proof}

As in \Cref{sec:NonlocalLaplaceI}, the second term \(\RbdII(v)\) in \eqref{eq:halfspace-Rbd} can be estimated by the same argument as \Cref{lem:halfspace-Rint-bound}. We record the result here for later use. 
    \begin{equation}\label{eq:halfspace-RbdII-bound}
    \abs{\RbdII(v)}\le C\del^{\min\{2,m\}}\norm{G}_{H^{m+1}(\widehat\Om_{\del_0})}
    \norm{v}_{\cV_\del},
    \qquad \forall v\in \cV_\del.
\end{equation}
Indeed, apply \Cref{lem:halfspace-Fourier-consistency} to the whole-space representation of \(\RbdII\) with \(q=\bar G\), \(\psi= q_\del^{\mathrm{bd}} E_0v\), and \(j=m\). Then the extension bound for \(\bar G\) and  the first estimate in \Cref{lem:halfspace-cutoff-energy} give \eqref{eq:halfspace-RbdII-bound}.

\subsubsection{Boundary residual and energy estimates} 
We obtain a similar result to \Cref{prop:nonlocal-boundary-layer} and postpone its proof to Appendix~\ref{sec:halfspace-analytic-foundations}.

\begin{proposition}
\label{prop:halfspace-boundary-layer}
For every fixed \(a\ge1\), there is \(C_a>0\) such
that for any \(0<\del\le\del_0/2\),
\begin{equation*}
    \norm{v}_{L^2(\Omega^-_{a\del})}
    \le C_a\del
       \abs{v}_{\cV_\del},
    \qquad \forall v\in\cV_\del.
\end{equation*}
\end{proposition}

\begin{lemma}
\label{lem:halfspace-cutoff-energy}
For \(v\in\cV_\del\), let \(E_0 v\) be its zero extension to \(\R^d\). Let \( q_\del\) denote either \( q_\del^{\mathrm{int}}\) or \( q_\del^{\mathrm{bd}}\). Then there is \(C>0\), independent of \(v\), such that
\begin{equation*}
\begin{aligned}
\norm{\fG_{w_\del,\R^d}^{\nub}
       ( q_\del E_0v)}_{L^2(\R^d,\R^d)}&\le C \norm{\fG_{w_\del,\R^d}^{\nub}(E_0 v)}_{L^2(\R^d,\R^d)}+ C\del^{-1} \norm{v}_{L^2(\Om_{5\del}^-)}\\
&\le
 C\norm{v}_{\cV_\del}.
\end{aligned}
\end{equation*}
If \(q\in H^1_0(\Om)\), then
\begin{equation*}
\norm{\fG_{w_\del,\R^d}^{\nub}
       ( q_\del^{\mathrm{bd}}E_0q)}_{L^2(\R^d,\R^d)}
\le C\bigl(\norm{\nabla q}_{L^2(\Omega^-_{4\del},\R^d)}
+\del^{-1}\norm{q}_{L^2(\Omega^-_{5\del})}\bigr).
\end{equation*}
\end{lemma}
\begin{proof}
    The first part of \Cref{lem:halfspace-localization-bridge} gives
    \[\norm{\fG_{w_\del,\R^d}^{\nub}(E_0 v)}_{L^2(\R^d,\R^d)}= \norm{\fG_{w_\del}^{\nub} v}_{L^2(\Omh,\R^d)} = \abs{v}_{\cV_\del}.\]
    For \( q_\del= q_\del^{\mathrm{int}}\), the product rule in \Cref{lem:halfspace-product-rule} gives
    \begin{equation*}
        \fG^{\nub}_{w_\del,\R^d}( q_\del E_0 v)(\xb) =  q_\del (\xb)\fG^{\nub}_{w_\del,\R^d}E_0 v(\xb) + \Rb( q_\del, E_0 v)(\xb),
    \end{equation*}
    where \begin{equation*}
        \Rb( q_\del, E_0 v) (\xb)= \int_{\R^d} \chi_{\nub}(\yb-\xb)\frac{\yb-\xb}{\abs{\yb-\xb}} w_\del(\yb- \xb)( q_\del(\yb) -  q_\del(\xb))E_0 v(\yb)d\yb.
    \end{equation*}
    For \( q_\del= q_\del^{\mathrm{bd}}\), the same identity follows by subtracting the formula for \( q_\del^{\mathrm{int}} E_0v\) from that for \(E_0v\); in particular, \(\Rb( q_\del^{\mathrm{bd}},E_0v)=-\Rb( q_\del^{\mathrm{int}},E_0v)\). Since \(\supp \nabla  q_\del \subset \overline{\Om_{3\del}^-\cap \Om_{2\del}^\circ}\) and \(\supp w_\del \subset \overline{B_\del(\bm 0)}\), the remainder \(\Rb( q_\del, E_0 v)\) vanishes a.e. outside \(\overline{(\partial\Om)_\del}\), where \((\partial\Om)_\del:=\Om_{4\del}^- \cap \Om_\del^\circ\).
    Thus, by the first-moment condition in \Cref{assu:halfspace-kernel} and the Lipschitz bound \(\abs{\nabla  q_\del}\le C\del^{-1}\),
    \[\begin{split}
        \int_{\R^d}\abs{\Rb( q_\del, E_0 v)(\xb) }^2d\xb  &\le C\del^{-2}\int_{(\partial\Om)_\del}\left(\int_{\R^d} \abs{\yb - \xb}w_\del(\yb - \xb) \abs{  v(\yb)} d\yb \right)^2 d\xb\\
        &\le C(2d)^2 \del^{-2} \norm{v}_{L^2(\Om_{5\del}^-)}^2.
    \end{split}\]
    Thus, taking the \(L^2(\R^d)\) norm on both sides of the product identity above and applying \Cref{prop:halfspace-boundary-layer}, we obtain the first estimate in \Cref{lem:halfspace-cutoff-energy}. It remains to show the second estimate of the lemma. Repeating the preceding estimate, the remainder is bounded by \(C\del^{-1}\norm{q}_{L^2(\Om_{5\del}^-)}\). The first term \( q_\del^{\mathrm{bd}}\fG^{\nub}_{w_\del,\R^d} E_0 q\) vanishes a.e. outside \(\overline{(\partial\Om)^{3,1}_\del}\), and hence
    \[
    \norm{ q_\del^{\mathrm{bd}} \fG_{w_\del,\R^d}^{\nub} E_0 q}_{L^2(\R^d,\R^d)} \le \norm{\fG_{w_\del,\R^d}^{\nub} E_0 q}_{L^2((\partial\Om)^{3,1}_\del)}\le C\norm{\nabla E_0 q}_{L^2((\partial\Om)^{4,2}_\del)} = C\norm{\nabla q}_{L^2(\Om_{4\del}^-)}.
    \]
    Here the second inequality follows from Minkowski's inequality,
    the first moment of \(w_\del\), and \Cref{lem:halfspace-H1-bound-domain} since
    \((\partial\Om)_\del^{3,1}+B_\del(\bm0)
    \subset(\partial\Om)_\del^{4,2}\).
    This proves the second estimate.
\end{proof}

\begin{lemma}
\label{lem:halfspace-Rbd-bound}
Let \(u_0=U_0-G|_\Om\in H^{m+1}(\Om)\), and suppose that
\eqref{eq:inhom-trace-compatibility} holds. Then
\[
    \Tr(\normal^j u_0)=0,
    \qquad j=0,\ldots,m-1.
\]
Let \(\Rbd\) and \(\RbdI\) be given by \eqref{eq:halfspace-Rbd}. Then
\begin{equation*}
    |\RbdI(v)|
    \le C\del^{m-\frac{1}{2}}
       \norm{u_0}_{H^{m+1}(\Om)}
       \abs{ q_\del^{\mathrm{bd}} E_0 v}_{\cS^{\nub}_{w_\del}(\R^d)}
    \le C\del^{m-\frac{1}{2}}
       \norm{u_0}_{H^{m+1}(\Om)}
       \norm{v}_{\cV_\del},
    \qquad \forall v\in\cV_\del.
\end{equation*}
\end{lemma}
\begin{proof}
    Decompose \(\RbdI\) into local and nonlocal parts:
\[
    \RbdI(v)=I(v)+II(v),
\]
where
\[
    I(v):=\int_\Om -\Delta u_0(\xb) q_\del^{\mathrm{bd}}(\xb)v(\xb)\,d\xb=\int_{\Omega^-_{3\del}} -\Delta u_0(\xb) q_\del^{\mathrm{bd}}(\xb)E_0 v(\xb)\,d\xb
\]
and
\[
II(v):=-a_\del(E_0 u_0,  q_\del^{\mathrm{bd}} v)=-\int_{\Omh}\fG_{w_\del}^{\nub} E_0 u_0(\xb)\cdot \fG_{w_\del}^{\nub}( q_\del^{\mathrm{bd}} v)(\xb)\,d\xb.
\]
We first estimate the nonlocal part \(II(v)\). Using
\eqref{eq:halfspace-RbdI-gradient-localization} and Cauchy--Schwarz,
followed by the same estimate from the first moment and \Cref{lem:halfspace-H1-bound-domain}, with
\((\partial\Om)_\del^{4,1}+B_\del(\bm0)
\subset(\partial\Om)_\del^{5,2}\), gives
\[\begin{split}
    \abs{II(v)}&= \abs{\int_{(\partial\Om)^{4,1}_\del} \fG_{w_\del,\R^d}^{\nub} E_0 u_0(\xb)\cdot \fG_{w_\del,\R^d}^{\nub}( q_\del^{\mathrm{bd}} E_0 v)(\xb)\,d\xb}\\
    &\le \norm{\fG_{w_\del,\R^d}^{\nub} E_0 u_0}_{L^2((\partial\Om)^{4,1}_\del,\R^d)} \norm{\fG_{w_\del,\R^d}^{\nub}( q_\del^{\mathrm{bd}} E_0 v)}_{L^2((\partial\Om)^{4,1}_\del,\R^d)}\\
    &\le C\norm{\nabla E_0 u_0}_{L^2((\partial\Om)^{5,2}_\del,\R^d)} \norm{\fG_{w_\del,\R^d}^{\nub}( q_\del^{\mathrm{bd}} E_0 v)}_{L^2((\partial\Om)^{4,1}_\del,\R^d)}.
\end{split}\]
 Since \(E_0 u_0\) has zero normal traces up to order \(m-1\) by
\eqref{eq:inhom-trace-compatibility},
\Cref{cor:local-boundary-double-layer-est}, with \(a=5\), \(b=2\),
\(S^-=u_0\), and \(S^+=0\), applies because \(2\del\le\del_0\).
It gives
\[\norm{\nabla E_0 u_0}_{L^2((\partial\Om)^{5,2}_\del,\R^d)}\le  C\del^{m-\frac{1}{2}}\norm{u_0}_{H^{m+1}(\Om)}.\]
Thus
\[
    |II(v)|
    \le C\del^{m-\frac{1}{2}}\norm{u_0}_{H^{m+1}(\Om)}
    \abs{ q_\del^{\mathrm{bd}} E_0 v}_{\cS^{\nub}_{w_\del}(\mathbb R^d)}.
\]
For the local part \(I(v)\),
\Cref{lem:halfspace-localization-bridge,prop:halfspace-boundary-layer} give
\[
    \norm{ q_\del^{\mathrm{bd}} v}_{L^2(\Omega^-_{3\del})}
    \le C\del\abs{ q_\del^{\mathrm{bd}} E_0v}_{\cS_{w_\del}^{\nub}(\R^d)}.
\]
With this bound and the zero traces above, the local calculation in the
proof of \Cref{lem:Rbd-bound} gives
\[
    \abs{I(v)}
    \le C\del^{m-\frac{1}{2}}\norm{u_0}_{H^{m+1}(\Om)}
       \abs{ q_\del^{\mathrm{bd}} E_0v}_{\cS_{w_\del}^{\nub}(\R^d)}.
\]
This and the estimate for \(II(v)\) prove the first inequality in
\Cref{lem:halfspace-Rbd-bound}. The second follows from
\Cref{lem:halfspace-cutoff-energy}. 
\end{proof} 

\begin{proof}[Energy estimate in \Cref{thm:LapII-rates-continuous}]
Note that \Cref{lem:halfspace-Rint-bound,lem:halfspace-Rbd-bound} yield
\[
|\Rbd(v)|\le C\del^{\min\{2,m-\frac{1}{2}\}}\left(\norm{u_0}_{H^{m+1}(\Om)} + \norm{G}_{H^{m+1}(\widehat\Om_{\del_0})}\right) \norm{v}_{\cV_\del},
    \qquad \forall v\in\cV_\del.
\] 
Then \Cref{thm:abstract-energy-rate} gives \eqref{eq:LapII-energy-rate}.
\end{proof}

\subsubsection{\(L^2\) estimate}

For \(\phi\in L^2(\Om)\), let \(z_\del[\phi]\in\cV_\del\) solve
\begin{equation*}
    a_\del(v,z_\del[\phi])
    =\int_\Om v(\xb)\phi(\xb)\,d\xb,
    \qquad \forall v\in\cU_\del,
\end{equation*}
and let \(z\) solve \eqref{eq:local-dual}.  The local regularity estimate
\eqref{eq:local-dual-reg} is used below.

Applying the energy estimate \eqref{eq:LapII-energy-rate} to the problem for \(z_\del[\phi]\) with \(m=1\) and \(G=0\), as well as \eqref{eq:local-dual-reg}, gives
\begin{equation}\label{eq:halfspace-dual-energy}
    \norm{z_\del[\phi]-E_0z}_{\cV_\del}
    \le C\del^{\frac{1}{2}}\norm{\phi}_{L^2(\Om)}.
\end{equation}

\begin{proposition}
\label{prop:halfspace-dual-boundary}
Extend \(z_\del[\phi]\) by zero to \(\mathbb R^d\). Then
\begin{equation*}
    \abs{ q_\del^{\mathrm{bd}} E_0z_\del[\phi]}_{\cS_{w_\del}^{\nub}(\R^d)}
    \le C\del^{1/2}\norm{\phi}_{L^2(\Om)}.
\end{equation*}
\end{proposition}
\begin{proof}
Using \eqref{eq:halfspace-dual-energy},
\Cref{lem:halfspace-cutoff-energy}, and \eqref{eq:local-dual-reg}, one may
estimate the two terms in
\(z_\del[\phi]=E_0z+(z_\del[\phi]-E_0z)\) as in the proof of
\Cref{prop:dual-boundary-layer} and obtain the desired estimate.
\end{proof}

\begin{proof}[\(L^2\) estimate in
\Cref{thm:LapII-rates-continuous}]
Applying the first inequality in \Cref{lem:halfspace-Rbd-bound} with \(v=z_\del[\phi]\) and \Cref{prop:halfspace-dual-boundary} gives
\begin{equation*}
\begin{aligned}
    \abs{\RbdI(z_\del[\phi])}
    &\le C\del^{m-\frac{1}{2}}\norm{u_0}_{H^{m+1}(\Om)}
       \abs{ q_\del^{\mathrm{bd}} E_0z_\del[\phi]}_{\cS_{w_\del}^{\nub}(\R^d)} \\
    &\le C\del^m\norm{u_0}_{H^{m+1}(\Om)}
       \norm{\phi}_{L^2(\Om)}.
\end{aligned}
\end{equation*}
Then \Cref{lem:halfspace-Rint-bound}
and \eqref{eq:halfspace-RbdII-bound}, together with dual stability,
now give the \(L^2\) estimate in \Cref{thm:LapII-rates-continuous} by \Cref{thm:abstract-L2-rate}.
\end{proof}

\subsubsection{Discrete estimates}

\begin{proof}[Proof of \Cref{thm:LapII-rates-discrete}]
Cauchy--Schwarz and \eqref{eq:halfspace-Poincare} give uniform continuity
and coercivity, hence unique solvability. The inclusion
\(V_{h,0}^{k}\subset\cU_{\del,h}=\cV_{\del,h}\), the uniform embedding
\eqref{eq:halfspace-uniform-H1-bound}, and
\eqref{eq:FE-H1-approximation} give the primal approximation factor
\(h^{\min\{m,k\}}\). Combining these estimates with
\eqref{eq:halfspace-dual-energy} gives the dual approximation factor
\(\del^{1/2}+h\). The residual estimates established above then yield
\eqref{eq:LapII-energy-rate-discrete} and
\eqref{eq:LapII-L2-rate-discrete} by the Galerkin argument in the proof of
\Cref{thm:LapI-rates-discrete}.
\end{proof}

\begin{remark}\label{rem:LapII-improved-rates-integrable-kernel}
If \(w\in L^1(\R^d)\), the convergence rates in    \eqref{eq:LapII-energy-rate-discrete} and \eqref{eq:LapII-L2-rate-discrete} can be improved by
\[
 \del^{\min\{2,m-1/2\}}+h^{\min\{m,k\}} \min\{1,\del^{-1}h\}
\]
and
\[
 \del^{\min\{2,m\}}+
 \bigl(\del^{1/2}+ h \min\{1,\del^{-1}h\}\bigr) h^{\min\{m,k\}} 
 \min\{1,\del^{-1}h\},
\]
respectively. This gives better estimates when $h < \delta$. The constant \(C\) may additionally depend on \(\norm{w}_{L^1}\), but is
independent of \(\del\) and \(h\).
\end{remark}

\begin{proof}[Proof of \Cref{prop:LapII-piecewise-constant}]
We first establish the half-space approximation estimate. For
\(0\le\alpha<\frac{1}{2}\), zero extension is bounded from
\(H^\alpha(\Om)\) to \(H^\alpha(\R^d)\). Minkowski's inequality and the
translation estimate used in \cite[Lemma~3.1]{tian2014asymptotically} give,
for \(v\in H^\alpha(\Om)\),
\[
\begin{aligned}
 \abs{E_0v}_{\cS_{w_\del}^{\nub}(\R^d)}
 &\le \int_{\R^d}\chi_{\nub}(\zb)w_\del(\zb)
    \norm{E_0v(\,\cdot+\zb)-E_0v}_{L^2(\R^d)}\,d\zb\\
 &\le C\del^{\alpha-1}
    \left(\int_{B_1(\bm0)}\abs{\zb}^\alpha w(\zb)\,d\zb\right)
    \norm{v}_{H^\alpha(\Om)}.
\end{aligned}
\]
For \(\alpha=0\), use
\(\norm{E_0v(\,\cdot+\zb)-E_0v}_{L^2(\R^d)}
\le2\norm{v}_{L^2(\Om)}\).
By the moment assumption in \Cref{prop:LapII-piecewise-constant}, adding the \(L^2\) term gives
\begin{equation*} 
 \norm{E_0v}_{\cU_\del}
 \le C\del^{\alpha-1}\norm{v}_{H^\alpha(\Om)},
 \qquad 0\le\alpha<\frac{1}{2}.
\end{equation*}
Since \(V_h^0\subset H^\alpha(\Om)\), this also proves
\(V_h^0\subset\cU_\del\).

For the  projection \(Q_h\) from \eqref{eq:Qh-definition},
applying \eqref{eq:fractional-approximation} yields, for
\(u\in H^1(\Om)\),
\begin{equation*}  
 \norm{E_0(u-Q_hu)}_{\cU_\del}
 \le C\del^{\alpha-1}\norm{u-Q_hu}_{H^\alpha(\Om)}
 \le C\left(\frac{h}{\del}\right)^{1-\alpha}
 \norm{u}_{H^1(\Om)}.
\end{equation*}
The rest of the proof follows as in the proof of \Cref{prop:LapI-piecewise-constant}.
\end{proof}

\section{Applications to nonlocal Laplacians with heterogeneous localization}
\label{sec:NonlocalLaplaceHeteroLoc}

In this section we consider nonlocal diffusion operators whose interaction
radius shrinks to zero as that point approaches the
boundary.  As a result, the boundary condition of this type of nonlocal model is in classical sense, compared to the nonlocal volumetric boundary conditions used in type (I) and type (II) nonlocal Laplacians. Earlier studies of such heterogeneous localization can be found in \cite{du2022fractional,tian2017trace,tao2019nonlocal,scott2024nonlocal,scott2025nonlocal}.

We first state the assumptions and variational problems, following the framework in \cite{scott2024nonlocal,scott2025nonlocal}.  We then give the
continuous and discrete convergence theorems.  The remainder of the section
is devoted to the consistency estimates and their proofs.

\subsection{Kernel, localization, and energy spaces}
\label{subsec:hetero-setting}

Throughout this section, \(\Om\subset\R^d\) is a bounded \(C^{2,1}\) domain.
Unless stated otherwise, \(C\) may depend on
\(\Om,d,\ka,\del_0,\ga\) and the fixed Robin, lifting, trace, Poincar\'e,
and elliptic regularity constants and, in the finite element estimates,
on \(m,k\) and the fixed mesh regularity constants, but not on \(\del\),
\(h\), or the functions being estimated. Throughout this section we assume
that \(\ga\) satisfies \Cref{assu:kernel} and belongs to \(L^\infty(\R^d)\).

\begin{assumption}
\label{ass:hetero-localization}
There is a constant \(\ka\ge1\) and a function
\(\eta:\overline\Om\to[0,\infty)\) such that
\begin{equation*}
    \frac1\ka\mathrm{dist}(\xb,\p\Om)
    \le \eta(\xb)
    \le \ka\mathrm{dist}(\xb,\p\Om),
    \qquad \xb\in\overline\Om,
\end{equation*}
and
\begin{equation*}
    \eta\in C^0(\overline\Om)\cap W^{2,\infty}(\Om),
    \qquad
    \norm{\eta}_{W^{2,\infty}(\Om)}\le \ka.
\end{equation*}
We also require the global Lipschitz bound
\begin{equation}\label{eq:hetero-localization-Lipschitz}
    \abs{\eta(\xb)-\eta(\yb)}
    \le\ka\abs{\xb-\yb},
    \qquad \xb,\yb\in\overline\Om.
\end{equation}
Moreover, the interior trace of its gradient satisfies
\begin{equation*}
    \Tr(\grad\eta)(\xb)=-\nb(\xb)
    \quad\text{for \(\sH^{d-1}\)-almost every \(\xb\in\p\Om\)},
\end{equation*}
where \(\nb\) is the outward unit normal.  
\end{assumption}

Such a localization function can be chosen for any Lipschitz domain; see for example the construction in
\cite[Section~8, pp.~4219]{scott2024nonlocal}.

\begin{figure}[htbp]
\centering
\begin{tikzpicture}[>=Stealth,line join=round]
\def\bl#1#2{2.0+0.21*cos(2*(#1))+0.12*sin(3*(#1))+(#2)}

\fill[blue!10]
  plot[variable=\t,domain=0:360,samples=200,smooth cycle]
  (\t:{\bl{\t}{0}});
\draw[very thick]
  plot[variable=\t,domain=0:360,samples=200,smooth cycle]
  (\t:{\bl{\t}{0}});

\foreach \a/\rr/\ba in {40/0.65/0.46, 15/1.40/0.28, -15/1.80/0.10}{
  \fill[orange!30,draw=orange!80!black] (\a:\rr) circle (\ba);
  \fill[black] (\a:\rr) circle (1.1pt);
}

\node[font=\footnotesize,anchor=west] at (2.75,1.75)
  {$B\bigl(\xb,\eta_{\delta}(\xb)\bigr)$};
\draw[thin] (2.72,1.72) -- (40:0.65);

% \draw[dashed,thin,black!60]
%   (-15:1.80) -- (-15:{\bl{-15}{0}});
% \node[font=\footnotesize,anchor=north west] at (-15:2.02)
%   {$\operatorname{dist}(x,\partial\Omega)$};
\node[align=center] at (0,-2.9)
{$\eta_{\delta}(\xb)=\delta\eta(\xb)
\asymp\delta\operatorname{dist}(\xb,\partial\Omega)$};
% \node[align=center] at (0,-2.9)
%   {$\eta_{\delta}(x)=\delta\eta(x)
%     \asymp\delta\operatorname{dist}(x,\partial\Omega)$:\\
%    every interaction ball stays inside $\Omega$, so no exterior collar is used};

\end{tikzpicture}
\caption{Heterogeneous localization with shrinking interaction balls.}
\label{fig:hetero-localization}
\end{figure}
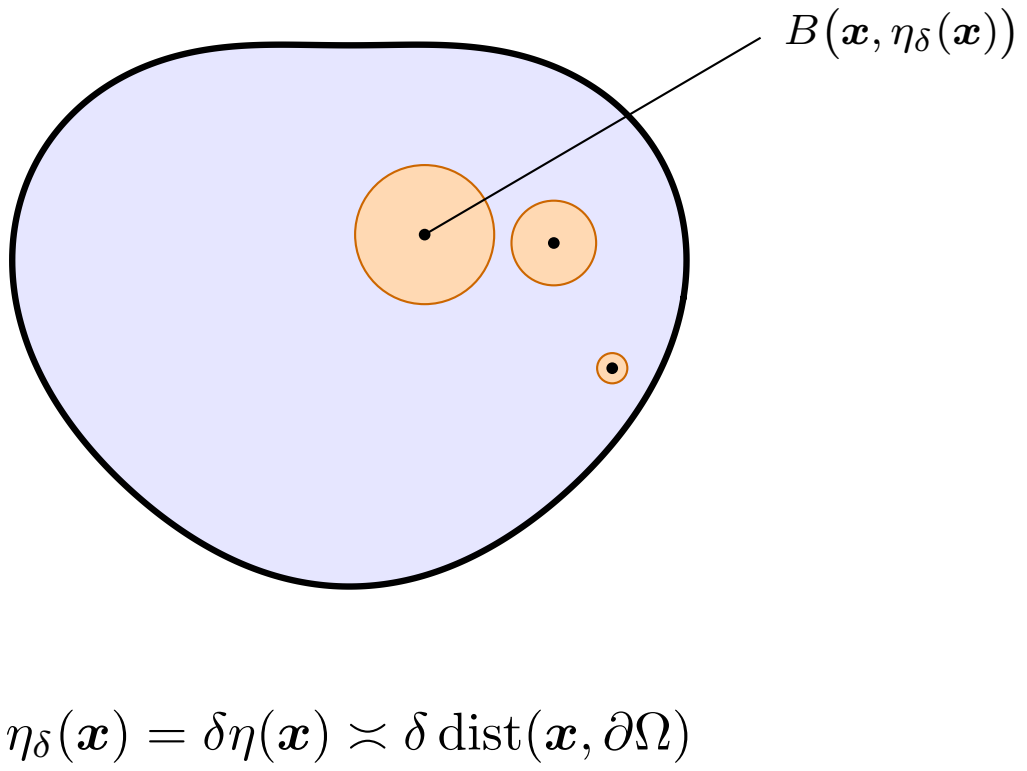

For \(0<\del<\del_0\), where
\begin{equation}\label{eq:hetero-admissible-delta}
    0<\del_0<\min\left\{1,\frac1{6\ka}\right\},
\end{equation}
set
\begin{equation*}
    \eta_{\delta}(\xb):=\del\eta(\xb).
\end{equation*}
Then
\begin{equation}\label{eq:hetero-ball-contained}
    \eta_{\delta}(\xb)<\dist(\xb,\p\Om),
    \qquad
    B(\xb,\eta_{\delta}(\xb))\subset\Om,
    \qquad \xb\in\Om.
\end{equation}
See \Cref{fig:hetero-localization}.

% =====================================================================
%  FIGURE 5 --- heterogeneous localization (Section 6).
% =====================================================================

These assumptions are specializations of the framework in \cite{scott2024nonlocal,scott2025nonlocal} where -- in the notation of those works --
\(p=2\), \(\beta=0\), \(q(r)=r\), \(\lambda=\eta\), with radial profile
\(r\mapsto\ga(r\bm e_1)\).
The nonlocal kernel and its symmetrization are
\begin{equation*}
\begin{aligned}
    \ga_{\del,\eta}(\xb,\yb)
    &:=\frac{1}{\eta_{\delta}(\xb)^{d+2}}
      \ga\left(\frac{\yb-\xb}{\eta_{\delta}(\xb)}\right),\\
    \widehat\ga_{\del,\eta}(\xb,\yb)
    &:=\ga_{\del,\eta}(\xb,\yb)
      +\ga_{\del,\eta}(\yb,\xb).
\end{aligned}
\end{equation*}
The associated operator and bilinear form are  
\begin{align*}
    \cL_{\del,\eta} U(\xb)
    &:=\int_\Om\widehat\ga_{\del,\eta}(\xb,\yb)
       \bigl(U(\yb)-U(\xb)\bigr)\,d\yb,
       \\
    a_{\del,\eta}(U,v)
    &:=\int_\Om\int_\Om\ga_{\del,\eta}(\xb,\yb)
       \bigl(U(\yb)-U(\xb)\bigr)
       \bigl(v(\yb)-v(\xb)\bigr)\,d\yb\,d\xb.
\end{align*}
The bilinear form can be rewritten using the symmetrized kernel as
\[
 a_{\del,\eta}(U,v)
 =\frac12\int_\Om\int_\Om
   \widehat\ga_{\del,\eta}(\xb,\yb)
   \bigl(U(\yb)-U(\xb)\bigr)
   \bigl(v(\yb)-v(\xb)\bigr)\,d\yb\,d\xb .
\]

Define 
\begin{equation}\label{eq:hetero-energy-space}
\begin{aligned}
    \cS_{\del}(\Om)
    &:={\left\{v\in L^2(\Om):[v]_{\cS_{\del}(\Om)}<\infty\right\}},\\
    [v]_{\cS_{\del}(\Om)}^2
    &:=\int_{\Om}\int_{B(\xb,\del\dist(\xb,\p\Om))}\frac{(v(\yb)-v(\xb))^2}{\left(\del \dist(\xb,\p\Om)\right)^{d+2}}\,d\yb\,d\xb,\\
    \norm{v}_{\cS_{\del}(\Om)}^2
    &:=\norm{v}_{L^2(\Om)}^2+[v]_{\cS_{\del}(\Om)}^2.
\end{aligned}
\end{equation}

By \cite[Theorem~2.2]{scott2024nonlocal}, there are constants \(c,C>0\),
independent of \(\del\), such that
\[
 c[v]_{\cS_{\del}(\Om)}^2
 \le a_{\del,\eta}(v,v)
 \le C[v]_{\cS_{\del}(\Om)}^2,
 \qquad v\in\cS_{\del}(\Om).
\]
By \cite[Corollary~5.3]{scott2024nonlocal}, for
\(0<\del<\del_0\),
\begin{equation}\label{eq:hetero-H1-embedding}
    \norm{v}_{\cS_{\del}(\Om)}
    \le C\norm{v}_{H^1(\Om)},
    \qquad v\in H^1(\Om).
\end{equation}
The trace theorem given by \cite[Theorem~1.3]{scott2024nonlocal} gives the
boundary trace on the heterogeneous space. It agrees with the Sobolev trace
on \(H^1(\Om)\), and we use \(\Tr\) for both operators:
\begin{equation}\label{eq:hetero-trace}
    \Tr:\cS_{\del}(\Om)\longrightarrow H^{1/2}(\p\Om),
    \qquad
    \norm{\Tr v}_{H^{1/2}(\p\Om)}
    \le C_{\mathrm{tr}}\norm{v}_{\cS_{\del}(\Om)}.
\end{equation}
The zero-trace, zero-mean, and boundary Poincar\'e inequalities used below
are \cite[Theorems~5.7, 5.8, and~5.10]{scott2024nonlocal}, respectively.

For \(U\in H^2(\Om)\) and every \(v\in\cS_{\del}(\Om)\), the weak
\(H^2\) extension of the nonlocal Green identity
\cite[Theorem~5.7, Remark~5.8]{scott2025nonlocal} is
\begin{equation}\label{eq:GreensIdentity}
    a_{\del,\eta}(U,v)
    =-\int_\Om \cL_{\del,\eta} U(\xb)v(\xb)\,d\xb
      +Q_\del\int_{\p\Om}
       \p_{\nb}U(\xb)(\Tr v)(\xb)\,d\sigma(\xb),
\end{equation}
where
\begin{equation}\label{eq:BoundaryConstant}
    Q_\del
    :=\int_{B_1(\bm0)}
       \frac{\operatorname{arctanh}(\del z_d)}{\del z_d}
       z_d^2\ga(\zb)\,d\zb.
\end{equation}

\subsection{Dirichlet, Neumann, and Robin problems}
\label{subsec:hetero-problems}

We follow the notation of
\Cref{sec:NonlocalLaplaceI,sec:NonlocalLaplaceII}.  The solutions are
denoted by \(U_\del\) and \(U_0\).  For the Dirichlet problem, choose a lifting \(G\in H^2(\Om)\) of the boundary datum and set
\begin{equation*}
    u_\del:=U_\del-G,
    \qquad
    u_0:=U_0-G.
\end{equation*}
For the Neumann and Robin problems, set \(G=0\), so that
\(u_\del=U_\del\) and \(u_0=U_0\).

The trial and test spaces are
\begin{equation}\label{eq:hetero-spaces}
    \cU_\del=\cV_\del:=
    \begin{cases}
       \cS_{\del,0}(\Om)
       :=\overline{C_c^\infty(\Om)}^{\,\cS_{\del}(\Om)},
           &\text{Dirichlet},\\
       \mathring{\cS}_{\del}(\Om)
       :=\{v\in\cS_{\del}(\Om):(v)_\Om=0\},
           &\text{Neumann},\\
       \cS_{\del}(\Om),&\text{Robin}.
    \end{cases}
\end{equation}
By \cite[Theorem~5.6]{scott2024nonlocal},
\(\cS_{\del,0}(\Om)=\ker\Tr\).
The limiting spaces are
\begin{equation}\label{eq:hetero-local-spaces}
    \cU_0=\cV_0:=
    \begin{cases}
       H_0^1(\Om),&\text{Dirichlet},\\
       \mathring H^1(\Om):=\{v\in H^1(\Om):(v)_\Om=0\},
           &\text{Neumann},\\
       H^1(\Om),&\text{Robin}.
    \end{cases}
\end{equation}

For the Robin problem, assume \(\mu>0\),
\(b\in W^{2,\infty}(\p\Om)\), \(b\ge0\), and \(b\not\equiv0\).  Define
\begin{equation*}
    a_\del(w,v):=
    \begin{cases}
       a_{\del,\eta}(w,v),&\text{Dirichlet or Neumann},\\
       a_{\del,\eta}(w,v)
       +Q_\del\displaystyle\int_{\p\Om}
          \frac{b(\xb)}{\mu}(\Tr w)(\xb)(\Tr v)(\xb)\,d\sigma(\xb),
           &\text{Robin},
    \end{cases}
\end{equation*}
and
\begin{equation*}
    a_0(w,v):=
    \begin{cases}
       \displaystyle\int_\Om
          \grad w(\xb)\cdot\grad v(\xb)\,d\xb,
           &\text{Dirichlet or Neumann},\\
       \displaystyle\int_\Om
          \grad w(\xb)\cdot\grad v(\xb)\,d\xb
       +\displaystyle\int_{\p\Om}
          \frac{b(\xb)}{\mu}(\Tr w)(\xb)(\Tr v)(\xb)\,d\sigma(\xb),
           &\text{Robin}.
    \end{cases}
\end{equation*}

We use the full norms on the spaces in
\eqref{eq:hetero-spaces}--\eqref{eq:hetero-local-spaces}. Let \(C_P\) denote the applicable
uniform constant in the zero-trace or zero-mean Poincar\'e inequality
\[
    \norm{v}_{L^2(\Om)}
    \le C_P[v]_{\cS_{\del}(\Om)}.
\]
For the Robin space, let \(C_{P,R}\) be the uniform constant in
\[
    \norm{v}_{L^2(\Om)}
    \le C_{P,R}a_\del(v,v)^{1/2}.
\]
The zero-trace and zero-mean inequalities are
\cite[Theorems~5.7 and~5.8]{scott2024nonlocal}.  For the Robin case,
choose \(\epsilon_b>0\) so that
\(\Gamma_b:=\{\xb\in\p\Om:b(\xb)\ge\epsilon_b\}\) has positive surface
measure; such a choice is possible because \(b\ge0\), \(b\not\equiv0\),
and \(b\) is continuous.  The boundary Poincar\'e inequality
\cite[Theorem~5.10]{scott2024nonlocal}, applied on \(\Gamma_b\), and
\[
  a_\del(v,v)
  \ge c[v]_{\cS_{\del}(\Om)}^2
      +\frac{\epsilon_b}{\mu}
       \norm{\Tr v}_{L^2(\Gamma_b)}^2
\]
give the stated Robin inequality, with a constant allowed to depend on
\(b\) and \(\mu\).  Since
\(\operatorname{arctanh}(s)/s\ge1\) and
\eqref{eq:gamma:2ndMoment} holds,
\[
    1\le Q_\del\le
    Q_{\max}:=\frac{\operatorname{arctanh}(\del_0)}{\del_0}.
\]
Consequently, there is a constant \(\alpha_a>0\), independent of \(\del\), such that
\begin{equation}\label{eq:hetero-full-norm-coercivity}
    a_\del(v,v)
    \ge \alpha_a\norm{v}_{\cS_{\del}(\Om)}^2,
    \qquad
    v\in\cU_\del.
\end{equation}
The energy comparison and the trace estimate \eqref{eq:hetero-trace} also give
\begin{equation}\label{eq:hetero-full-norm-continuity}
    \abs{a_\del(w,v)}
    \le M_a\norm{w}_{\cS_{\del}(\Om)}
             \norm{v}_{\cS_{\del}(\Om)},
\end{equation}
for a constant \(M_a>0\).
All these constants are independent of \(\del\).  The analogous local
Poincar\'e inequalities make \(a_0\) coercive in the full \(H^1(\Om)\)
norm on the corresponding local spaces.

The nonlocal problem is
\begin{equation}\label{eq:hetero-nonlocal-weak}
    \left\{
    \begin{aligned}
       &\text{Find }u_\del\in\cU_\del\text{ such that}\\
       &a_\del(u_\del,v)=\cF_\del(v),
       \qquad v\in\cV_\del,
    \end{aligned}
    \right.
\end{equation}
where
\begin{equation*}
    \cF_\del(v):=
    \begin{cases}
       (f,v)_{L^2(\Om)}-a_\del(G,v),&\text{Dirichlet},\\
       Q_\del\bigl((f,v)_{L^2(\Om)}+\langle g,\Tr v\rangle_{\p\Om}\bigr),
           &\text{Neumann},\\
       (f,v)_{L^2(\Om)}
       +Q_\del\langle g/\mu,\Tr v\rangle_{\p\Om},&\text{Robin}.
    \end{cases}
\end{equation*}
For the Neumann problem we assume
\begin{equation*}
    \int_\Om f(\xb)\,d\xb+\langle g,1\rangle_{\p\Om}=0.
\end{equation*}
Formally, the nonlocal Green identity \eqref{eq:GreensIdentity} relates
\eqref{eq:hetero-nonlocal-weak} to the equation
\(-\cL_{\del,\eta}U_\del=f_\del\) in \(\Om\), where
\begin{equation*} 
    f_\del:=
    \begin{cases}
       f,&\text{Dirichlet or Robin},\\
       Q_\del f,&\text{Neumann}.
    \end{cases}
\end{equation*}

The local problem is
\begin{equation}\label{eq:hetero-local-weak}
    a_0(u_0,v)=\cF_0(v),
    \qquad v\in\cV_0,
\end{equation}
with
\begin{equation*}
    \cF_0(v):=
    \begin{cases}
       (f,v)_{L^2(\Om)}-(\grad G,\grad v)_{L^2(\Om)},
           &\text{Dirichlet},\\
       (f,v)_{L^2(\Om)}+\langle g,\Tr v\rangle_{\p\Om},
           &\text{Neumann},\\
       (f,v)_{L^2(\Om)}+\langle g/\mu,\Tr v\rangle_{\p\Om},
           &\text{Robin}.
    \end{cases}
\end{equation*}
Thus the local solution satisfies
\begin{equation*}
    -\Delta U_0=f\quad\text{in }\Om,
\end{equation*}
together with \(\Tr U_0=g\) in the Dirichlet case,
\(\p_{\nb} U_0=g\) in the Neumann case, and
\(\mu\p_{\nb} U_0+b\Tr U_0=g\) in the Robin case.

The well-posedness results
\cite[Theorems~2.2, 2.3, and~2.4]{scott2025nonlocal} apply to the
Dirichlet, Neumann, and Robin problems, respectively; the uniform constants
are the ones displayed in
\eqref{eq:hetero-full-norm-coercivity}--\eqref{eq:hetero-full-norm-continuity}.
For the continuous result, we assume 
\begin{equation}\label{eq:hetero-data-assumptions}
    f\in L^2(\Om),
    \qquad
    g\in
    \begin{cases}
       H^{3/2}(\p\Om),&\text{Dirichlet},\\
       H^{1/2}(\p\Om),&\text{Neumann or Robin}.
    \end{cases}
\end{equation}
For the Dirichlet problem, the lifting can be chosen so that
\[
    \norm{G}_{H^2(\Om)}
    \le C\norm{g}_{H^{3/2}(\p\Om)}.
\]
The regularity of elliptic boundary problems gives
\begin{equation}\label{eq:hetero-local-H2-regularity}
 \norm{U_0}_{H^2(\Om)}+\norm{u_0}_{H^2(\Om)}
 \le C\left(
  \norm{f}_{L^2(\Om)}+
  \begin{cases}
   \norm{g}_{H^{3/2}(\p\Om)},&\text{Dirichlet},\\
   \norm{g}_{H^{1/2}(\p\Om)},&\text{Neumann or Robin}
  \end{cases}
 \right);
\end{equation}
see \cite[Section~2.5.1]{grisvard2011elliptic}.
We also use the local dual regularity estimate
\begin{equation}\label{eq:hetero-local-dual-H2-regularity}
    \norm{z_0[\phi]}_{H^2(\Om)}
    \le C\norm{\phi}_{L^2(\Om)},
\end{equation}
for the homogeneous local dual problem associated with \(a_0\).

\subsection{Main convergence results}
\label{subsec:hetero-main-results}

The error and residual follow the convention of the abstract framework:
\begin{equation*}
    \edel:=u_\del-u_0=U_\del-U_0,
    \qquad
    \Rdel(v):=\cF_\del(v)-a_\del(u_0,v).
\end{equation*}

Let
\(E:H^{1/2}(\p\Om)\to H^1(\Om)\) be an extension operator such that \(\Tr\circ E\) is the identity on \(H^{1/2}(\p\Om)\), where \(\Tr\) is the trace map.    For
\(\varphi\in H^{1/2}(\p\Om)\), define the approximate normal derivative  as in \cite[(9.1)]{scott2025nonlocal}:
\begin{equation}\label{eq:NormalDerivDef}
    \left\langle
       \left(\p_{\nb}U_\del\right)_\del,
       \varphi
    \right\rangle_{\p\Om}
    :=\frac{1}{Q_\del}
      \left(
        a_{\del,\eta}(U_\del,E\varphi)
        -(f_\del,E\varphi)_{L^2(\Om)}
      \right).
\end{equation} 

\begin{theorem}
\label{thm:hetero-rates-continuous}
Suppose that \Cref{assu:kernel,ass:hetero-localization} and
\eqref{eq:hetero-data-assumptions} hold, and that \(\ga\in L^\infty(\R^d)\).
Let \(u_\del\) and \(u_0\) solve
\eqref{eq:hetero-nonlocal-weak} and \eqref{eq:hetero-local-weak},
respectively.  Let \(m\ge1\) be an integer and suppose that
\(U_0\in H^{m+1}(\Om)\).  Then there is a constant \(C>0\), independent
of \(\del\), such that
\begin{align}
    \norm{u_\del-u_0}_{\cS_{\del}(\Om)}
    &\le C\del^{\min\{2,m\}}
       \norm{U_0}_{H^{m+1}(\Om)},
       \label{eq:hetero-energy-rate}\\
    \norm{U_\del-U_0}_{L^2(\Om)}
    &\le C\del^2
       \norm{U_0}_{H^{2}(\Om)},
       \label{eq:hetero-L2-rate}\\
    \norm{
       (\p_{\nb} U_\del)_\del-\p_{\nb} U_0
    }_{H^{-1/2}(\p\Om)}
    &\le C\del^{\min\{2,m\}}
       \norm{U_0}_{H^{m+1}(\Om)}.
       \label{eq:hetero-normal-rate}
\end{align} 
\end{theorem}

\begin{remark}[Sufficient data regularity]
\label{rem:hetero-sufficient-data-regularity}
Let \(m\ge1\) be an integer and suppose that \(\Om\) is a
\(C^{m,1}\) domain. In the Robin case, assume in addition that
\(b\in W^{m,\infty}(\p\Om)\). If
\[
 f\in H^{m-1}(\Om),
 \qquad
 g\in
 \begin{cases}
  H^{m+1/2}(\p\Om),&\text{Dirichlet},\\
  H^{m-1/2}(\p\Om),&\text{Neumann or Robin},
 \end{cases}
\]
and the Neumann compatibility condition holds when applicable, regularity gives
\[
 \norm{U_0}_{H^{m+1}(\Om)}
 \le C\left(
  \norm{f}_{H^{m-1}(\Om)}+
  \begin{cases}
   \norm{g}_{H^{m+1/2}(\p\Om)},&\text{Dirichlet},\\
   \norm{g}_{H^{m-1/2}(\p\Om)},&\text{Neumann or Robin}.
  \end{cases}
 \right)
\]
See \cite[Sections~2.2.2 and~2.5.1]{grisvard2011elliptic}.  Thus the
right-hand side in each estimate of
\Cref{thm:hetero-rates-continuous} may be bounded by the displayed data
norms whenever they are finite.  In the Dirichlet case, a bounded right
inverse of the trace may also be used to choose
\(G\in H^{m+1}(\Om)\) with
\[
 \norm{G}_{H^{m+1}(\Om)}
 \le C\norm{g}_{H^{m+1/2}(\p\Om)}.
\]
Consequently, \(u_0=U_0-G\) belongs to \(H^{m+1}(\Om)\), and the sum of
the \(H^{m+1}\)-norms of \(U_0\) and \(u_0\) is bounded by the same data
norms.  For the Neumann and Robin problems this conclusion is immediate
because \(G=0\).
\end{remark}

Let \(m\ge1\) be the integer in
\Cref{thm:hetero-rates-continuous} and assume in addition
that \(\Om\) is a \(C^{m,1}\) domain. Fix a polynomial
degree \(k\ge1\). Let \(\{\mathcal T_h\}_{0<h\le h_0}\)
and \(F_K:\widehat K\to K\) be given as in \Cref{subsec:LapI-discrete-results}. Define
\[
 V_h^{k}
 :=
 \left\{
 v_h\in H^1(\Om):
 (v_h|_K)\circ F_K\in\P_k(\widehat K)
 \text{ for every }K\in\mathcal T_h
 \right\}.
\]
With \(\cU_0\) defined in \eqref{eq:hetero-local-spaces},
the corresponding constrained spaces are
\[
 V_h^{k}\cap\cU_0
 =
 \begin{cases}
 V_{h,0}^{k},&\text{Dirichlet},\\
 \{v_h\in V_h^{k}:(v_h)_\Om=0\},&\text{Neumann},\\
 V_h^{k},&\text{Robin}.
 \end{cases}
\]
Let \(\cU_{\del,h}=\cV_{\del,h}\subset\cU_\del\)
be conforming finite element spaces satisfying
\[
 V_h^{k}\cap\cU_0\subset\cU_{\del,h}.
\]

For the Dirichlet problem, use the interpolant in
\eqref{eq:FE-H1-approximation}. For the Robin and Neumann
problems,
\cite[Theorem~4.1 and Corollary~4.1]{bernardi1989optimal},
with the same extension to
\(C^{m,1}\) domains as in
\Cref{subsec:LapI-discrete-results}, provide an operator
\(J_h:H^1(\Om)\to V_h^{k}\) which satisfies \eqref{eq:FE-H1-approximation}.
For the Robin problem, set \(I_hv=J_hv\).
For the Neumann problem, set
\[
 I_hv=J_hv-(J_hv)_\Om.
\]
Then \(I_h v\in V^k_h\cap \cU_0\) if \(v\in \cU_0\). Since \(v-I_h v=v-J_h v - (v - J_h v)_{\Om}\) for \((v)_{\Om}=0\), \(\norm{v-I_h v}_{L^2(\Om)}\le \norm{v-J_h v}_{L^2(\Om)}\) and \(\norm{\nabla(v-I_h v)}_{L^2(\Om)}=\norm{\nabla(v-J_h v)}_{L^2(\Om)}\). 
Thus, in all three cases, there is an operator
\(I_h:\cU_0\to V_h^{k}\cap\cU_0\subset\cU_{\del,h}\)
such that
\begin{equation*}
 \norm{v-I_hv}_{H^1(\Om)}
 \le Ch^r\norm{v}_{H^{r+1}(\Om)},
 \qquad
 v\in H^{r+1}(\Om)\cap\cU_0,
 \quad r=1,\ldots,\min\{m,k\},
\end{equation*}
where \(C\) is independent of \(h\) and \(\del\).

The discrete problem is
\begin{equation*}
    a_\del(u_{\del,h},v_h)=\cF_\del(v_h),
    \qquad v_h\in\cV_{\del,h}.
\end{equation*}
Denote \(U_{\del,h}=u_{\del,h}+G\).

\begin{theorem}
\label{thm:hetero-rates-discrete}
Let \(m\ge1\) be the integer in \Cref{thm:hetero-rates-continuous}, and
suppose in addition that \(u_0\in H^{m+1}(\Om)\).  Under the preceding
finite element assumptions, there is a constant \(C>0\), independent of
\(\del\) and \(h\), such that
\begin{align}
    \norm{u_{\del,h}-u_0}_{\cS_{\del}(\Om)}
    &\le C\bigl(\del^{\min\{2,m\}}+h^{\min\{k,m\}}\bigr)
       \bigl(\norm{U_0}_{H^{m+1}(\Om)}
             +\norm{u_0}_{H^{m+1}(\Om)}\bigr),
       \label{eq:hetero-energy-rate-discrete}\\
    \norm{U_{\del,h}-U_0}_{L^2(\Om)}
    &\le C\bigl(\del^{2}+h^{\min\{k,m\}+1}\bigr)
       \bigl(\norm{U_0}_{H^{m+1}(\Om)}
             +\norm{u_0}_{H^{m+1}(\Om)}\bigr).
       \label{eq:hetero-L2-rate-discrete}
\end{align}
Under the sufficient data assumptions in
\Cref{rem:hetero-sufficient-data-regularity}, the last factor in both
estimates may be replaced by the data norms displayed there.
\end{theorem}

\subsection{Consistency estimates}
\label{subsec:hetero-consistency}

We begin with two elementary consequences of the assumptions.

\begin{lemma}
\label{lem:hetero-Q-expansion}
There is a constant \(C\), independent of \(\del\), such that
\begin{equation}\label{eq:hetero-Q-expansion}
    Q_\del
    =1+\frac{\del^2}{3}
       \int_{B_1(\bm0)}z_d^4\ga(\zb)\,d\zb+O(\del^4),
    \qquad
    \abs{Q_\del-1}\le C\del^2.
\end{equation}
In particular, \eqref{eq:hetero-admissible-delta} gives
\(1\le Q_\del\le 3/2\).
\end{lemma}

\begin{proof}
For \(\abs{s}<1\),
\[
    \frac{\operatorname{arctanh}s}{s}
    =1+\frac{s^2}{3}+O(s^4).
\]
Placing this expansion into the definition of $Q_\delta$ in \eqref{eq:BoundaryConstant} and using  
\eqref{eq:gamma:2ndMoment}
proves the result.
\end{proof}

For \(t\in[0,1]\) and \(\zb\in B_1(\bm0)\) define
\begin{equation*}
    \bsPhi_{t,\zb}(\xb):=\xb+t\eta_{\delta}(\xb)\zb.
\end{equation*}
The following change-of-variables facts will be used repeatedly.

\begin{lemma}
\label{lem:hetero-coordinate-map}
The map \(\bsPhi_{t,\zb} : \overline{\Om} \to \overline{\Om} \) is bi-Lipschitz,
with inverse denoted by \(\bsPsi_{t,\zb} : \overline{\Om} \to \overline{\Om}\), and
\begin{equation}\label{eq:hetero-Jacobian}
    J_{t,\zb}(\xb)
    :=\det \grad \bsPhi_{t,\zb}(\xb)
    =1+t\grad \eta_{\delta}(\xb)\cdot\zb,
    \qquad \frac23\le J_{t,\zb}(\bx) \le \frac43, \qquad \forall \bx \in \Omega.
\end{equation}
Consequently, for \(1\le p\le\infty\),
\begin{equation}\label{eq:hetero-composition-bound}
    \norm{w\circ\bsPhi_{t,\zb}}_{L^p(\Om)}
    +\norm{w\circ\bsPsi_{t,\zb}}_{L^p(\Om)}
    \le C\norm{w}_{L^p(\Om)}.
\end{equation}
All constants are uniform in \(t,\zb\), and \(\del\).
\end{lemma}

\begin{proof}
By \eqref{eq:hetero-localization-Lipschitz},
\(\eta_{\delta}\) is Lipschitz on \(\overline\Om\) with constant
\(\del\ka<1/3\).  Hence
\[
 \left| \abs{\bsPhi_{t,\zb}(\xb)-\bsPhi_{t,\zb}(\yb)}
 -\abs{\xb-\yb} \right|
 \leq t\abs{\zb}\del\ka
      \abs{\xb-\yb}
 \leq \frac13\abs{\xb-\yb}, \qquad \bx, \by \in \overline{\Omega}.
\]
A direct computation with the same type of estimate proves \eqref{eq:hetero-Jacobian}.
Thus $\bsPhi_{t,\zb}$ is bi-Lipschitz on $\overline{\Omega}$,  hence injective on $\overline{\Omega}$, and $\bsPhi_{t,\zb}(\bx) = \bx$ on $\p \Omega$, hence is invertible on $\p \Omega$. Moreover, the Jacobian of $\bsPhi_{t,\zb}$ is positive on $\Omega$. Therefore $\bsPhi_{t,\zb}$ is globally invertible as a map from $\overline{\Om}$ onto its image $\bsPhi_{t,\zb}(\overline{\Om})$.

It remains to show that $\bsPhi_{t,\zb}(\overline{\Om}) = \overline{\Om}$. Note that \(\bsPhi_{t,\zb}(\overline{\Om})\subset \overline{\Om}\) thanks to \eqref{eq:hetero-ball-contained}. First, $\bsPhi_{t,\zb}(\overline{\Omega})$ is nonempty, and it is closed since $\bsPhi_{t,\zb}$ is continuous. Moreover, because $\bsPhi_{t,\zb}$ is continuous, injective, and has uniformly positive Jacobian, $\bsPhi_{t,\zb}$ maps open sets to open sets, hence $\bsPhi_{t,\zb}(\Omega)$ is open. By extending $\bsPhi_{t,\zb}$ to a continuous function on $\R^d$ via $\bsPhi_{t,\zb}(\bx) = \bx$ for $\bx \in \R^d \setminus \overline{\Omega}$, we see that $\bsPhi_{t,\zb}(\overline{\Omega})$ is relatively open in $\overline{\Omega}$. Therefore since $\Omega$ is connected it must be that \(\bsPhi_{t,\zb}(\overline{\Om}) = \overline{\Om}\).

Finally, a straightforward change of variables gives \eqref{eq:hetero-composition-bound}.
\end{proof}

We define another change of variables; fix \(\xb\in\Om\) and define
\begin{equation*}
 \bsupsilon_{\xb}^{\del}(\yb)
 :=\frac{\xb-\yb}{\eta_{\delta}(\yb)},
 \qquad
 B_{\xb}^{\del}
 :=\left\{\yb\in\Om:
          \abs{\bsupsilon_{\xb}^{\del}(\yb)}<1\right\}
  =\left\{\yb\in\Om:
          \abs{\xb-\yb}<\eta_{\delta}(\yb)\right\}.
\end{equation*}
This is the negative of the map used in
\cite[Lemmas~3.3--3.4]{scott2025nonlocal}, hence has similar properties. We record them below for completeness.

\begin{lemma}\label{lem:hetero-adjoint-coordinate-map}
For each \(\xb\in\Om\),
\[
 \bsupsilon_{\xb}^{\del}:B_{\xb}^{\del}\longrightarrow B_1(\bm0)
\]
is a bi-Lipschitz bijection, with inverse
\[
 \zb\longmapsto\bsPsi_{1,\zb}(\xb).
\]
For almost every \(\yb\in B_{\xb}^{\del}\),
\begin{equation}\label{eq:hetero-adjoint-coordinate-Jacobian}
\begin{split}
 \left|\det \grad \bsupsilon_{\xb}^{\del}(\yb)\right|
 =\frac{\eta_{\delta}(\yb)+\grad \eta_{\delta}(\yb)\cdot(\xb-\yb)}
         {\eta_{\delta}(\yb)^{d+1}}, %
 \qquad
 \frac{2}{3\eta_{\delta}(\yb)^d}
 \le \left|\det \grad \bsupsilon_{\xb}^{\del}(\yb)\right|
 \le \frac{4}{3\eta_{\delta}(\yb)^d}.
\end{split}
\end{equation}
\end{lemma}

\begin{proof}
If \(\yb\in B_{\xb}^{\del}\) and
\(\zb=\bsupsilon_{\xb}^{\del}(\yb)\), then
\[
 \xb=\yb+\eta_{\delta}(\yb)\zb=\bsPhi_{1,\zb}(\yb),
\]
so \(\yb=\bsPsi_{1,\zb}(\xb)\).  Conversely, for every
\(\zb\in B_1(\bm0)\), \Cref{lem:hetero-coordinate-map} gives a unique
\(\yb=\bsPsi_{1,\zb}(\xb)\in\Om\).  The preceding identity then gives
\(\bsupsilon_{\xb}^{\del}(\yb)=\zb\), and hence
\(\yb\in B_{\xb}^{\del}\).  This proves that $\bsupsilon_{\bx}^\delta$ is a bijection.

Now, for almost every $\by \in B_{\bx}^\delta$
\[
 \grad \bsupsilon_{\xb}^{\del}(\yb)
 =-\eta_{\delta}(\yb)^{-1}
   \left(I+\bsupsilon_{\xb}^{\del}(\yb)
                   \otimes\grad \eta_{\delta}(\yb)\right).
\]
Then computing the determinant and using \(\norm{\grad \eta_{\delta}}_{L^\infty(\Om)}<1/3\) gives \eqref{eq:hetero-adjoint-coordinate-Jacobian}.
For \(\yb_1,\yb_2\in B_{\xb}^{\del}\), we have
\[
 \bsupsilon_{\xb}^{\del}(\yb_1)-\bsupsilon_{\xb}^{\del}(\yb_2)
 =\frac{\yb_2-\yb_1}{\eta_\delta(\yb_1)}
 +\frac{\eta_\delta(\yb_2)-\eta_\delta(\yb_1)}{\eta_\delta(\yb_1)}
   \bsupsilon_{\xb}^{\del}(\yb_2).
\]
The Lipschitz bound for \(\eta_\delta\), together with
\(\abs{\bsupsilon_{\xb}^{\del}(\yb_i)}<1\) and
\((1-\ka\del)\eta_\delta(\yb_i)\le\eta_\delta(\xb)
\le(1+\ka\del)\eta_\delta(\yb_i)\), gives
\[
 \frac{(1-\ka\del)^2}{\eta_\delta(\xb)}\abs{\yb_1-\yb_2}
 \le\abs{\bsupsilon_{\xb}^{\del}(\yb_1)-\bsupsilon_{\xb}^{\del}(\yb_2)}
 \le\frac{(1+\ka\del)^2}{\eta_\delta(\xb)}\abs{\yb_1-\yb_2}.
\]
Since \(\eta_\delta(\xb)>0\), this proves the bi-Lipschitz property.
\end{proof}

We next recall the boundary localized convolution.  Choose an even, radial
function
\(\psi\in C_c^\infty(B_{c_\ga/2}(\bm0))\), with \(\psi\ge0\), bounded
below by a positive constant on \(B_{c_\ga/4}(\bm0)\), and satisfying
\(\int_{B_1(\bm0)}\psi(\zb)\,d\zb=1\).  Set
\begin{equation}\label{eq:hetero-K-definition}
    K_\del v(\xb)
    :=\int_{B_1(\bm0)}
       \psi(\zb)v\bigl(\xb+\eta_{\delta}(\xb)\zb\bigr)\,d\zb.
\end{equation}

\begin{lemma}
\label{lem:hetero-boundary-convolution}
Under \Cref{assu:kernel,ass:hetero-localization}, there is a linear
operator \(K_\del:\cS_{\del}(\Om)\to H^1(\Om)\), defined by
\eqref{eq:hetero-K-definition}, such that
\begin{align}
    \norm{K_\del v}_{H^1(\Om)}
    &\le C\norm{v}_{\cS_{\del}(\Om)},
       \label{eq:hetero-boundary-convolution-H1}\\
    \norm{v-K_\del v}_{L^2(\Om)}
    &\le C\del\norm{v}_{\cS_{\del}(\Om)},
       \label{eq:hetero-boundary-convolution-L2}\\
    \Tr(K_\del v)&=\Tr v.
       \notag
\end{align}
\end{lemma}

\begin{proof}
See 
\cite[Theorem~1.1, Theorem~2.2, Corollary~4.8 and Theorem~5.4]{scott2024nonlocal}.
\end{proof}

The next lemma records the consistency between the nonlocal operator and its local limit.

\begin{lemma}
\label{lem:hetero-strong-consistency}
For every \(v \in H^2(\Om)\),
\begin{equation}\label{eq:hetero-operator-H2}
    \norm{\cL_{\del,\eta}v}_{L^2(\Om)}
    \le C\norm{v}_{H^2(\Om)}.
\end{equation}
For every \(v\in H^3(\Om)\),
\begin{equation}\label{eq:hetero-strong-first-order}
    \norm{\Delta v-\cL_{\del,\eta}v}_{L^2(\Om)}
    \le C\del\norm{v}_{H^3(\Om)}.
\end{equation}
\end{lemma}
\begin{proof}
The continuity $\cL_{\delta,\eta} : H^2(\Omega) \to L^2(\Omega)$, including its extension from
smooth functions and its identification with an absolutely convergent integral, is shown in
\cite[Theorem~5.7 and Remark~5.8]{scott2025nonlocal}.
It remains to prove
\eqref{eq:hetero-strong-first-order}, which we do first for
\(v\in C^3(\overline\Om)\). Write $\widehat\ga_{\del,\eta}(\xb,\yb) = \ga_{\del,\eta}(\xb,\yb) + \ga_{\del,\eta}(\yb,\xb)$, and split the operator
\(\cL_{\del,\eta}=\cL_{\del,\eta}^{1}+\cL_{\del,\eta}^{2}\) accordingly.

\emph{Step 1:}
For fixed \(\xb\), the change of variables $\bz = \frac{\by-\bx}{\eta_{\delta}(\bx)}$ gives
\begin{equation*}
    \cL_{\del,\eta}^{1}v(\xb)
    =\frac{1}{\eta_{\delta}(\xb)^2}
      \int_{B_1(\bm0)}\ga(\zb)
      \bigl(v(\xb+\eta_{\delta}(\xb)\zb)-v(\xb)\bigr)\,d\zb .
\end{equation*}
Taylor's formula for $v(\bx+\eta_{\delta}(\bx)\bz)$ gives
\[
\begin{split}
 \frac{v(\xb+\eta_{\delta}(\xb)\zb)-v(\xb)}{\eta_{\delta}(\xb)^2}
 ={}&\frac{\grad v(\xb)\cdot\zb}{\eta_{\delta}(\xb)}
    +\frac12D^2v(\xb)[\zb,\zb]\\
 &+\frac{\eta_{\delta}(\xb)}2\int_0^1(1-t)^2
 D^3v(\bsPhi_{t,\bz}(\bx))[\zb,\zb,\zb]\,dt.
\end{split}
\]
Inserting this into the integral gives three terms. The first term vanishes by the radial symmetry of $\gamma$, and the normalization property in \eqref{eq:gamma:2ndMoment} implies that the second term evaluates as
\(\tfrac12\Delta v(\xb)\). Thus
\[
\begin{split}
 \cL_{\del,\eta}^{1}v(\xb)-\frac12\Delta v(\xb)
 =\frac{\eta_{\delta}(\xb)}2
 \int_{B_1(\bm0)}\int_0^1
 &(1-t)^2\ga(\zb)
 D^3v(\bsPhi_{t,\bz}(\bx))[\zb,\zb,\zb]
 \,dt\,d\zb .
\end{split}
\]
Since \(\eta_{\delta}\le\ka\del\), integrating in $\bx$ and using
\eqref{eq:hetero-composition-bound} gives
\begin{equation*}
    \norm{\cL_{\del,\eta}^{1}v-\tfrac12\Delta v}_{L^2(\Om)}
    \le C\del\norm{v}_{H^3(\Om)}.
\end{equation*}

\emph{Step 2:}
Taylor expand $v(\by) - v(\bx)$ to get
\begin{equation*}
\begin{split}
    \cL_{\del,\eta}^{2}v(\bx) &= - \int_{B_{\bx}^\delta} \frac{ \gamma( \bsupsilon_{\bx}^\delta(\by) ) }{ \eta_{\delta}(\by)^{d+1} } \bsupsilon_{\bx}^\delta(\by) \, d \by \cdot \grad v(\bx) \\
    &\quad + \frac{1}{2} \int_{B_{\bx}^\delta} \frac{ \gamma( \bsupsilon_{\bx}^\delta(\by) ) }{ \eta_{\delta}(\by)^{d} } D^2 v(\bx) [ \bsupsilon_{\bx}^\delta(\by),\bsupsilon_{\bx}^\delta(\by) ] \, d \by \\
    &\quad - \int_{B_{\bx}^\delta} \int_0^1 \frac{(1-t)^2}{2} \frac{ \gamma( \bsupsilon_{\bx}^\delta(\by) ) }{ \eta_{\delta}(\by)^{d-1} } D^3 v(\bx + t(\by-\bx) ) [ \bsupsilon_{\bx}^\delta(\by),\bsupsilon_{\bx}^\delta(\by),\bsupsilon_{\bx}^\delta(\by) ] \, d t \, d \by := I + II + III.
\end{split}
\end{equation*}
We estimate each term separately. First, since $1 = \frac{\eta_{\delta}(\bx) - \eta_{\delta}(\by) - \grad \eta_{\delta}(\by) \cdot (\bx-\by)}{\eta_{\delta}(\bx)} + \frac{\eta_{\delta}(\by) + \grad \eta_{\delta}(\by) \cdot (\bx-\by)}{\eta_{\delta}(\bx)}$,
\begin{equation*}
\begin{split}
    I &=  - \int_{B_{\bx}^\delta} \gamma( \bsupsilon_{\bx}^\delta(\by) ) \frac{\eta_{\delta}(\bx) - \eta_{\delta}(\by) - \grad \eta_{\delta}(\by) \cdot (\bx-\by)}{\eta_{\delta}(\bx) \eta_{\delta}(\by)^{d+1} } \bsupsilon_{\bx}^\delta(\by) \, d \by \cdot \grad v(\bx) \\
    &\quad - \frac{1}{\eta_{\delta}(\bx)} \int_{B_{\bx}^\delta} \gamma( \bsupsilon_{\bx}^\delta(\by) ) \frac{\eta_{\delta}(\by) + \grad \eta_{\delta}(\by) \cdot (\bx-\by)}{ \eta_{\delta}(\by)^{d+1} } \bsupsilon_{\bx}^\delta(\by) \, d \by \cdot \grad v(\bx).
\end{split}
\end{equation*}
By changing variables $\bz =  \bsupsilon_{\bx}^\delta(\by)$ and noting the formula for its Jacobian in \eqref{eq:hetero-adjoint-coordinate-Jacobian}, the second term is seen to be zero. As for the first term, by the twice-differentiable assumption on $\eta_{\delta}$ we have $|\eta_{\delta}(\bx) - \eta_{\delta}(\by) - \grad \eta_{\delta}(\by) \cdot (\bx-\by)| \leq C \delta |\bx-\by|^2$. Then using the comparability of $\eta_{\delta}(\bx)$ and $\eta_{\delta}(\by)$ in the integral, as well as the lower bound on the Jacobian \eqref{eq:hetero-adjoint-coordinate-Jacobian}
\begin{equation}\label{eq-hetero-uniform-pf1}
    |I| \leq C(\kappa) \delta |\grad v(\bx)| \int_{ B_{\bx}^\delta } \gamma( \bsupsilon_{\bx}^\delta(\by) ) |\bsupsilon_{\bx}^\delta(\by)|^3 |\det \grad \bsupsilon_{\bx}^\delta(\by)| \, d \by  \leq C(\kappa,\gamma) \delta |\grad v(\bx)|,
\end{equation}
where we changed coordinates $\bz =  \bsupsilon_{\bx}^\delta(\by)$ in the last inequality to obtain the integral's independence of $\bx$.

We next estimate $II - \frac{1}{2} \Delta v$. Using \eqref{eq:gamma:2ndMoment} and the change of coordinates $\by = (\bsupsilon_{\bx}^\delta)^{-1}(\bz)$, we write 
\begin{equation*}
    \frac{1}{2} \Delta v(\bx) = \frac{1}{2} \int_{B_1({\bf 0})} \gamma(\bz) D^2 v(\bx)[\bz,\bz] \, d \bz =  \frac{1}{2} (-1)^d \int_{B_{\bx}^\delta} \gamma(\bsupsilon_{\bx}^\delta(\by)) D^2 v(\bx)[\bsupsilon_{\bx}^\delta(\by),\bsupsilon_{\bx}^\delta(\by)] \det \grad \bsupsilon_{\bx}^\delta(\by) \, d \by.
\end{equation*}
Note that $\eta_{\delta}(\by)^{-d} = \frac{(-1)^d}{1 + \grad \eta_{\delta}(\by) \cdot \bsupsilon_{\bx}^\delta(\by) } \det \grad \bsupsilon_{\bx}^\delta(\by)$ for all $\by \in B_{\bx}^\delta$ and $\bx \in \Omega$; we place this in II to get
\begin{equation*}
    \begin{split}
        &II - \frac{1}{2} \Delta v(\bx)\\
        =& \frac{(-1)^d }{2} \int_{B_{\bx}^\delta} \gamma(\bsupsilon_{\bx}^\delta(\by)) \left( \frac{1}{ 1 + \grad \eta_{\delta}(\by) \cdot \bsupsilon_{\bx}^\delta(\by) } - 1 \right) D^2 v(\bx)[\bsupsilon_{\bx}^\delta(\by),\bsupsilon_{\bx}^\delta(\by)] \det \grad \bsupsilon_{\bx}^\delta(\by) \, d \by.
    \end{split}
\end{equation*}
Now, $\left| \frac{1}{ 1 + \grad \eta_{\delta}(\by) \cdot \bsupsilon_{\bx}^\delta(\by) } - 1 \right| \leq \frac{\kappa \delta}{1-\kappa \delta}$ for all $ \by \in B_{\bx}^\delta$ and $\bx \in \Omega$, so
\begin{equation*}
    \begin{split}
        \left| II - \frac{1}{2} \Delta v(\bx) \right| \leq \frac{\delta \kappa }{2 (1-\kappa \delta)} \int_{B_{\bx}^\delta} \gamma(\bsupsilon_{\bx}^\delta(\by)) |\bsupsilon_{\bx}^\delta(\by)|^2  |\det \grad \bsupsilon_{\bx}^\delta(\by)| \, d \by |D^2 v(\bx)|.
    \end{split}
\end{equation*}
Changing coordinates again $\bz = \bsupsilon_{\bx}^\delta(\by)$ we obtain
\begin{equation}\label{eq-hetero-uniform-pf2}
    \left| II - \frac{1}{2} \Delta v(\bx) \right|  \leq C(\kappa,\gamma) \delta  |D^2 v(\bx)|.
\end{equation}

To begin estimating the third term, integrate in $\bx$ and use H\"older's inequality, along with the bound $\frac{1}{\eta_{\delta}(\by)^{d-1}} \leq \frac{\kappa \delta}{1-\kappa \delta} |\det \grad \bsupsilon_{\bx}^\delta(\by)|$:
\begin{equation*}
    \begin{split}
        &\vnorm{III}_{L^2(\Omega)}^2\\
        \leq& \frac{C \delta^2}{(1-\kappa \delta)^2} \int_{\Omega} \left( \int_{B_{\bx}^\delta} \gamma(\bsupsilon_{\bx}^\delta(\by)) |\bsupsilon_{\bx}^\delta(\by)|^4 |\det \grad \bsupsilon_{\bx}^\delta(\by)| \, d \by \right) \\
        &\qquad \qquad \cdot \left( \int_{\Omega} \1_{B_{\bx}^\delta}(\by)  \gamma(\bsupsilon_{\bx}^\delta(\by))  |\bsupsilon_{\bx}^\delta(\by)|^2 \left|  \int_0^1 \frac{(1-t)^2}{2} |D^3 v(\bx+t(\by-\bx))| \, d t \right|^2 |\det \grad \bsupsilon_{\bx}^\delta(\by)| \, d \by \right) \, d \bx\\
        &\leq C \delta^2 \int_{B_1({\bf 0})} |\bz|^4 \gamma(\bz) \, d \bz  \int_\Omega \int_\Omega \int_0^1 \1_{B_{\bx}^\delta}(\by) \gamma(\bsupsilon_{\bx}^\delta(\by))  |\bsupsilon_{\bx}^\delta(\by)|^2  |D^3 v(\bx+t(\by-\bx))|^2 \, |\det \grad \bsupsilon_{\bx}^\delta(\by)| \, d t \, d \bx \, d \by.
    \end{split}
\end{equation*}
Now we change variables in the $\bx$-integral; let $\bw = \bsupsilon_{\bx}^\delta(\by)$; then the comparability of the Jacobian to $\eta_{\delta}(\by)^{-d}$ results in the bound
\begin{equation*}
    \begin{split}
        \vnorm{III}_{L^2(\Omega)}^2 &\leq C(\kappa,\gamma) \delta^2 \int_\Omega \int_{B_1({\bf 0})} \int_0^1 |\bw|^2  \gamma(\bw) |D^3 v(\by + (1-t) \eta_{\delta}(\by) \bw)|^2 \, d t \, d \bw \, d \by \\
        &= C(\kappa,\gamma) \delta^2 \int_0^1 \int_{B_1({\bf 0})} |\bw|^2 \ga(\bw)\vnorm{ D^3 v \circ \bsPhi_{1-t,\bw} }_{L^2(\Omega)}^2 \, d \bw \, d t.
    \end{split}
\end{equation*}
Applying \eqref{eq:hetero-composition-bound} gives $\vnorm{III}_{L^2(\Omega)}^2 \leq C \delta^2 \vnorm{v}_{H^3(\Omega)}^2$. Finally, integrating \eqref{eq-hetero-uniform-pf1} and \eqref{eq-hetero-uniform-pf2} in $\bx$ gives the result for $v \in C^3(\overline{\Omega})$. Approximating a general function in $H^3(\Omega)$ with $C^3$ functions finishes the proof.

\end{proof}

\begin{lemma}
\label{lem:hetero-H1-test-consistency}
For every \(v\in H^2(\Om)\) and \(\varphi\in H^1(\Om)\),
\begin{equation}\label{eq:hetero-H1-test-consistency-H2}
 \abs{
   a_{\del,\eta}(v,\varphi)
   -(\grad v,\grad\varphi)_{L^2(\Om)}
 }
 \le C\del\norm{v}_{H^2(\Om)}
              \norm{\varphi}_{H^1(\Om)}.
\end{equation}
For every \(v\in H^3(\Om)\) and \(\varphi\in H^1(\Om)\),
\begin{equation}\label{eq:hetero-H1-test-consistency}
 \abs{
   a_{\del,\eta}(v,\varphi)
   -(\grad v,\grad\varphi)_{L^2(\Om)}
 }
 \le C\del^2\norm{v}_{H^3(\Om)}
               \norm{\varphi}_{H^1(\Om)}.
\end{equation}
\end{lemma}

\begin{proof}
We provide the proof only for \eqref{eq:hetero-H1-test-consistency}. The estimate \eqref{eq:hetero-H1-test-consistency-H2} follows from a very similar argument, and we remark on the differences at the end of the proof.

We first assume \(v,\varphi\in C^\infty(\overline\Om)\).
The substitution
\(\zb=(\yb-\xb)/\eta_{\delta}(\xb)\) in the bilinear form gives
\[
\begin{split}
 a_{\del,\eta}(v,\varphi)
 =\int_\Om\int_{B_1(\bm0)}
 \ga(\zb)
 &\frac{v(\xb+\eta_{\delta}(\xb)\zb)-v(\xb)}
       {\eta_{\delta}(\xb)}\\
 {}&\quad
 \frac{\varphi(\xb+\eta_{\delta}(\xb)\zb)-\varphi(\xb)}
       {\eta_{\delta}(\xb)}
 \,d\zb\,d\xb .
\end{split}
\]
By the fundamental theorem of calculus,
\[
 \frac{\varphi(\xb+\eta_{\delta}(\xb)\zb)-\varphi(\xb)}
      {\eta_{\delta}(\xb)}
 =
 \int_0^1
 \zb\cdot\grad\varphi(\bsPhi_{t,\zb}(\xb))\,dt .
\]
Consequently,
\begin{equation}\label{eq:hetero-consistency-identity}
\begin{split}
 a_{\del,\eta}(v,\varphi)
 =\int_\Om\int_{B_1(\bm0)}\int_0^1
 \ga(\zb)
 \frac{v(\xb+\eta_{\delta}(\xb)\zb)-v(\xb)}
      {\eta_{\delta}(\xb)}
 \zb\cdot\grad\varphi(\bsPhi_{t,\zb}(\xb))
 \,dt\,d\zb\,d\xb.
\end{split}
\end{equation}
The integral is absolutely convergent by the Cauchy--Schwarz inequality,
\eqref{eq:hetero-H1-embedding}, and \eqref{eq:hetero-composition-bound}.
We need to choose a Taylor expansion of $\frac{v(\xb+\eta_{\delta}(\xb)\zb)-v(\xb)}
      {\eta_{\delta}(\xb)}$ so that lower-order terms vanish.
To this end, for each $t \in (0,1)$, Taylor expansion of both $v(\xb+\eta_{\delta}(\xb)\zb)$ and $v(\xb)$ about $\bsPhi_{t,\zb}(\xb)=\xb+t\eta_{\delta}(\xb)\zb$ gives
\begin{align*}
 v(\xb+\eta_{\delta}(\xb)\zb)
 ={}&v(\bsPhi_{t,\zb}(\xb))
 +(1-t)\eta_{\delta}(\xb)
   \grad v(\bsPhi_{t,\zb}(\xb))\cdot\zb\\
 &+\frac{(1-t)^2\eta_{\delta}(\xb)^2}{2}
   D^2v(\bsPhi_{t,\zb}(\xb))[\zb,\zb]\\
 &+\frac{\eta_{\delta}(\xb)^3}{2}
   \int_0^{1-t}(1-t-r)^2
   D^3v(\bsPhi_{t+r,\zb}(\xb))[\zb,\zb,\zb]\,dr
\end{align*}
and
\begin{align*}
 v(\xb)
 ={}&v(\bsPhi_{t,\zb}(\xb))
 -t\eta_{\delta}(\xb)
   \grad v(\bsPhi_{t,\zb}(\xb))\cdot\zb\\
 &+\frac{t^2\eta_{\delta}(\xb)^2}{2}
   D^2v(\bsPhi_{t,\zb}(\xb))[\zb,\zb]\\
 &-\frac{\eta_{\delta}(\xb)^3}{2}
   \int_{-t}^{0}(t+r)^2
   D^3v(\bsPhi_{t+r,\zb}(\xb))[\zb,\zb,\zb]\,dr.
\end{align*}

Subtracting the two and dividing by \(\eta_{\delta}(\xb)\) yields
\begin{align*}
 &\frac{v(\xb+\eta_{\delta}(\xb)\zb)-v(\xb)}
       {\eta_{\delta}(\xb)}
 =\grad v(\bsPhi_{t,\zb}(\xb))\cdot\zb
   +\left(\frac12-t\right)
      \eta_{\delta}(\xb)
      D^2v(\bsPhi_{t,\zb}(\xb))[\zb,\zb]\\
 &\quad
 +\frac{\eta_{\delta}(\xb)^2}{2}
   \int_{-t}^{0}(t+r)^2
   D^3v(\bsPhi_{t+r,\zb}(\xb))[\zb,\zb,\zb]\,dr\\
 &\quad
 +\frac{\eta_{\delta}(\xb)^2}{2}
   \int_{0}^{1-t}(1-t-r)^2
   D^3v(\bsPhi_{t+r,\zb}(\xb))[\zb,\zb,\zb]\,dr .
\end{align*}
After multiplying by $\gamma(\bz) \grad\varphi(\bsPhi_{t,\zb}(\xb))\cdot\zb$ and integrating over $t \in (0,1)$, $\bz \in B_1({\bf 0})$ and $\bx \in \Omega$, \eqref{eq:hetero-composition-bound} and the Cauchy-Schwarz inequality imply that the last two expressions in the sum contribute the following majorant to the estimate of \eqref{eq:hetero-consistency-identity}:
\[
 C\del^2
 \left(\int_{B_1(\bm0)}\ga(\zb)\abs{\zb}^4\,d\zb\right)
 \norm{v}_{H^3(\Om)}\norm{\varphi}_{H^1(\Om)}.
\]
Therefore \eqref{eq:hetero-H1-test-consistency} will be established for smooth functions if we show that
\begin{equation}\label{eq:hetero-H1-test-consistency-Pf1}
    \begin{split}
     I := \left| \int_\Omega \int_{B_1({\bf 0})} \int_0^1 \ga(\zb) (\zb \otimes \zb) 
     [ \grad v(\bsPhi_{t,\zb}(\xb)), \grad\varphi(\bsPhi_{t,\zb}(\xb))] \,dt\,d\zb\,d\xb 
    - (\grad v, \grad \varphi)_{L^2(\Omega)} \right| 
    \end{split}
\end{equation}
and
\begin{equation*}
    II := \left| \int_\Omega \int_{B_1({\bf 0})} \int_0^1 \ga(\zb)\left(\frac12-t\right)
      \eta_{\delta}(\xb)
      D^2v(\bsPhi_{t,\zb}(\xb))[\zb,\zb]
\grad\varphi(\bsPhi_{t,\zb}(\xb))\cdot \zb
 \,dt\,d\zb\,d\xb \right| 
\end{equation*}
satisfy
\begin{equation}\label{eq:hetero-H1-test-consistency-Pf3}
    I \leq C \delta^2 \| v \|_{H^1(\Omega)} \| \varphi \|_{H^1(\Omega)} \qquad \text{ and } \qquad II \leq C \delta^2 \| v \|_{H^2(\Omega)} \| \varphi \|_{H^1(\Omega)}.
\end{equation}

To this end, we use the normalization condition on $\gamma$ and the coordinate change $\bx = \bsPsi_{t,\bz}(\by)$ to write
\begin{equation*}
    \begin{split}
        (\grad v, \grad \varphi)_{L^2(\Omega)} &= \int_\Omega \int_{B_1({\bf0})} \int_0^1 \gamma(\bz) (\bz \otimes \bz) [ \grad v(\by),\grad \varphi(\by) ] \,dt\,d\bz\,d\by \\
        &= \int_\Omega \int_{B_1({\bf0})} \int_0^1 \gamma(\bz) (\bz \otimes \bz) [ \grad v(\bsPhi_{t,\bz}(\bx)),\grad \varphi(\bsPhi_{t,\bz}(\bx)) ] J_{t,\bz}(\bx) \,dt\,d\bz\,d\bx, 
    \end{split}
\end{equation*}
where the identities of Lemma \ref{lem:hetero-coordinate-map} were additionally used.
Next, since $(\bz \otimes \bz) [ \grad v(\by),\grad \varphi(\by) ] (\grad \eta_{\delta}(\by) \cdot \bz)$ is an odd function of $\bz$, a similar calculation shows that
\begin{equation*}
    0 =  \int_\Omega \int_{B_1({\bf0})} \int_0^1 \gamma(\bz) (\bz \otimes \bz) [ \grad v(\bsPhi_{t,\bz}(\bx)),\grad \varphi(\bsPhi_{t,\bz}(\bx)) ] t (\grad \eta_{\delta}(\bsPhi_{t,\bz}(\bx)) \cdot \bz) J_{t,\bz}(\bx) \,dt\,d\bz\,d\bx.
\end{equation*}
Placing these two formulas into \eqref{eq:hetero-H1-test-consistency-Pf1}, 
\begin{equation*}
    I = \left| \int_\Omega \int_{B_1({\bf 0})} \int_0^1 \ga(\zb) (\zb \otimes \zb) [ \grad v(\bsPhi_{t,\zb}(\xb)), \grad\varphi(\bsPhi_{t,\zb}(\xb))] \mathcal{K}(\bx,\bz,t,\delta) \,dt\,d\zb\,d\xb \right|,
\end{equation*}
where $\mathcal{K}(\bx,\bz,t,\delta) := 1 - J_{t,\bz}(\bx) + t J_{t,\bz}(\bx) \grad \eta_{\delta}(\bsPhi_{t,\bz}(\bx)) \cdot \bz$. Thus, the estimate \eqref{eq:hetero-H1-test-consistency-Pf3} for $I$ will follow from the assumptions on $\gamma$, the Cauchy-Schwarz inequality, and \eqref{eq:hetero-composition-bound} if we show that
\begin{equation}\label{eq:hetero-H1-test-consistency-Pf4}
    |\mathcal{K}(\bx,\bz,t,\delta)| \leq C \delta^2 \qquad \text{ on } \Omega \times B_1({\bf 0}) \times [0,1].
\end{equation}
To see this, use that $|\grad \eta_{\delta}| \leq \kappa \delta$ and that $\grad \eta_{\delta}$ is Lipschitz with Lipschitz constant $\kappa \delta$:
\begin{equation*}
    \begin{split}
        |\mathcal{K}(\bx,\bz,t,\delta)| &= |t \bz \cdot \big( \grad \eta_{\delta}(\bsPhi_{t,\bz}(\bx)) - \grad \eta_{\delta}(\bx) \big) + t^2 (\bz \otimes \bz)[\grad \eta_{\delta}(\bx),\grad \eta_{\delta}(\bsPhi_{t,\bz}(\bx))]| \\
        &\leq \kappa \delta |\bsPhi_{t,\bz}(\bx) - \bx| + \kappa^2 \delta^2 \\
        &= t \kappa \delta \eta_{\delta}(\bx) |\bz| + \kappa^2 \delta^2 \leq C \delta^2.
    \end{split}
\end{equation*}
Therefore \eqref{eq:hetero-H1-test-consistency-Pf4}, and hence \eqref{eq:hetero-H1-test-consistency-Pf3} for $I$, holds.

It remains to show the estimate \eqref{eq:hetero-H1-test-consistency-Pf3} for $II$. We use a similar argument of subtracting an odd function and estimating the difference. Since
\begin{equation*}
    0 = \int_\Omega \int_{B_1({\bf 0})} \int_0^1 \ga(\zb)\left(\frac12-t\right)
      \eta_{\delta}(\bsPhi_{t,\zb}(\xb))
      D^2v(\bsPhi_{t,\zb}(\xb))[\zb,\zb]
(\grad\varphi(\bsPhi_{t,\zb}(\xb))\cdot \zb) J_{t,\bz}(\bx)
 \,dt\,d\zb\,d\xb,
\end{equation*}
we can subtract this from the integral defining $II$, and so \eqref{eq:hetero-H1-test-consistency-Pf3} for $II$ holds if 
\begin{equation*}
    |\eta_{\delta}(\bsPhi_{t,\zb}(\xb)) J_{t,\bz}(\bx) - \eta_{\delta}(\xb)| \leq C\delta^2 \qquad \text{ for } (\bx,\bz,t) \in \Omega \times B_1({\bf 0}) \times [0,1].
\end{equation*}
The estimate is similar to that of $\cK$, and uses the same properties of $\eta_{\delta}$:
\begin{equation*}
    \begin{split}
        |\eta_{\delta}(\bsPhi_{t,\zb}(\xb)) J_{t,\bz}(\bx) - \eta_{\delta}(\xb)| &= |\eta_{\delta}(\bsPhi_{t,\zb}(\xb)) - \eta_{\delta}(\bx) + t  \eta_{\delta}(\bsPhi_{t,\zb}(\xb)) \grad \eta_{\delta}(\xb) \cdot \bz| \\
        &\leq |\eta_{\delta}(\bsPhi_{t,\zb}(\xb)) - \eta_{\delta}(\bx)| + \kappa^2 \delta^2 \\
        &\leq \kappa \delta |\bsPhi_{t,\zb}(\xb) -\bx| + \kappa^2 \delta^2 \leq C \delta^2.
    \end{split}
\end{equation*}
Thus \eqref{eq:hetero-H1-test-consistency-Pf3} for $II$ is established, and so \eqref{eq:hetero-H1-test-consistency} holds for $v$ and $\varphi$ in $C^\infty(\overline{\Omega})$.

Finally, to see that \eqref{eq:hetero-H1-test-consistency} holds for $v \in H^3(\Omega)$ and $\varphi \in H^1(\Omega)$, note that both bilinear forms are continuous on \(H^1(\Om)\times H^1(\Om)\) by
\eqref{eq:hetero-H1-embedding}. Then the result follows by density of smooth functions in $H^k(\Omega)$ for $k = 1,2,3$.

The proof of \eqref{eq:hetero-H1-test-consistency-H2} follows a similar thread of reasoning; the main difference is that we use Taylor expansion up to the second order so that the remainder term is controlled by \(O(\delta)\) involving \(\norm{v}_{H^2(\Om)}\). The estimate of the term $I$ remains the same. The detailed proof is thus omitted.

\end{proof}

\begin{corollary}
\label{cor:hetero-H2-H2-consistency}
Under the assumptions of \Cref{lem:hetero-H1-test-consistency},
there is a constant $C>0$, independent of $\delta$, such that
for every $v,\varphi\in H^2(\Omega)$,
\begin{equation}
\label{eq:hetero-H2-H2-consistency}
  \bigl|a_{\delta,\eta}(v,\varphi)-(\nabla v,\nabla \varphi)_{L^2(\Omega)}\bigr|
  \leq C\delta^2\|v\|_{H^2(\Omega)}\|\varphi\|_{H^2(\Omega)}.
\end{equation}
\end{corollary}
\begin{proof}
The corollary can be proved in the same way as \eqref{eq:hetero-H1-test-consistency}; Taylor expanding both $v$ and $\varphi$ to second order presents an additional power of $\delta$. We omit the details.
\end{proof}

The next proposition uses the estimates just proved to control parts of the residual error.
\begin{proposition}
\label{prop:hetero-full-residual}
For every \(v\in H^2(\Om)\) and \(\varphi\in\cS_{\del}(\Om)\),
\begin{equation}\label{eq:hetero-full-residual-H2}
 \abs{(\Delta v-\cL_{\del,\eta}v,\varphi)_{L^2(\Om)}}
 \le C\del\norm{v}_{H^2(\Om)}
              \norm{\varphi}_{\cS_{\del}(\Om)}.
\end{equation}
For every \(v\in H^3(\Om)\) and \(\varphi\in\cS_{\del}(\Om)\),
\begin{equation}\label{eq:hetero-full-residual}
 \abs{(\Delta v-\cL_{\del,\eta}v,\varphi)_{L^2(\Om)}}
 \le C\del^2\norm{v}_{H^3(\Om)}
               \norm{\varphi}_{\cS_{\del}(\Om)}.
\end{equation}
Consequently, for every integer \(m\ge1\), every
\(v\in H^{m+1}(\Om)\), and every
\(\varphi\in\cS_{\del}(\Om)\),
\begin{equation}\label{eq:hetero-full-residual-regularity}
 \abs{(\Delta v-\cL_{\del,\eta}v,\varphi)_{L^2(\Om)}}
 \le C\del^{\min\{2,m\}}\norm{v}_{H^{m+1}(\Om)}
               \norm{\varphi}_{\cS_{\del}(\Om)}.
\end{equation}
The same estimates hold on each constrained test space.
\end{proposition}

\begin{proof}
\emph{Step 1: the boundary localized convolution term.}
We first estimate the part tested against \(K_\del \varphi\).  The local and
nonlocal Green identities give
\begin{align*}
 &(\Delta v-\cL_{\del,\eta}v,K_\del \varphi)_{L^2(\Om)}
 \\
 &\quad=
 a_{\del,\eta}(v,K_\del \varphi)
 -(\grad v,\grad K_\del \varphi)_{L^2(\Om)}
 +(1-Q_\del)
   \langle\p_{\nb}v,\Tr(K_\del \varphi)\rangle_{\p\Om}.
\end{align*}
For \(v\in H^2(\Om)\), use
\eqref{eq:hetero-H1-test-consistency-H2}; for
\(v\in H^3(\Om)\), use
\eqref{eq:hetero-H1-test-consistency}.  The trace theorem,
\eqref{eq:hetero-Q-expansion}, and
\eqref{eq:hetero-boundary-convolution-H1} then give, respectively,
\[
 \abs{(\Delta v-\cL_{\del,\eta}v,K_\del \varphi)_{L^2(\Om)}}
 \le C\del\norm{v}_{H^2(\Om)}
              \norm{\varphi}_{\cS_{\del}(\Om)}
\]
and
\[
 \abs{(\Delta v-\cL_{\del,\eta}v,K_\del \varphi)_{L^2(\Om)}}
 \le C\del^2\norm{v}_{H^3(\Om)}
               \norm{\varphi}_{\cS_{\del}(\Om)}.
\]

\emph{Step 2: the remaining test-function term.}
We next estimate the part paired with \(\varphi-K_\del \varphi\).  For \(v\in H^2(\Om)\),
\eqref{eq:hetero-operator-H2} and
\eqref{eq:hetero-boundary-convolution-L2} give
\[
 \abs{(\Delta v-\cL_{\del,\eta}v,\varphi-K_\del \varphi)_{L^2(\Om)}}
 \le C\del\norm{v}_{H^2(\Om)}
              \norm{\varphi}_{\cS_{\del}(\Om)}.
\]
For \(v\in H^3(\Om)\),
\eqref{eq:hetero-strong-first-order} supplies one power of \(\del\), and
\eqref{eq:hetero-boundary-convolution-L2} supplies the other:
\[
 \abs{(\Delta v-\cL_{\del,\eta}v,\varphi-K_\del \varphi)_{L^2(\Om)}}
 \le C\del^2\norm{v}_{H^3(\Om)}
               \norm{\varphi}_{\cS_{\del}(\Om)}.
\]

\emph{Step 3: conclusion.}
We finally add the two parts of
\[
 (\Delta v-\cL_{\del,\eta}v,\varphi)_{L^2(\Om)}
 =
 (\Delta v-\cL_{\del,\eta}v,K_\del \varphi)_{L^2(\Om)}
 +
 (\Delta v-\cL_{\del,\eta}v,\varphi-K_\del \varphi)_{L^2(\Om)}.
\]
This proves \eqref{eq:hetero-full-residual-H2} and
\eqref{eq:hetero-full-residual}.  When \(m=1\),
\eqref{eq:hetero-full-residual-regularity} is the first of these estimates.
When \(m\ge2\), \eqref{eq:hetero-full-residual-regularity} follows from the second estimate and the embedding
\(H^{m+1}(\Om)\hookrightarrow H^3(\Om)\).  Restriction gives the same
estimates on each constrained test space.
\end{proof}

\subsection{Proof of the continuous convergence theorem}
\label{subsec:hetero-continuous-proof}

\begin{proof}[Proof of \Cref{thm:hetero-rates-continuous}]
We first identify the residual in each problem.  Applying the local and
nonlocal Green identities to \(U_0\), and using the corresponding boundary
condition, gives
\begin{equation}\label{eq:hetero-residual-cases}
    \Rdel(v)=
    \begin{cases}
       -(\Delta U_0-\cL_{\del,\eta}U_0,v)_{L^2(\Om)},
          &\text{Dirichlet},\\
       (-(\Delta U_0-\cL_{\del,\eta}U_0)+(Q_\del-1)f,v)_{L^2(\Om)},
          &\text{Neumann},\\
       -(\Delta U_0-\cL_{\del,\eta}U_0,v)_{L^2(\Om)},
          &\text{Robin}.
    \end{cases}
\end{equation}
For example, in the Neumann case
\[
 \begin{split}
 \Rdel(v)
 &=Q_\del(f,v)_{L^2(\Om)}
      +Q_\del\langle g,\Tr v\rangle_{\p\Om}
      -a_{\del,\eta}(U_0,v)\\
 &=(Q_\del f+\cL_{\del,\eta} U_0,v)_{L^2(\Om)}\\
 &=(-(\Delta U_0-\cL_{\del,\eta}U_0)
       +(Q_\del-1)f,v)_{L^2(\Om)}.
 \end{split}
\]
For Robin data, the two boundary terms cancel because
\(\p_{\nb} U_0=(g-b\Tr U_0)/\mu\).  In the Dirichlet case \(\Tr v=0\), and
\(U_0=u_0+G\), so the same calculation produces the first line of
\eqref{eq:hetero-residual-cases}.

In all three cases, the difference of operators in \eqref{eq:hetero-residual-cases} is controlled by
\eqref{eq:hetero-full-residual-regularity}, and additionally in the Neumann case, the
remaining term can be seen to have the same order; using that
\(f=-\Delta U_0\), \(\del_0<1\), and 
\eqref{eq:hetero-Q-expansion}
\begin{equation}\label{eq:NeumannCorrection}
 \begin{split}
 \abs{(Q_\del-1)(f,v)_{L^2(\Om)}}
 &\le C\del^2\norm{f}_{L^2(\Om)}\norm{v}_{L^2(\Om)}\\
 &\le C\del^{\min\{2,m\}}\norm{U_0}_{H^{m+1}(\Om)}
       \norm{v}_{\cS_{\del}(\Om)}.
 \end{split}
\end{equation}
Thus in all three cases
\begin{equation}\label{eq:hetero-residual-dual-bound}
    \abs{\Rdel(v)}
    \le C\del^{\min\{2,m\}}
       \norm{U_0}_{H^{m+1}(\Om)}
       \norm{v}_{\cS_{\del}(\Om)},
    \qquad v\in\cV_\del.
\end{equation}
Coercivity and \Cref{thm:abstract-energy-rate} imply
\eqref{eq:hetero-energy-rate}.

We next prove \eqref{eq:hetero-L2-rate}.
Since $m\ge1$, we have $U_0\in H^2(\Om)$.
The residual estimate \eqref{eq:hetero-residual-dual-bound}
with $m=1$ gives
\[
  |\mathcal{R}_\delta(v)|
  \leq C\delta\|U_0\|_{H^2(\Omega)}
              \|v\|_{\mathcal{S}_\delta(\Omega)}.
\]
For $\phi\in L^2(\Omega)$, let $z_\delta=z_\delta[\phi]$
and $z_0=z_0[\phi]$ solve the nonlocal and local dual problems
\[
  a_\delta(v,z_\delta)=(v,\phi)_{L^2(\Omega)},
  \quad v\in\mathcal{U}_\delta,
  \qquad
  a_0(v,z_0)=(v,\phi)_{L^2(\Omega)},
  \quad v\in\mathcal{U}_0.
\]
In the Neumann case, replace $\phi$ by $\phi-(\phi)_\Omega$. By \eqref{eq:hetero-local-dual-H2-regularity},
$\|z_0\|_{H^2(\Omega)}\leq C\|\phi\|_{L^2(\Omega)}$. The nonlocal Green's identity and the homogeneous local boundary condition for $z_0$ give, for
$v\in\mathcal{V}_\delta$, \[
  a_\delta(z_\delta-z_0,v) = (\phi +\mathcal{L}_{\delta,\eta}z_0, v )_{L^2(\Omega)}
  =-(\Delta z_0-\mathcal{L}_{\delta,\eta}z_0,v)_{L^2(\Omega)}.
\]
\Cref{prop:hetero-full-residual} and coercivity therefore
give
\begin{equation}
    \label{eq:hetero-dual-energy-H2}
  \|z_\delta-z_0\|_{\mathcal{S}_\delta(\Omega)}
  \leq C\delta\|z_0\|_{H^2(\Omega)}
  \leq C\delta\|\phi\|_{L^2(\Omega)}.
\end{equation}
We now estimate $\mathcal{R}_\delta(z_0)$. Set $B_\delta(v,\varphi) :=a_{\delta,\eta}(v,\varphi)-(\nabla v,\nabla \varphi)_{L^2(\Omega)}$.
Then we have
\[
  \mathcal{R}_\delta(z_0)
  =-B_\delta(U_0,z_0)+
  \begin{cases}
    0, &\text{Dirichlet},\\
    (Q_\delta-1)(\nabla U_0,\nabla z_0)_{L^2(\Omega)},
      &\text{Neumann},\\
    (Q_\delta-1)
      \langle\partial_n U_0,\operatorname{Tr}z_0\rangle_{\partial\Omega},
      &\text{Robin}.
  \end{cases}
\]
For Neumann data, the correction follows from
$(f,z_0)_{L^2(\Omega)}+\langle g,\operatorname{Tr}z_0\rangle_{\partial\Omega}
=(\nabla U_0,\nabla z_0)_{L^2(\Omega)}$.
For Robin data, it follows from
$(g-b\operatorname{Tr}U_0)/\mu=\partial_n U_0$. \Cref{cor:hetero-H2-H2-consistency},
$|Q_\delta-1|\leq C\delta^2$, and the trace estimates now yield
\begin{equation}
    \label{eq:hetero-local-dual-residual-H2} |\mathcal{R}_\delta(z_0)|
  \leq C\delta^2\|U_0\|_{H^2(\Omega)}\|z_0\|_{H^2(\Omega)}
  \leq C\delta^2\|U_0\|_{H^2(\Omega)}\|\phi\|_{L^2(\Omega)}.
\end{equation}
Consequently,
\[
\begin{aligned}
  |\mathcal{R}_\delta(z_\delta)|
  &\leq |\mathcal{R}_\delta(z_0)|
       +|\mathcal{R}_\delta(z_\delta-z_0)|\leq C\delta^2\|U_0\|_{H^2(\Omega)}\|\phi\|_{L^2(\Omega)}
+ C\delta\|U_0\|_{H^2(\Omega)}
          \|z_\delta-z_0\|_{\mathcal{S}_\delta(\Omega)}\\
  &\leq C\delta^2\|U_0\|_{H^2(\Omega)}\|\phi\|_{L^2(\Omega)}.
\end{aligned}
\]
\Cref{thm:abstract-L2-rate} now implies \eqref{eq:hetero-L2-rate}.

It remains to estimate the approximate normal derivative.  Let
\(\varphi\in H^{1/2}(\p\Om)\).  The extension operator satisfies
\(\norm{E\varphi}_{H^1(\Om)}\le C\norm{\varphi}_{H^{1/2}(\p\Om)}\).
By \eqref{eq:NormalDerivDef}, \eqref{eq:GreensIdentity}, and the local Green's identity, we have
\begin{equation}\label{eq:NormalDerivativeEstimate}
 Q_\del\left\langle
    (\p_{\nb} U_\del)_\del-\p_{\nb} U_0,\varphi
 \right\rangle
 =a_{\del,\eta}(\edel,E\varphi)
   -(f_\del+\cL_{\del,\eta}U_0,E\varphi)_{L^2(\Om)}.
\end{equation}
Here
\[
 f_\del+\cL_{\del,\eta}U_0
 =\begin{cases}
   -(\Delta U_0-\cL_{\del,\eta}U_0),
      &\text{Dirichlet or Robin},\\
   -(\Delta U_0-\cL_{\del,\eta}U_0)+(Q_\del-1)f,
      &\text{Neumann}.
 \end{cases}
\]
The first term of \eqref{eq:NormalDerivativeEstimate} is bounded by
\(C\norm{\edel}_{\cS_{\del}(\Om)}
\norm{E\varphi}_{H^1(\Om)}\).  The second term of \eqref{eq:NormalDerivativeEstimate} is bounded
by \eqref{eq:hetero-full-residual-regularity} and \eqref{eq:NeumannCorrection}.  Since \(Q_\del\) is uniformly bounded away
from zero, taking the supremum over \(\varphi\) proves
\eqref{eq:hetero-normal-rate}.
\end{proof}

\subsection{Proof of the discrete estimates}
\label{subsec:hetero-discrete-proof}

\begin{proof}[Proof of \Cref{thm:hetero-rates-discrete}]
Set \(r:=\min\{k,m\}\) and \(s:=\min\{2,m\}\). To establish \eqref{eq:hetero-energy-rate-discrete}, we have by
\eqref{eq:hetero-H1-embedding} and the quasi-interpolation estimate,
\begin{equation}\label{eq:hetero-primal-best-approximation}
 \inf_{v_h\in\cU_{\del,h}}
    \norm{u_0-v_h}_{\cS_{\del}(\Om)}
 \le C h^r\norm{u_0}_{H^{r+1}(\Om)}
 \le Ch^r\norm{u_0}_{H^{m+1}(\Om)}.
\end{equation}
The interpolant is chosen to preserve zero trace (in the Dirichlet case) or zero mean (in the Neumann case).
Combining
\eqref{eq:hetero-residual-dual-bound} and
\eqref{eq:hetero-primal-best-approximation} with the discrete
energy estimate of \Cref{thm:abstract-energy-discrete} proves
\eqref{eq:hetero-energy-rate-discrete}.

For the $L^2$ estimate, let $z_\delta=z_\delta[\phi]$ and
$z_0=z_0[\phi]$ be the dual solutions used above, and let
$z_{\delta,h}=z_{\delta,h}[\phi]\in\mathcal{V}_{\delta,h}$ satisfy
\[
  a_\delta(v_h,z_{\delta,h})
  =(v_h,\phi)_{L^2(\Omega)},
  \qquad v_h\in\mathcal{U}_{\delta,h}.
\] 
As before, we replace $\phi$ by $\phi-(\phi)_\Omega$ in the Neumann case.
The dual energy estimate \eqref{eq:hetero-dual-energy-H2}, finite element approximation of $z_0\in H^2(\Omega)$  and \eqref{eq:hetero-H1-embedding} give
\[
  \inf_{v_h\in\mathcal{V}_{\delta,h}}
    \|z_\delta-v_h\|_{\mathcal{S}_\delta(\Omega)}
  \leq C(\delta+h)\|\phi\|_{L^2(\Omega)}.
\]
Thus the function \(\rho\) in
\Cref{assu:abstract-dual-discrete} can be taken to be \(C(\del+h)\). Using this, the triangle inequality and \eqref{eq:hetero-dual-energy-H2}, we have
\begin{equation}
    \label{eq:hetero-discrete-dual-energy-H2}
  \|z_{\delta,h}-z_0\|_{\mathcal{S}_\delta(\Omega)}
  \leq C(\delta+h)\|\phi\|_{L^2(\Omega)}.
\end{equation}

If $m\geq2$, 
discrete dual stability 
\(\norm{z_{\del,h}[\phi]}_{\cS_{\del}(\Om)}
\le C\norm{\phi}_{L^2(\Om)}\) and \eqref{eq:hetero-residual-dual-bound}  yield
\begin{equation*}
    \abs{\Rdel(z_{\del,h}[\phi])}
    \le C\del^s\norm{U_0}_{H^{m+1}(\Om)}
       \norm{\phi}_{L^2(\Om)}.
\end{equation*}
Thus
\(\beta(\del,h)=C\del^s\norm{U_0}_{H^{m+1}(\Om)}\).  \Cref{thm:abstract-L2-rate-discrete} and
\eqref{eq:hetero-primal-best-approximation} now give 
\begin{equation}\label{eq:hetero-discrete-L2-preliminary}
 \norm{U_{\del,h}-U_0}_{L^2(\Om)}
 \le C\bigl(\del^s+(\del+h)h^r\bigr)
       \bigl(\norm{U_0}_{H^{m+1}(\Om)}
             +\norm{u_0}_{H^{m+1}(\Om)}\bigr).
\end{equation}
Since $r\geq1$ and $0<h\leq1$, 
Young's inequality gives
\[
    \del h^r
    \le \frac12\del^2+\frac12h^{2r}
    \le \frac12\del^2+\frac12h^{r+1}.
\]
This proves \eqref{eq:hetero-L2-rate-discrete} for $m\geq2$.

If $m=1$, then $r=1$.
Split the residual against the discrete dual solution as
\[
  \mathcal{R}_\delta(z_{\delta,h})
  =\mathcal{R}_\delta(z_0)
   +\mathcal{R}_\delta(z_{\delta,h}-z_0).
\]
Use \eqref{eq:hetero-local-dual-residual-H2} for the first term
and \eqref{eq:hetero-residual-dual-bound} at $m=1$ for the second.
Combining these estimates with
\eqref{eq:hetero-discrete-dual-energy-H2} gives
\[
\begin{aligned}
  |\mathcal{R}_\delta(z_{\delta,h})|
  &\leq C\delta^2\|U_0\|_{H^2(\Omega)}\|\phi\|_{L^2(\Omega)}+ C\delta\|U_0\|_{H^2(\Omega)}
      \|z_{\delta,h}-z_0\|_{\mathcal{S}_\delta(\Omega)}\\
  &\leq C(\delta^2+\delta h)
      \|U_0\|_{H^2(\Omega)}\|\phi\|_{L^2(\Omega)}.
\end{aligned}
\]
Thus we may take
$\beta(\delta,h)=C(\delta^2+\delta h)\|U_0\|_{H^2(\Omega)}$.
Theorem~\ref{thm:abstract-L2-rate-discrete} and
\eqref{eq:hetero-primal-best-approximation} now give
\[
  \|U_{\delta,h}-U_0\|_{L^2(\Omega)}
  \leq C(\delta^2+\delta h+h^2)\bigl(\|U_0\|_{H^{2}(\Omega)}
      +\|u_0\|_{H^{2}(\Omega)}\bigr)
  \leq C(\delta^2+h^2)\bigl(\|U_0\|_{H^{2}(\Omega)}
      +\|u_0\|_{H^{2}(\Omega)}\bigr).
\]
This is \eqref{eq:hetero-L2-rate-discrete} for $m=1$ and completes
the proof.
\end{proof}

\section{Numerical experiments}\label{sec:Numerical}

\subsection{Experiments for the type (I) nonlocal Laplacian}
\label{subsec:fullspace-numerics}
We consider \(\Om=(0,1)\) with zero exterior values and use the
manufactured local solution
\[
    u_0(x)=x^m(1-x)^m
    \left[1+0.4\sin(2\pi x)+0.2\cos(3\pi x)\right],
    \qquad f(x)=-u_0''(x),\quad x\in(0,1).
\]
For each \(m\), we fix \(f(x)\) and compare \(u_{\del,h}\) with
\(E_0u_0\), which equals \(u_0\) in \(\Om\) and zero in the exterior
collar.  We first examine the modeling error as \(\del\) decreases,
then study the discrete error along paths \(\del=h^\beta\).

For the general estimates, we use the fractional kernel
\[
    \ga_{\del}(s)
    =(1-\alpha)\del^{2\alpha-2}|s|^{-1-2\alpha}
      \mathbf{1}_{\{0<|s|\le\del\}},\qquad \alpha=0.95,
\]
where \(\alpha\) is the fractional order.  This kernel satisfies
\(\int_{-\del}^{\del}s^2\ga_{\del}(s)\,\dd s=1\).
We solve \eqref{eq:fullSpace:Poisson:weak} using continuous Galerkin
(CG) spaces of degree \(k\) on uniform meshes with \(h=1/N\), imposing
zero values at both endpoints.
For each kernel, we use its own full nonlocal norm,
\[
    \|v\|_{\cS_{\del}(\Om)}^2
    =\|v\|_{L^2(\Om)}^2+
    \int_{\Omh}\int_{\Omh}
       \ga_{\del}(x-y)\bigl(v(x)-v(y)\bigr)^2\,\dd y\,\dd x.
\]
All rates are computed from consecutive refinements of \(\del\) or \(h\).

\subsubsection{Convergence of the continuous problem}
\Cref{thm:LapI-rates-continuous} gives
\[
    \|u_{\del}-E_0u_0\|_{\cS_{\del}(\Om)}
    =O\!\left(\del^{\min\{2,m-1/2\}}\right),\qquad
    \|u_{\del}-u_0\|_{L^2(\Om)}
    =O\!\left(\del^{\min\{2,m\}}\right).
\]
We take \(k=3\), \(m=1,2,3\), and
\(\del=0.2,0.1,0.05,0.025\).
For each fixed \(\del\), we refine \(N=128,256,512\), adding
\(N=1024\) for \(\del=0.025\).  The final spatial refinement changes
both errors by less than \(0.46\%\).
\Cref{tab:tiandu-continuous-rates-revised} reports the rates between
\(\del=0.05\) and \(0.025\), using the errors on the finest meshes.
The observed rates agree with the predicted dependence on \(m\) in
both norms.

\begin{table}[H]
\centering
\caption{Modeling rates for the fractional kernel with \(\alpha=0.95\).}
\label{tab:tiandu-continuous-rates-revised}
\begin{tabular}{c|cc|cc}
\hline
\(m\) & \multicolumn{2}{c|}{Full norm} &
\multicolumn{2}{c}{\(L^2(\Om)\)}\\
\cline{2-5}
& theory & observed & theory & observed\\
\hline
1 & 0.500 & 0.508 & 1.000 & 1.005\\
2 & 1.500 & 1.583 & 2.000 & 1.983\\
3 & 2.000 & 1.974 & 2.000 & 1.999\\
\hline
\end{tabular}
\end{table}

\subsubsection{Convergence of the discrete problem}\label{subsec:ti-discrete-numerics}
Let \(r=\min\{m,k\}\) and set \(\del=h^\beta\).
\Cref{thm:LapI-rates-discrete} gives
\begin{align}
    \|u_{\del,h}-E_0u_0\|_{\cS_{\del}(\Om)}
    &=O\!\left(h^{\beta\min\{2,m-1/2\}}+h^r\right),\notag\\
    \|u_{\del,h}-u_0\|_{L^2(\Om)}
    &=O\!\left(h^{\beta\min\{2,m\}}+h^{r+\beta/2}+h^{r+1}\right).
    \label{eq:ti-L2-bound-alt}
\end{align}

We organize the experiments according to whether the mixed term
$h^{r+\beta/2}$ affects the predicted $L^2$ convergence rate.

\medskip
\noindent\textbf{Case 1: Estimates without a mixed term.} By \Cref{rem:LapI-L2-separated}, the mixed term can be absorbed when $m=1$ or $r\geq3$.
For $\delta=h^\beta$, it can also be absorbed when $\beta\geq2$, since $h^{r+\beta/2}\leq h^{r+1}$ for $0<h\leq1$.
For $m\geq2$, the same holds when $\beta\leq2r/3$, since then
$h^{r+\beta/2}\leq h^{2\beta}$.
In all these cases, the $L^2$ estimate reduces to
\[
  \|u_{\delta,h}-u_0\|_{L^2(\Omega)}
  =O\left(h^{\beta\min\{2,m\}}+h^{r+1}\right).
\]
We first take \(m=3\), \(\del=h^2\), and \(N=16,32,64,128\),
using CG spaces of degrees \(k=1,2,3\). The choice $\beta=2$ ensures that the mixed term can be absorbed. 
\Cref{tab:tiandu-degree-rates-revised} shows that the observed
rates are close to $k$ in the full norm and $k+1$ in $L^2(\Omega)$,
as predicted.
The rates are computed from $N=64$ and $N=128$.

\begin{table}[H]
\centering
\caption{CG rates for the fractional kernel with \(\alpha=0.95\),
\(m=3\), and \(\del=h^2\).}
\label{tab:tiandu-degree-rates-revised}
\begin{tabular}{c|cc|cc}
\hline
\(k\) & \multicolumn{2}{c|}{Full norm} &
\multicolumn{2}{c}{\(L^2(\Om)\)}\\
\cline{2-5}
& theory & observed & theory & observed\\
\hline
1 & 1.000 & 0.999 & 2.000 & 1.998\\
2 & 2.000 & 1.998 & 3.000 & 2.998\\
3 & 3.000 & 2.996 & 4.000 & 3.995\\
\hline
\end{tabular}
\end{table}

We next consider \((m,k,\beta)=(1,2,0.75)\) and \((3,3,1.25)\).
Both choices fall under \Cref{rem:LapI-L2-separated}, where the
mixed term in the general \(L^2\) estimate can be absorbed.
The observed \(L^2\) rates in
\Cref{tab:tiandu-complete-rates-revised} agree with the resulting
predictions.  The full norm rate for \(m=1\) is \(0.361\), close to
the predicted order \(0.375\).  For \(m=k=3\), the full norm rate
decreases from \(2.696\) at \(N=256\) to \(2.625\) at \(N=512\),
approaching the predicted order \(2.5\).

\begin{table}[H]
\centering
\caption{Additional CG rates for the fractional kernel.
Rates use \(N=256,512\) for \(m=1\) and \(N=128,256\) for \(m=3\).}
\label{tab:tiandu-complete-rates-revised}
\begin{tabular}{ccc|cc|cc}
\hline
\(m\) & \(k\) & \(\beta\) & \multicolumn{2}{c|}{Full norm} &
\multicolumn{2}{c}{\(L^2(\Om)\)}\\
\cline{4-7}
 & & & theory & observed & theory & observed\\
\hline
1 & 2 & 0.75 & 0.375 & 0.361 & 0.750 & 0.749\\
3 & 3 & 1.25 & 2.500 & 2.696 & 2.500 & 2.500\\
\hline
\end{tabular}
\end{table}

\medskip
\noindent\textbf{Case 2: Estimates with a mixed term.} The remaining cases are
\[
  m\geq2,\qquad r\in\{1,2\},\qquad \frac{2r}{3}<\beta<2.
\]
Here,
\[
  r+\frac{\beta}{2}<\min\{2\beta,r+1\},
\]
so the mixed term dominates the total error. 

We take \(m=2\) and consider \((k,\beta)=(1,0.75)\) with
\(N=32,64,128,256,512,1024\), and \((k,\beta)=(2,1.5)\) with
\(N=32,64,128,256\).
\Cref{tab:tiandu-beta-rates-revised} shows that the full norm
errors converge at approximately the predicted orders, while the
$L^2(\Omega)$ errors converge faster than predicted.
For \(k=1\), the last three consecutive \(L^2\) rates are
\(1.756,1.612,1.524\), compared with the predicted value \(1.375\).
For \(k=2\), the observed \(L^2\) rate is $3.000$, compared with the
predicted value \(2.75\).

\begin{table}[H]
\centering
\caption{CG rates for the fractional kernel with \(m=2\).
Rates use \(N=512,1024\) for \(k=1\) and \(N=128,256\) for \(k=2\).}
\label{tab:tiandu-beta-rates-revised}
\begin{tabular}{ccc|cc|cc}
\hline
\(m\) & \(k\) & \(\beta\) & \multicolumn{2}{c|}{Full norm} &
\multicolumn{2}{c}{\(L^2(\Om)\)}\\
\cline{4-7}
 & & & theory & observed & theory & observed\\
\hline
2 & 1 & 0.75 & 1.000 & 1.013 & 1.375 & 1.524\\
2 & 2 & 1.50 & 2.000 & 1.999 & 2.750 & 3.000\\
\hline
\end{tabular}
\end{table}

To sum up, in Case 1 the theoretical estimates in $L^2$ without a mixed term agree with the numerical results, and in Case 2 the observed $L^2$ rates are also close to those predicted by omitting the mixed term.  These observations suggest that the mixed term may not be needed
in the $L^2$ error bound. We therefore conjecture a 
sharper $L^2$ estimate under the assumptions of \Cref{thm:LapI-rates-discrete}:
\begin{equation}
\label{eq:ti-L2-sharp-alt}
    \|u_{\del,h}-u_0\|_{L^2(\Om)}
    =O\!\left(\del^{\min\{2,m\}}+h^{\min\{m,k\}+1}\right).
\end{equation}
Along $\delta=h^\beta$, the conjecture predicts the order
$\min\{\beta\min\{2,m\},\min\{m,k\}+1\}$.
For the two experiments in \Cref{tab:tiandu-beta-rates-revised}, these orders are $1.5$ and $3$, close to the observed rates $1.524$
and $3.000$. The general estimate instead predicts orders $1.375$ and $2.75$.
Establishing the conjectured $L^2$ rate in the remaining cases
is an open question.  

\subsubsection{Improved rates for integrable kernels.}\label{subsec:ti-integrable-numerics} We now test the improved estimates in \Cref{rem:improved-rates-integrable-kernel}, which assume the
additional condition $\gamma\in L^1(\mathbb{R})$. 
We use discontinuous Galerkin (DG) 
\(P_1\) spaces and the constant kernel
\[
    \ga_{\del}(s)=\frac{3}{2\del^3}\mathbf{1}_{\{|s|\le\del\}},
    \qquad \int_{-\del}^{\del}s^2\ga_{\del}(s)\,\dd s=1.
\]
We take \(m=2\), \(\del=h^\beta\), \(\beta=0.75,0.9\), and
\(N=32,64,128,256\).
For \(0<\beta<1\), \Cref{rem:improved-rates-integrable-kernel} gives
the \(L^2\) bound
\[
    O\!\left(\del^2+\del^{-1/2}h^2+\del^{-2}h^4\right)
    =O\!\left(h^{2\beta}+h^{2-\beta/2}+h^{4-2\beta}\right).
\]
Its predicted exponent is \(\min\{2\beta,2-\beta/2\}\), compared with
\(\min\{2\beta,1+\beta/2\}\) in the general estimate.
\Cref{tab:integrable-kernel-rates-revised} compares these
predictions with the rates computed from $N=128,256$.
At $\beta=0.75$, the improved prediction $1.5$ agrees with the
observed rate $1.518$.
At $\beta=0.9$, integrability raises the predicted exponent from
$1.45$ to $1.55$, while the observed rate is $1.811$.

\begin{table}[H]
\centering
\caption{Predicted and observed convergence rates for DG $P_1$
    approximations with the constant kernel, $m=2$, and $\delta=h^\beta$.  Predicted $L^2$ rates from both the general and improved bounds are listed.
    Observed rates are computed from $N=128,256$.}
\label{tab:integrable-kernel-rates-revised}
\begin{tabular}{c|cc|ccc}
\hline
\(\beta\) & \multicolumn{2}{c|}{Full norm} &
\multicolumn{3}{c}{\(L^2(\Om)\)}\\
\cline{2-6}
 & integrable & observed & general & integrable & observed\\
\hline
0.75 & 1.125 & 1.137 & 1.375 & 1.500 & 1.518\\
0.90 & 1.100 & 1.206 & 1.450 & 1.550 & 1.811\\
\hline
\end{tabular}
\end{table}

These experiments also provide further evidence for
\eqref{eq:ti-L2-sharp-alt}.
For $m=2$, $k=1$, and $0<\beta<1$, the conjectured exponent is $2\beta$.
The bound for integrable kernels reaches this exponent for $\beta\leq4/5$;
for $4/5<\beta<1$, it improves the general estimate but still gives
a smaller exponent than the conjecture.
Thus, at $\beta=0.75$, the improved estimate already gives the
conjectured order $1.5$, whereas at $\beta=0.9$ a gap remains between
the improved prediction $1.55$ and the conjectured order $1.8$.
The observed rates $1.518$ and $1.811$ are close to the conjectured
orders in both cases.

\subsection{Experiments for the type (II) nonlocal Laplacian}
\label{subsec:halfspace-numerics}
We use \(\Om=(0,1)\), \(\nub=+1\), and the same manufactured family
as in \Cref{subsec:fullspace-numerics}. The exterior values are zero on $\Omega_{2\delta}^{+}$.   For each \(m\), we fix
\(f(x)=-u_0''(x)\) and compare \(u_{\del,h}\) with \(E_0u_0\),
first as \(\del\) decreases and then along paths \(\del=h^\beta\).
For the general estimates, we use the fractional kernel
\[
    w_{\del}(s)=(1-\alpha)\del^{\alpha-1}s^{-1-\alpha}
       \mathbf{1}_{(0,\del]}(s),\qquad \alpha=0.95,
\]
where \(\alpha\) is the fractional order.  The normalization is
\(\int_0^\del s w_{\del}(s)\,\dd s=1\).
The scalar functions vanish
in \(\Om_{2\del}^+\).
For this kernel, we solve \eqref{eq:halfspace-Poisson-weak} using CG
spaces of degree \(k\) on uniform meshes with \(h=1/N\) and zero
endpoint values.  The experiments for integrable kernels below use DG spaces.

For each kernel, we use its own full norm,
\[
    \|v\|_{\cS_{w_{\del}}^{\nub}(\Omhh)}^2
    =\|v\|_{L^2(\Omhh)}^2+
       \|\cG_{w_{\del}}^{\nub}v\|_{L^2(\Omh)}^2.
\]
Since the errors vanish in the exterior collar, the scalar term is
integrated over $\Omega$, 
whereas the gradient term also includes
the collar $\Omega_\delta^{+}$.

\subsubsection{Convergence of the continuous problem}
\Cref{thm:LapII-rates-continuous} gives
\[
    \|u_{\del}-E_0u_0\|_{\cS_{w_{\del}}^{\nub}(\Omhh)}
    =O\!\left(\del^{\min\{2,m-1/2\}}\right),\qquad
    \|u_{\del}-u_0\|_{L^2(\Om)}
    =O\!\left(\del^{\min\{2,m\}}\right).
\]
We take \(k=3\), \(m=1,2,3\), and
\(\del=0.2,0.1,0.05,0.025\).
For each fixed \(\del\), we refine \(N=128,256,512\), adding
\(N=1024\) for \(\del=0.025\).  The final spatial refinement changes
both errors by less than \(0.30\%\).
\Cref{tab:halfspace-continuous-rates} reports the rates between
\(\del=0.05\) and \(0.025\), using the errors on the finest meshes.
The observed rates agree with the predicted dependence on \(m\) in
both norms.

\begin{table}[H]
\centering
\caption{Modeling rates for the fractional kernel with \(\alpha=0.95\).}
\label{tab:halfspace-continuous-rates}
\begin{tabular}{c|cc|cc}
\hline
\(m\) & \multicolumn{2}{c|}{Full norm} &
\multicolumn{2}{c}{\(L^2(\Om)\)}\\
\cline{2-5}
& theory & observed & theory & observed\\
\hline
1 & 0.500 & 0.512 & 1.000 & 1.007\\
2 & 1.500 & 1.595 & 2.000 & 1.980\\
3 & 2.000 & 1.968 & 2.000 & 1.998\\
\hline
\end{tabular}
\end{table}

\subsubsection{Convergence of the discrete problem}
Let \(r=\min\{m,k\}\) and set \(\del=h^\beta\).
\Cref{thm:LapII-rates-discrete} gives
\begin{align}
    \|u_{\del,h}-E_0u_0\|_{\cS_{w_{\del}}^{\nub}(\Omhh)}
    &=O\!\left(h^{\beta\min\{2,m-1/2\}}+h^r\right),\notag\\
    \|u_{\del,h}-u_0\|_{L^2(\Om)}
    &=O\!\left(h^{\beta\min\{2,m\}}+h^{r+\beta/2}+h^{r+1}\right).
    \label{eq:tii-L2-bound}
\end{align}
These bounds have the same exponents as those in \Cref{subsec:ti-discrete-numerics}, so we use the same
classification according to whether the mixed term affects the
predicted $L^2$ rate.

\medskip
\noindent\textbf{Case 1: Estimates without a mixed term.} The conditions for absorbing the mixed term are given in
\Cref{subsec:ti-discrete-numerics}.
We first take \(m=3\), \(\del=h^2\), and \(N=16,32,64,128\),
using CG spaces of degrees \(k=1,2,3\).
\Cref{tab:halfspace-degree-rates} reports the rates computed
from $N=64$ and $N=128$.
The observed orders are close to $k$ in the full norm and $k+1$
in $L^2(\Omega)$, as in \Cref{tab:tiandu-degree-rates-revised}.

\begin{table}[H]
\centering
\caption{CG rates for the fractional kernel with \(\alpha=0.95\),
\(m=3\), and \(\del=h^2\).}
\label{tab:halfspace-degree-rates}
\begin{tabular}{c|cc|cc}
\hline
\(k\) & \multicolumn{2}{c|}{Full norm} &
\multicolumn{2}{c}{\(L^2(\Om)\)}\\
\cline{2-5}
& theory & observed & theory & observed\\
\hline
1 & 1.000 & 0.998 & 2.000 & 1.997\\
2 & 2.000 & 1.998 & 3.000 & 2.998\\
3 & 3.000 & 2.995 & 4.000 & 3.994\\
\hline
\end{tabular}
\end{table}

We next consider \((m,k,\beta)=(1,2,0.75)\) and \((3,3,1.25)\).
For both choices, the mixed term in \eqref{eq:LapII-L2-rate-discrete}
can be absorbed.
The observed \(L^2\) rates in \Cref{tab:halfspace-complete-rates}
agree with the predicted orders \(0.75\) and \(2.5\).  The full norm
rate for \(m=1\) is \(0.367\), close to the predicted order \(0.375\).
For \(m=k=3\), the full norm rate decreases from \(2.574\) at
\(N=256\) to \(2.541\) at \(N=512\), approaching the predicted
order \(2.5\).

\begin{table}[H]
\centering
\caption{Additional CG rates for the fractional kernel.
Rates use \(N=256,512\) for \(m=1\) and \(N=128,256\) for \(m=3\).}
\label{tab:halfspace-complete-rates}
\begin{tabular}{ccc|cc|cc}
\hline
\(m\) & \(k\) & \(\beta\) & \multicolumn{2}{c|}{Full norm} &
\multicolumn{2}{c}{\(L^2(\Om)\)}\\
\cline{4-7}
 & & & theory & observed & theory & observed\\
\hline
1 & 2 & 0.75 & 0.375 & 0.367 & 0.750 & 0.750\\
3 & 3 & 1.25 & 2.500 & 2.574 & 2.500 & 2.500\\
\hline
\end{tabular}
\end{table}

\medskip
\noindent\textbf{Case 2: Estimates with a mixed term.} Here $m\geq2$, $r\in\{1,2\}$, and $2r/3<\beta<2$.
We take \(m=2\) with $(k,\beta)=(1,0.75)$ and $(2,1.5)$, using the
same mesh sequences as in the corresponding type (I) experiments
in Section~\ref{subsec:ti-discrete-numerics}.
Table~\ref{tab:halfspace-beta-rates} shows the same pattern as
Table~\ref{tab:tiandu-beta-rates-revised}: the full norm errors
converge at approximately the predicted orders, while the
$L^2(\Omega)$ errors converge faster than predicted.
For \(k=1\), the last three consecutive \(L^2\) rates are
\(1.533,1.490,1.479\), compared with the predicted value \(1.375\).
For \(k=2\), the observed \(L^2\) rate is $3.000$, compared with the
predicted value \(2.75\).

\begin{table}[H]
\centering
\caption{CG rates for the fractional kernel with \(m=2\).
Rates use \(N=512,1024\) for \(k=1\) and \(N=128,256\) for \(k=2\).}
\label{tab:halfspace-beta-rates}
\begin{tabular}{ccc|cc|cc}
\hline
\(m\) & \(k\) & \(\beta\) & \multicolumn{2}{c|}{Full norm} &
\multicolumn{2}{c}{\(L^2(\Om)\)}\\
\cline{4-7}
 & & & theory & observed & theory & observed\\
\hline
2 & 1 & 0.75 & 1.000 & 1.013 & 1.375 & 1.479\\
2 & 2 & 1.50 & 2.000 & 1.998 & 2.750 & 3.000\\
\hline
\end{tabular}
\end{table}

These results also support the estimate without a mixed term
proposed in \eqref{eq:ti-L2-sharp-alt}.
For the type (II) model, we conjecture that, under the assumptions of
Theorem~\ref{thm:LapII-rates-discrete},
\begin{equation}
\label{eq:tii-L2-sharp}
  \|u_{\delta,h}-u_0\|_{L^2(\Omega)}
  =O\!\left(\delta^{\min\{2,m\}}+h^{\min\{m,k\}+1}\right).
\end{equation}
For the two experiments in Table~\ref{tab:halfspace-beta-rates},
the conjectured orders are $1.5$ and $3$, close to the observed
rates $1.479$ and $3.000$.
The general estimate \eqref{eq:tii-L2-bound} instead predicts
orders $1.375$ and $2.75$.

\subsubsection{Improved rates for integrable kernels.}
We now test Remark~\ref{rem:LapII-improved-rates-integrable-kernel}
using DG
\(P_1\) spaces and the constant kernel
\[
    w_{\del}(s)=\frac{2}{\del^2}\mathbf{1}_{(0,\del]}(s),\qquad
    \int_0^\del s w_{\del}(s)\,\dd s=1.
\]
We take $m=2$, $\delta=h^\beta$, $\beta=0.75,0.9$, and
$N=32,64,128,256$.
For $0<\beta<1$, the remark gives the same improved $L^2$ bound as
in Section~\ref{subsec:ti-integrable-numerics}, with predicted
exponent $\min\{2\beta,2-\beta/2\}$.
\Cref{tab:halfspace-integrable-rates} compares this prediction
with the general exponent $\min\{2\beta,1+\beta/2\}$ and the
observed rates computed from $N=128,256$.
The full norm predictions also use the improved bound.

\begin{table}[H]
\centering
\caption{Predicted and observed convergence rates for DG $P_1$
    approximations with the constant kernel, $m=2$, and $\delta=h^\beta$.
    Predicted $L^2$ rates from both the general and improved bounds are listed.
    Observed rates are computed from $N=128,256$.}
\label{tab:halfspace-integrable-rates}
\begin{tabular}{c|cc|ccc}
\hline
\(\beta\) & \multicolumn{2}{c|}{Full norm} &
\multicolumn{3}{c}{\(L^2(\Om)\)}\\
\cline{2-6}
 & integrable & observed & general & integrable & observed\\
\hline
0.75 & 1.125 & 1.135 & 1.375 & 1.500 & 1.522\\
0.90 & 1.100 & 1.245 & 1.450 & 1.550 & 1.810\\
\hline
\end{tabular}
\end{table}

The comparison with the improved bound and the conjecture follows
the same pattern as in Section~\ref{subsec:ti-integrable-numerics}.
At $\beta=0.75$, the observed $L^2$ rate $1.522$ is close to the
improved prediction $1.5$.
At $\beta=0.9$, the observed rate $1.810$ is close to the conjectured
order $1.8$ in \eqref{eq:tii-L2-sharp}, while the improved bound
predicts only $1.55$.

\subsection{Experiments for the heterogeneous localization model}
\label{subsec:heterogeneous-localization-numerics}

We consider \(\Om=(0,1)\) with
\[
    \eta_{\del}(x)=\del\eta(x),\qquad \eta(x)=x(1-x),\qquad
    \ga(z)=\frac32\mathbf{1}_{\{|z|\le1\}}.
\]
The kernel is even and satisfies
\(\int_{-1}^{1}z^2\ga(z)\,\dd z=1\).  The corresponding kernel and
Dirichlet bilinear form are
\begin{align*}
    \ga_{\del,\eta}(x,y)
    &=\frac{1}{\eta_{\del}(x)^3}
      \ga\!\left(\frac{y-x}{\eta_{\del}(x)}\right),\\
    a_{\del}(u,v)=a_{\del,\eta}(u,v)
    &=\int_{\Om}\int_{\Om}\ga_{\del,\eta}(x,y)
      \bigl(u(y)-u(x)\bigr)\bigl(v(y)-v(x)\bigr)\,\dd y\,\dd x.
\end{align*}
We impose the boundary condition \(u(0)=u(1)=0\).

To test the dependence on regularity, set
\[
    A(x)=x(1-x)\left[1+0.4\sin(2\pi x)+0.2\cos(3\pi x)\right],
    \qquad x_c=3/8,
\]
and use, for \(m=1,2,3\),
\[
    u_{0,m}(x)=
    \begin{cases}
        0,&x\le x_c,\\[2pt]
        \displaystyle\frac{A(x)(x-x_c)^{m+1/2}}{\log(e/(x-x_c))},
          &x>x_c.
    \end{cases}
\]
For each \(m\), we fix \(f(x)=-u_{0,m}''(x)\) and compare the nonlocal
solutions with \(u_0=u_{0,m}\).

This family has exactly \(H^{m+1}(\Om)\) regularity.
The derivatives through order \(m\) are continuous at \(x_c\), so no
point mass occurs in the weak derivative of order \(m+1\).
Since \(A(x_c)>0\), for some \(c_m>0\) and sufficiently small \(a>0\),
\[
    u_{0,m}^{(m+1)}(x_c+t)
    \sim\frac{c_m}{\sqrt t\,\log(e/t)}\quad(t\downarrow0),
    \qquad
    \int_0^a\frac{\dd t}{t\log^2(e/t)}
    =\frac{1}{\log(e/a)}<\infty.
\]
Thus \(u_{0,m}\in H^{m+1}(\Om)\cap H^1_0(\Om)\), and
\(f\in H^{m-1}(\Om)\).  For \(0<\varepsilon<1\), the fractional
seminorm of \(u_{0,m}^{(m+1)}\), restricted to pairs on opposite sides
of \(x_c\), is bounded below by a positive multiple of
\[
    \int_0^a\frac{\dd t}{t^{1+2\varepsilon}\log^2(e/t)}=\infty.
\]
Thus \(u_{0,m}\notin H^{m+1+\varepsilon}(\Om)\) for every
\(\varepsilon>0\).

We report the \(L^2(\Om)\) norm and the full norm
\[
    \|v\|_{\cS_{\del,\eta}(\Om)}
    =\left(\|v\|_{L^2(\Om)}^2+a_{\del,\eta}(v,v)\right)^{1/2}
\]
on the space \(\cS_\del(\Om)\) defined in \eqref{eq:hetero-energy-space}.
The energy comparison in \Cref{subsec:hetero-setting} gives
\[
    c\|v\|_{\cS_\del(\Om)}
    \le \|v\|_{\cS_{\del,\eta}(\Om)}
    \le C\|v\|_{\cS_\del(\Om)},
    \qquad v\in\cS_\del(\Om),
\]
with constants \(c,C>0\) independent of \(\del\) and \(h\) for the
fixed \(\eta\) and \(\ga\).  Thus the convergence estimates in
\Cref{sec:NonlocalLaplaceHeteroLoc} also hold for this norm.

The discrete experiments use conforming
continuous Galerkin Lagrange spaces of degree \(k\) on uniform meshes,
with \(h=1/N\) and the two endpoint degrees of freedom eliminated.
Rates are computed from consecutive refinements.

\subsubsection{Convergence of the continuous problem}
The energy estimate in \Cref{thm:hetero-rates-continuous}, together
with the improved \(L^2\) argument based on
\eqref{eq:hetero-H2-H2-consistency}, gives
\[
    \|u_{\del}-u_0\|_{\cS_{\del,\eta}(\Om)}
    =O\!\left(\del^{\min\{2,m\}}\right),\qquad
    \|u_{\del}-u_0\|_{L^2(\Om)}=O(\del^2).
\]
We take \(k=3\), \(m=1,2,3\), and
\(\del=0.2,0.1,0.05,0.025\).
For each fixed \(\del\), we retain the analytic singular part of the
nonlocal solution and approximate the remaining function in the cubic
CG space on a uniform mesh.  We refine \(N=128,256,512\), adding
\(N=1024\) for \(\del=0.025\).  The final spatial refinement changes
both errors by less than \(0.01\%\).
\Cref{tab:heterogeneous-localization-modeling-rates} reports the
rates between \(\del=0.05\) and \(0.025\), using the errors on the
finest meshes.

\begin{table}[H]
\centering
\caption{Modeling rates for the heterogeneous localization model.}
\label{tab:heterogeneous-localization-modeling-rates}
\begin{tabular}{c|cc|cc}
\hline
\(m\) & \multicolumn{2}{c|}{Full norm} &
\multicolumn{2}{c}{\(L^2(\Om)\)}\\
\cline{2-5}
& theory & observed & theory & observed\\
\hline
1 & 1.000 & 1.076 & 2.000 & 2.013\\
2 & 2.000 & 2.000 & 2.000 & 2.000\\
3 & 2.000 & 2.000 & 2.000 & 2.000\\
\hline
\end{tabular}
\end{table}

For \(m=1\), the consecutive full norm rates are
\(1.120,1.090,1.076\), approaching the predicted order one, while
the \(L^2\) rates \(2.021,2.016,2.013\) are close to two.
For \(m=2,3\), both norms converge at approximately second order,
as predicted.

\subsubsection{Convergence of the discrete problem}
Let \(r=\min\{m,k\}\).  The energy estimate and the improved
\(L^2\) estimate in the proof of \Cref{thm:hetero-rates-discrete} are
\begin{align*}
    \|u_{\del,h}-u_0\|_{\cS_{\del,\eta}(\Om)}
    &=O\!\left(\del^{\min\{2,m\}}+h^r\right),\\
    \|u_{\del,h}-u_0\|_{L^2(\Om)}
    &=O\!\left(\del^2+h^{r+1}\right).
\end{align*}
We first test the spatial orders with \(m=3\) and \(\del=h^2\),
using \(N=8,16,32,64,128\) and \(k=1,2,3\).
\Cref{tab:heterogeneous-localization-degree-rates} shows rates close
to \(k\) in the full norm and \(k+1\) in \(L^2(\Om)\).

\begin{table}[H]
\centering
\caption{CG rates for \(m=3\) and \(\del=h^2\).
Rates use \(N=64,128\).}
\label{tab:heterogeneous-localization-degree-rates}
\begin{tabular}{c|cc|cc}
\hline
\(k\) & \multicolumn{2}{c|}{Full norm} &
\multicolumn{2}{c}{\(L^2(\Om)\)}\\
\cline{2-5}
& theory & observed & theory & observed\\
\hline
1 & 1.000 & 0.996 & 2.000 & 1.995\\
2 & 2.000 & 1.997 & 3.000 & 2.997\\
3 & 3.000 & 2.998 & 4.000 & 3.998\\
\hline
\end{tabular}
\end{table}

We next test the improved \(L^2\) estimate at \(m=1\).
Here \(r=1\) for every \(k\ge1\), so the predicted bound is
\(O(\del^2+h^2)\).  We take \(k=1,2,3\), \(\del=h^\beta\), and
\(\beta=0.5,1,2\), giving the respective predicted orders
\(1,2,2\).  We use \(N=32,64,128,256,512\), extending the
\((k,\beta)=(2,2)\) case to \(N=1024,2048,4096,8192\).
\Cref{tab:heterogeneous-localization-m1-L2-rates} reports the final
\(L^2\) rates.

\begin{table}[H]
\centering
\caption{\(L^2(\Om)\) rates for \(m=1\) and \(\del=h^\beta\).
Rates use \(N=256,512\), except for \((k,\beta)=(2,2)\), which uses
\(N=4096,8192\).}
\label{tab:heterogeneous-localization-m1-L2-rates}
\begin{tabular}{c|cc|cc|cc}
\hline
\(k\) & \multicolumn{2}{c|}{\(\beta=0.5\)} &
\multicolumn{2}{c|}{\(\beta=1\)} & \multicolumn{2}{c}{\(\beta=2\)}\\
\cline{2-7}
& theory & observed & theory & observed & theory & observed\\
\hline
1 & 1.000 & 1.014 & 2.000 & 2.000 & 2.000 & 1.999\\
2 & 1.000 & 0.997 & 2.000 & 2.019 & 2.000 & 2.084\\
3 & 1.000 & 0.998 & 2.000 & 2.000 & 2.000 & 2.087\\
\hline
\end{tabular}
\end{table}

The \(\beta=0.5\) and \(\beta=1\) paths give \(L^2\) rates close
to one and two, respectively, as predicted.  For \(k=2\) and
\(\beta=2\), the last three consecutive \(L^2\) rates are
\(2.115,2.092,2.084\), approaching the predicted order two.

Finally, we vary the relative sizes of the modeling and spatial errors.
With \(\del=Ch^\beta\), their predicted orders are
\[
    \min\{\beta\min\{2,m\},r\}\quad\text{in the full norm},\qquad
    \min\{2\beta,r+1\}\quad\text{in }L^2(\Om).
\]
We take \(m=2\), \(C=3\), and the eight choices of \(k\) and
\(\beta\) in \Cref{tab:heterogeneous-localization-beta-rates}.
The rates are computed from \(N=64\) and \(N=128\).

\begin{table}[H]
\centering
\caption{CG rates for \(m=2\) and \(\del=3h^\beta\).
Rates use \(N=64,128\).}
\label{tab:heterogeneous-localization-beta-rates}
\begin{tabular}{ccc|cc|cc}
\hline
\(m\) & \(k\) & \(\beta\) & \multicolumn{2}{c|}{Full norm} &
\multicolumn{2}{c}{\(L^2(\Om)\)}\\
\cline{4-7}
& & & theory & observed & theory & observed\\
\hline
2 & 1 & 0.50 & 1.000 & 1.145 & 1.000 & 1.017\\
2 & 1 & 0.67 & 1.000 & 1.076 & 1.340 & 1.421\\
2 & 1 & 1.00 & 1.000 & 0.997 & 2.000 & 1.996\\
2 & 1 & 2.00 & 1.000 & 0.996 & 2.000 & 1.995\\
\hline
2 & 2 & 1.00 & 2.000 & 1.996 & 2.000 & 2.006\\
2 & 2 & 1.33 & 2.000 & 1.997 & 2.660 & 2.908\\
2 & 2 & 1.50 & 2.000 & 1.997 & 3.000 & 2.996\\
2 & 2 & 2.00 & 2.000 & 1.998 & 3.000 & 2.997\\
\hline
\end{tabular}
\end{table}

For \((m,k,\beta)=(2,1,0.67)\) and \((2,2,1.33)\), the observed
\(L^2\) rates are \(1.421\) and \(2.908\), compared with the
predicted orders \(1.34\) and \(2.66\).  We therefore repeat
these two cases with \(C=5,10\), using the same mesh pair.

\begin{table}[H]
\centering
\caption{Effect of \(C\) on the CG rates for \(m=2\) and
\(\del=Ch^\beta\).  Rates use \(N=64,128\).}
\label{tab:heterogeneous-localization-C-sensitivity}
\begin{tabular}{ccc|cc|cc}
\hline
\(m\) & \(k\) & \(\beta,C\) & \multicolumn{2}{c|}{Full norm} &
\multicolumn{2}{c}{\(L^2(\Om)\)}\\
\cline{4-7}
& & & theory & observed & theory & observed\\
\hline
2 & 1 & 0.67,\ 3 & 1.000 & 1.076 & 1.340 & 1.421\\
2 & 1 & 0.67,\ 5 & 1.000 & 1.115 & 1.340 & 1.360\\
2 & 1 & 0.67,\ 10 & 1.000 & 1.283 & 1.340 & 1.374\\
\hline
2 & 2 & 1.33,\ 3 & 2.000 & 1.997 & 2.660 & 2.908\\
2 & 2 & 1.33,\ 5 & 2.000 & 1.995 & 2.660 & 2.748\\
2 & 2 & 1.33,\ 10 & 2.000 & 1.993 & 2.660 & 2.666\\
\hline
\end{tabular}
\end{table}

For \((m,k,\beta)=(2,1,0.67)\), increasing \(C\) from \(3\) to
\(5\) changes the \(L^2\) rate from \(1.421\) to \(1.360\), close
to the modeling order \(1.34\).
For \((m,k,\beta)=(2,2,1.33)\), the \(C=10\) path gives a rate of
\(2.666\), close to the modeling order \(2.66\), while the
\(C=3\) and \(C=5\) paths give \(2.908\) and \(2.748\).
The changes of \(C\) alter the relative size of the modeling and
spatial errors.

\appendix

\section{Foundations for nonlocal operators on half spaces}
\label{sec:halfspace-analytic-foundations}

\subsection{Pointwise and distributional operators}
Throughout this appendix, \(w_\del\) satisfies
\Cref{assu:halfspace-kernel}.  We define nonlocal half-space operators on general domains and list a few properties without proof. The proofs for half-space operators on \(\R^d\) can be found in \cite{han2023nonlocal}. For general domains in \(\R^d\) it is straightforward to prove these results using similar arguments.

\begin{lemma}\label{lem:shifted-H1-bound}
Let \(D,\cO\subset\R^d\) be open and suppose that
\(D+B_\del(\bm0)\subset\cO\).  Then, for every
\(u\in H^1(\cO)\) and \(\zb\in B_\del(\bm0)\),
\begin{equation*}
    \int_{D}|u(\xb+\zb)-u(\xb)|^2\,d\xb
    \le |\zb|^2\int_{\cO}|\nabla u(\xb)|^2\,d\xb.
\end{equation*}
\end{lemma}

For an open set \(\cO\subset\R^d\), write
\[
    \cOin:=\{\xb\in\cO:\dist(\xb,\R^d\setminus\cO)>\del\},
\]
with \((\R^d)^\circ_\del=\R^d\).
For \(\varphi\in C^\infty_c(\cO)\), the pointwise nonlocal gradient \(\cG_{w_\del,\cO}^{\nub}\varphi:\cOin\to\R^d\) is defined by
\begin{equation*}
\cG_{w_\del,\cO}^{\nub}\varphi(\xb):=\int_{\cO}\chi_{\nub}(\yb-\xb)\frac{\yb-\xb}{|\yb-\xb|}w_\del(\yb-\xb)\bigl(\varphi(\yb)-\varphi(\xb)\bigr)\,d\yb,\qquad \xb\in\cOin.
\end{equation*}

For \(\bsvarphi\in C_c^\infty(\cO;\R^d)\), the pointwise
nonlocal divergence  \(\cD_{w_\del,\cO}^{-\nub}\bsvarphi:\cOin\to\R\) is defined by 
\begin{equation*}
\cD_{w_\del,\cO}^{-\nub}\bsvarphi(\xb):=\int_{\cO}\chi_{-\nub}(\yb-\xb)\frac{\yb-\xb}{|\yb-\xb|}\cdot\bigl(\bsvarphi(\yb)-\bsvarphi(\xb)\bigr)w_\del(\yb-\xb)\,d\yb,\qquad \xb\in\cOin.
\end{equation*}

\begin{definition}[Distributional nonlocal gradient]
\label{def:distributional-gradient}
Let \(\cO\subset\R^d\) be open.  For \(u\in L^2(\cO)\), define
\(\fG_{w_\del,\cO}^{\nub}u\in\mathcal D'(\cOin;\R^d)\) by
\begin{equation*}
\left\langle\fG_{w_\del,\cO}^{\nub}u,\bsvarphi\right\rangle:=-\int_{\cO}u(\xb)\cD_{w_\del,\R^d}^{-\nub}(E_0\bsvarphi)(\xb)\,d\xb,\qquad \bsvarphi\in C_c^\infty(\cOin;\R^d).
\end{equation*}
Here \(E_0\bsvarphi\in C^\infty_c(\R^d;\R^d)\) is the zero extension of \(\bsvarphi\).
\end{definition}

The nonlocal \(H^1\) space is
\begin{equation*}
\begin{aligned}
\cS_{w_\del}^{\nub}(\cO)
 :=\{u\in L^2(\cO):
       \fG_{w_\del,\cO}^{\nub}u
       \in L^2(\cOin;\R^d)\},
\end{aligned}
\end{equation*}
with
\begin{equation*}
\begin{aligned}
\norm{u}_{\cS_{w_\del}^{\nub}(\cO)}^2
 &:=\norm{u}_{L^2(\cO)}^2
    +|u|_{\cS_{w_\del}^{\nub}(\cO)}^2,\\
|u|_{\cS_{w_\del}^{\nub}(\cO)}
 &:=\norm{\fG_{w_\del,\cO}^{\nub}u}_{L^2(\cOin;\R^d)}.
\end{aligned}
\end{equation*}

\begin{proposition}
\label{prop:halfspace-gradient-closed}
The operator in \Cref{def:distributional-gradient} is closed from
\(L^2(\cO)\) to \(L^2(\cOin;\R^d)\).  Consequently,
\(\cS_{w_\del}^{\nub}(\cO)\) is a Hilbert space.
\end{proposition}

\begin{lemma}
\label{lem:halfspace-H1-bound-domain}
Every
\(q\in H^1(\cO)\) belongs to
\(\cS_{w_\del}^{\nub}(\cO)\), its distributional gradient agrees
with the pointwise nonlocal gradient in \(L^2(\cOin;\R^d)\), and
\begin{equation*}
    \norm{\fG_{w_\del,\cO}^{\nub}q}_{L^2(\cOin;\R^d)}
    \le d\norm{\nabla q}_{L^2(\cO;\R^d)}.
\end{equation*}
\end{lemma}

\begin{lemma}
\label{lem:halfspace-finite-range-vanishing}
Let \(D\) and \(\cO\) be two open sets such that \(D \subset \cOin\) and let \(u\in\cS_{w_\del}^{\nub}(\cO)\).  If
\(u=0\) a.e. in \(\cO\setminus D\), then
\begin{equation*}
    \fG_{w_\del,\cO}^{\nub}u=0
    \quad\text{a.e. in }
    \cOin\setminus\overline{D+B_\del(\bm0)}.
\end{equation*}
\end{lemma}

\begin{lemma}
\label{lem:halfspace-locality}
Let \(D\) and \(\cO\) be two open sets such that \(D \subset \cOin\). Suppose that \(u\in L^2(\cO)\)  and
\(\widetilde u\in\cS_{w_\del}^{\nub}(\R^d)\)  satisfy
\(u=\widetilde u\) a.e. in \(D+B_\del(\bm 0)\).  Then
\begin{equation*}
    \left.\fG_{w_\del,\cO}^{\nub}u\right|_D
    = \left.\fG_{w_\del,\R^d}^{\nub}\widetilde u\right|_D.
\end{equation*}
\end{lemma}

\begin{definition}[Distributional nonlocal divergence]
\label{def:distributional-divergence}
Let \(\cO\subset\R^d\) be open.  For \(\vb\in L^2(\cO;\R^d)\), define
\(\fD_{w_\del,\cO}^{-\nub}\vb\in\mathcal D'(\cOin)\) by
\begin{equation*}
\left\langle\fD_{w_\del,\cO}^{-\nub}\vb,\varphi\right\rangle:=-\int_{\cO}\vb(\xb)\cdot\cG_{w_\del,\R^d}^{\nub}(E_0\varphi)(\xb)\,d\xb,\qquad \varphi\in C_c^\infty(\cOin).
\end{equation*}
Here \(E_0\varphi\in C^\infty_c(\R^d)\) is the zero extension of \(\varphi\).
Moreover, \Cref{lem:halfspace-H1-bound-domain} with \(\cO=\R^d\)
gives
\begin{equation*}
\left|\left\langle\fD_{w_\del,\cO}^{-\nub}\vb,\varphi\right\rangle\right|
\le d\norm{\vb}_{L^2(\cO;\R^d)}
       \norm{\nabla\varphi}_{L^2(\cOin;\R^d)}.
\end{equation*}
Thus the distributional divergence defines a continuous linear
functional on \(C_c^\infty(\cOin)\).
\end{definition}

For \(\vb\in H^1(\cO;\R^d)\), the same argument in
\Cref{lem:halfspace-H1-bound-domain} also shows that the pointwise
divergence extends to \(H^1(\cO;\R^d)\) and agrees with
\(\fD_{w_\del,\cO}^{-\nub}\vb\) in \(L^2(\cOin)\), with
\[
    \norm{\fD_{w_\del,\cO}^{-\nub}\vb}_{L^2(\cOin)}
    \le d\norm{\nabla\vb}_{L^2(\cO;\R^{d\times d})}.
\]

We now establish a product rule for the nonlocal operators, which shows that \(\cS_{w_\del}^{\nub}(\R^d)\) is closed under multiplication by smooth functions.  This result is a slight generalization of \cite[Proposition~3.1]{han2023nonlocal} and \cite[Lemma~3.1]{han2023nonlocal}.
\begin{lemma}
\label{lem:halfspace-product-rule}
Let
\(\eta\in C^\infty(\R^d)\cap W^{1,\infty}(\R^d)\) and
\(\bsvarphi\in C_c^\infty(\R^d;\R^d)\).  Then
\begin{equation*}
    \cD_{w_\del,\R^d}^{-\nub}(\eta\bsvarphi)(\xb)
    =\eta(\xb)\cD_{w_\del,\R^d}^{-\nub}\bsvarphi(\xb)
     +S(\eta,\bsvarphi)(\xb),
\end{equation*}
where
\begin{equation*}
S(\eta,\bsvarphi)(\xb):=\int_{\R^d}\chi_{\nub}(\xb-\yb)\frac{\yb-\xb}{|\yb-\xb|}w_\del(\yb-\xb)(\eta(\yb)-\eta(\xb))\cdot\bsvarphi(\yb)\,d\yb.
\end{equation*}
If \(v\in\cS_{w_\del}^{\nub}(\R^d)\), then
\(\eta v\in\cS_{w_\del}^{\nub}(\R^d)\) and
\begin{equation*}
    \fG_{w_\del,\R^d}^{\nub}(\eta v)(\xb)
    =\eta(\xb)\fG_{w_\del,\R^d}^{\nub}v(\xb)
     +\Rb(\eta,v)(\xb)
\end{equation*}
for a.e. \(\xb\in\R^d\), where
\begin{equation*}
\Rb(\eta,v)(\xb):=\int_{\R^d}\chi_{\nub}(\yb-\xb)\frac{\yb-\xb}{|\yb-\xb|}w_\del(\yb-\xb)(\eta(\yb)-\eta(\xb))v(\yb)\,d\yb.
\end{equation*}
Moreover,
\begin{equation*}
    \norm{\Rb(\eta,v)}_{L^2(\R^d;\R^d)}
    \le d\norm{\nabla\eta}_{L^\infty(\R^d;\R^d)}
          \norm{v}_{L^2(\R^d)}.
\end{equation*}
\end{lemma}

\begin{proof}
Expanding the difference in
\(\cD_{w_\del,\R^d}^{-\nub}(\eta\bsvarphi)\) gives
the stated product formula for the divergence. The Lipschitz bound for \(\eta\),
the first moment of \(w_\del\), and Cauchy--Schwarz justify Fubini.
Interchanging the variables and using radiality give
\[
    \int_{\R^d}v(\xb)S(\eta,\bsvarphi)(\xb)\,d\xb
    =\int_{\R^d}\Rb(\eta,v)(\xb)\cdot\bsvarphi(\xb)\,d\xb.
\]
The definition of the distributional gradient and this product formula give
\[
\begin{aligned}
 \left\langle\fG_{w_\del,\R^d}^{\nub}(\eta v),\bsvarphi\right\rangle
 &=-\int_{\R^d}\eta(\xb)v(\xb)
       \cD_{w_\del,\R^d}^{-\nub}\bsvarphi(\xb)\,d\xb\\
 &=-\int_{\R^d}v(\xb)
       \cD_{w_\del,\R^d}^{-\nub}(\eta\bsvarphi)(\xb)\,d\xb
       +\int_{\R^d}v(\xb)S(\eta,\bsvarphi)(\xb)\,d\xb\\
 &=\int_{\R^d}
       \bigl(\eta(\xb)\fG_{w_\del,\R^d}^{\nub}v(\xb)
             +\Rb(\eta,v)(\xb)\bigr)\cdot\bsvarphi(\xb)\,d\xb.
\end{aligned}
\]
Finally, Minkowski's inequality and the normalized first moment of
\(w_\del\) imply
\[
\begin{aligned}
 \norm{\Rb(\eta,v)}_{L^2(\R^d;\R^d)}
 &\le \norm{\nabla\eta}_{L^\infty(\R^d;\R^d)}
       \int_{\R^d}\chi_{\nub}(\zb)|\zb|w_\del(\zb)
       \norm{v(\,\cdot+\zb)}_{L^2(\R^d)}\,d\zb\\
 &=d\norm{\nabla\eta}_{L^\infty(\R^d;\R^d)}
       \norm{v}_{L^2(\R^d)}.
\end{aligned}
\]
Since \(\eta\) is bounded, the preceding distributional identity also
proves \(\eta v\in\cS_{w_\del}^{\nub}(\R^d)\).
\end{proof}

\subsection{Nonlocal operators on bounded domains}
We assume that \(\Om\subset\R^d\) is a bounded \(C^{1,1}\) domain,
let \(\rho_1\) be the radius in the exterior ball condition
\eqref{eq:def_rho_1}, and
take \(0<2\del<\rho_1\). This condition implies
\[
    \dist(\pb+t\nb(\pb),\Om)=t,
    \qquad \pb\in\partial\Om,\quad 0\le t\le2\del.
\]
For \(r\in\{\del,2\del\}\) and
\(\xb\in\widehat\Omega_r\setminus\Om\), choose a nearest boundary
point \(\pb\). Then \(\xb=\pb+t\nb(\pb)\), where
\(t=\dist(\xb,\Om)\). Since the distance function has Lipschitz
constant \(1\), comparison with the exterior point \(\pb+r\nb(\pb)\) gives
\[
    \dist(\xb,\R^d\setminus\widehat\Omega_r)=r-t.
\]
For \(\xb\in\Om\), this distance is at least
\(r+\dist(\xb,\partial\Om)>r\). Consequently,
\begin{equation*}
    \left(\Omhh\right)^\circ_\del=\Omh,
    \qquad
    \left(\Omh\right)^\circ_\del=\Om.
\end{equation*}

These identities give the
operators and spaces used in \Cref{sec:NonlocalLaplaceII}.
Throughout this subsection, use the same abbreviations
\(\fG_{w_\del}^{\nub}:=\fG_{w_\del,\Omhh}^{\nub}\) and
\(\fD_{w_\del}^{-\nub}:=\fD_{w_\del,\Omh}^{-\nub}\).
The graph space \(\cS_{w_\del}^{\nub}(\Omhh)\) and its subspace
\(\cS_{w_\del,0}^{\nub}(\Om)\), equipped with their full norms,
are those defined there.

\begin{lemma}
\label{lem:halfspace-zero-extension-isometry}
Let \(v\in\cS^{\nub}_{w_\del, 0}(\Om)\).  Then
\(E_0v\in\cS_{w_\del}^{\nub}(\R^d)\) and
\begin{equation*}
\fG_{w_\del,\R^d}^{\nub}(E_0v)(\xb)
=
\begin{cases}
\fG_{w_\del,\Omhh}^{\nub}v(\xb),
    &\text{for a.e. }\xb\in\Omh,\\
0,&\text{for a.e. }\xb\in\R^d\setminus\Omh.
\end{cases}
\end{equation*}
Consequently, zero extension \(E_0:\cS^{\nub}_{w_\del, 0}(\Om)\to \cS^{\nub}_{w_\del}(\R^d)\) is an isometric  embedding with 
\[|E_0v|_{\cS_{w_\del}^{\nub}(\R^d)}
 =|v|_{\cS_{w_\del}^{\nub}(\Omhh)},\quad \norm{E_0v}_{\cS_{w_\del}^{\nub}(\R^d)}
  =\norm{v}_{\cS_{w_\del}^{\nub}(\Omhh)}.\] 
\end{lemma}

\begin{proof}
Fix \(\bsvarphi\in C_c^\infty(\R^d;\R^d)\). For
\(0<\varepsilon<\del/2\), choose
\(\eta_\varepsilon\in C_c^\infty(\Omh)\) with
\(0\le\eta_\varepsilon\le1\) and
\(\eta_\varepsilon=1\) on \(\widehat\Omega_{\del-\varepsilon}\).
Since \(v=0\) a.e. in \(\Ompp\),
the definition of the distributional gradient gives
\begin{equation*}
 \int_{\Omh}\fG_{w_\del}^{\nub}v(\xb)\cdot
       \eta_\varepsilon(\xb)\bsvarphi(\xb)\,d\xb
 =-\int_\Om v(\xb)\cD_{w_\del,\R^d}^{-\nub}
       (\eta_\varepsilon\bsvarphi)(\xb)\,d\xb.
\end{equation*}
For \(\xb\in\Om\), the field
\((1-\eta_\varepsilon)\bsvarphi\) vanishes at \(\xb\). On the support of
\(w_\del\), its value at \(\xb+\zb\) can be nonzero only when
\(\del-\varepsilon<|\zb|<\del\). Hence Minkowski's inequality yields
\begin{equation*}
 \norm{\cD_{w_\del,\R^d}^{-\nub}
       ((1-\eta_\varepsilon)\bsvarphi)}_{L^2(\Om)}
 \le \norm{\bsvarphi}_{L^2(\R^d;\R^d)}
       \int_{\{\del-\varepsilon<|\zb|<\del\}}w_\del(\zb)\,d\zb
 \longrightarrow0.
\end{equation*}
Here \(w_\del\) is integrable away from the origin. Letting
\(\varepsilon\to0\) in the testing identity, using dominated
convergence on the left and the preceding bound on the right, gives
\[
 \int_{\Omh}\fG_{w_\del}^{\nub}v(\xb)\cdot\bsvarphi(\xb)\,d\xb
 =-\int_{\R^d}E_0v(\xb)
       \cD_{w_\del,\R^d}^{-\nub}\bsvarphi(\xb)\,d\xb.
\]
This identifies the gradient of \(E_0v\) on \(\R^d\) with the zero
extension of \(\fG_{w_\del}^{\nub}v\), as claimed. The norm identities follow,
since zero extension also preserves the \(L^2\) norm.
\end{proof}

The density theorem \cite[Theorem~3.2]{han2023nonlocal}, with its
assumption that the kernel is nonzero near the origin relaxed under
the present kernel assumption as in
\cite[Section~2, pp.~4--6]{han2024compactnessa}, gives
\begin{equation*}
\overline{C_c^\infty(\Om)}^{\,\cS_{w_\del}^{\nub}(\R^d)}
=\left\{u\in\cS_{w_\del}^{\nub}(\R^d):
u=0\ \text{a.e. in }\R^d\setminus\Om\right\}.
\end{equation*}
Together with \Cref{lem:halfspace-zero-extension-isometry}, this shows
that \(C_c^\infty(\Om)\) is dense in
\(\cS^{\nub}_{w_\del,0}(\Om)\) in the full graph norm.

Let \(u\in\cS_{w_\del}^{\nub}(\Omhh)\) and \(f\in L^2(\Om)\), and set
\(\cL_{w_\del}^{\nub}u:=\fD_{w_\del}^{-\nub}
\fG_{w_\del}^{\nub}u\in\mathcal D'(\Om)\).
By \Cref{def:distributional-divergence,lem:halfspace-locality},
the equation \(-\cL_{w_\del}^{\nub}u=f\) in \(\mathcal D'(\Om)\)
is equivalent to the following identity for
\(v\in C_c^\infty(\Om)\). The density statement above extends it to
every \(v\in\cS_{w_\del,0}^{\nub}(\Om)\):
\begin{equation*}
    \int_{\Omh}\fG_{w_\del}^{\nub}u(\xb)\cdot
       \fG_{w_\del}^{\nub}v(\xb)\,d\xb
    =\int_\Om f(\xb)v(\xb)\,d\xb,
    \qquad v\in\cS_{w_\del,0}^{\nub}(\Om).
\end{equation*}

The distributional divergence and gradient also satisfy an integration
by parts identity. If \(\vb\in L^2(\Omh;\R^d)\) and
\(\fD_{w_\del,\Omh}^{-\nub}\vb\) is represented by a function in
\(L^2(\Om)\), then density gives
\begin{equation*}
    \int_{\Omh}
       \fG_{w_\del,\Omhh}^{\nub}v(\xb)\cdot\vb(\xb)\,d\xb
    +\int_\Om v(\xb)
       \fD_{w_\del,\Omh}^{-\nub}\vb(\xb)\,d\xb
    =0,
    \qquad v\in\cS_{w_\del,0}^{\nub}(\Om).
\end{equation*}

\subsection{Fourier symbol estimates}
\label{subsec:fourier-symbol-estimates}

We collect the Fourier symbol estimates used for the type (I) and type (II) problems.

For a kernel $\ga$ satisfying \Cref{assu:kernel}, write
\begin{equation*}
    \lambda_\ga(\xib)
    :=2\int_{\R^d}\ga(\zb)
    \bigl(1-\cos(2\pi\xib\cdot\zb)\bigr)\,d\zb.
\end{equation*}
The same notation with subscript $\ga_\del$ denotes the symbol in
\eqref{eq:full-space-symbol} where the shorthand notation \(\lambda_\del := \lambda_{\ga_\del}\) is used. In fact, the symbol \(\lambda_\del\) is the Fourier multiplier of the nonlocal operator \(-\cL_\del\), as shown in the following lemma. 

\begin{lemma}\label{lem:fourier-bond-energy}
    For any \(u,v\in \cS_\del(\R^d)\),
    \begin{equation*}
        \int_{\R^d}\int_{\R^d} \ga_\del(\xb - \yb)(u(\xb) - u(\yb))(v(\xb) - v(\yb))d\xb d\yb = \int_{\R^d} \lambda_\del(\xib)\hat{u}(\xib) \overline{\hat{v}(\xib)} d\xib.
    \end{equation*}
\end{lemma}
\begin{proof}
For \(\zb\in \R^d\), denote \(f_{\zb}(\xb):=\ga_\del(\zb)^{1/2}(u(\xb+\zb) - u(\xb))\) and \(g_{\zb}(\xb):=\ga_\del(\zb)^{1/2}(v(\xb+\zb)- v(\xb))\). Then for a.e. \(\zb\in \R^d\), \(f_{\zb}, g_{\zb}\in L^2(\R^d)\) since \(\int_{\R^d}\int_{\R^d}\ga_\del(\zb)(u(\xb+\zb)-u(\xb))^2 d\xb d\zb = \abs{u}_{\cS_\del(\R^d)}^2<\infty\) implies that \(\int_{\R^d}\ga_\del(\zb)(u(\xb+\zb)-u(\xb))^2 d\xb<\infty\) for a.e. \(\zb\in \R^d\). Since \(\widehat{f_{\zb}}(\xib) =  \ga_\del(\zb)^{1/2}\hat{u}(\xib)(e^{2\pi \rmi \xib \cdot \zb} - 1)\), by Plancherel's theorem,
    \begin{align*}
        \int_{\R^d}\int_{\R^d} \ga_\del(\xb - \yb)(u(\xb) - u(\yb))(v(\xb) - v(\yb))d\xb d\yb &= \int_{\R^d} \int_{\R^d} f_{\zb}(\xb) g_{\zb}(\xb) \, d\xb d\zb \\
        &= \int_{\R^d} \ga_\del(\zb) \int_{\R^d} \hat{u}(\xib) (e^{2\pi \rmi \xib \cdot \zb} - 1) \overline{\hat{v}(\xib) (e^{2\pi \rmi \xib \cdot \zb} - 1)} \, d\xib d\zb \\
        &= \int_{\R^d}\lambda_\del(\xib) \hat{u}(\xib) \overline{\hat{v}(\xib)} \, d\xib,
    \end{align*}
    where we used \(|e^{2\pi \rmi \xib \cdot \zb} - 1|^2 = 2(1-\cos(2\pi \xib\cdot \zb))\) in the last step.
\end{proof}

\begin{lemma}\label{lem:bond-symbol-estimates}
Suppose that $\ga$ satisfies \Cref{assu:kernel}. Then there exist constants $c,C>0$, depending only on $d$ and $\ga$, such that
\begin{equation*}
    c|\xib|^2\le\lambda_\ga(\xib)\le C|\xib|^2,
    \quad (|\xib|\le1),
    \quad
    \lambda_\ga(\xib)\ge c,\quad (|\xib|\ge1),
\end{equation*}
and
\begin{equation*}
    \bigl|\lambda_\ga(\xib)-(2\pi)^2|\xib|^2\bigr|
    \le C|\xib|^4,\qquad |\xib|\le1.
\end{equation*}
Let $\ga_\del$ be as in \eqref{eq:gamma-scaling}. Then for every $\del>0$,
\begin{equation}\label{eq:bond-symbol-scaling}
    \lambda_{\ga_\del}(\xib)
    =\del^{-2}\lambda_\ga(\del\xib),\qquad \xib\in\R^d,
\end{equation}
and
\begin{equation*}
    0\le\lambda_{\ga_\del}(\xib)
    \le(2\pi)^2\left(\int_{\R^d}|\zb|^2
    \ga_\del(\zb)\,d\zb\right)|\xib|^2,
    \qquad \xib\in\R^d.
\end{equation*}
If $\ga\in L^1(\R^d)$, then also
\begin{equation*}
    \lambda_{\ga_\del}(\xib)
    \le4\|\ga_\del\|_{L^1(\R^d)},\qquad \xib\in\R^d.
\end{equation*}
\end{lemma}

\begin{proof}
A change of variables gives the scaling identity, while
\(0\le1-\cos t\le\min\{t^2/2,2\}\) gives the two upper bounds.
Radiality and \eqref{eq:gamma:2ndMoment} imply
\[
    \int_{\R^d}\ga(\zb)(\xib\cdot\zb)^2\,d\zb=|\xib|^2.
\]
The fourth moment is finite because \(\ga\) vanishes outside
\(B_1(\bm0)\). Taylor's estimate therefore gives
the approximation bound and the upper and lower bounds for
\(|\xib|\) close to zero.

For the remaining lower bounds, \eqref{eq:gamma-lower-bound} gives
\[
    \lambda_\ga(\xib)
    \ge2c_\ga\int_{B_{c_\ga}(\bm0)}
       \bigl(1-\cos(2\pi\xib\cdot\zb)\bigr)\,d\zb>0,
    \qquad \xib\ne\bm0.
\]
The integral on the right tends to \(|B_{c_\ga}(\bm0)|\) as
\(|\xib|\to\infty\), by the Riemann--Lebesgue lemma.
Continuity and compactness on bounded annuli then prove
the lower bound for \(|\xib|\ge1\) and extend the bounds near zero
to \(|\xib|\le1\).
\end{proof}

Now we turn to the Fourier symbol for the nonlocal gradient. For any \(\phi\in C^\infty_c(\R^d)\), taking the Fourier transform of \(\cG^{\nub}_{w_\del,\R^d} \phi\) yields (\cite{han2023nonlocal})
    \begin{equation*}
        \widehat{\cG^{\nub}_{w_\del,\R^d} \phi}(\xib) = \lab_{w_\del}^{\nub}(\xib) \widehat{\phi}(\xib),\quad \xib\in \R^d,
    \end{equation*}
    where \(\lab_{w_\del}^{\nub}(\xib)\) is the Fourier symbol of \(\cG^{\nub}_{w_\del,\R^d}\) defined by 
\begin{equation*}
    \lab_{w_\del}^{\nub}(\xib) := \int_{\R^d} \chi_{\nub}(\zb)\frac{\zb}{\abs{\zb}} w_\del(\zb) \left(e^{2\pi \rmi \xib\cdot \zb} - 1\right)d\zb.
\end{equation*}
Since \(w_\del(\zb)=\del^{-d-1}w(\zb/\del)\), it follows that
\begin{equation}\label{eq:scaling-Fourier-symbol}
      \lab_{w_\del}^{\nub}(\xib) = \del^{-1}\lab_{w}^{\nub}(\del\xib),
\end{equation}
where 
\begin{equation*}
      \lab_{w}^{\nub}(\xib)
    =
    \int_{B_{1}(\mathbf{0})}
    \chi_{\nub}(\zb)\frac{\zb}{|\zb|}w(\zb)
    \bigl(e^{2\pi i\xib\cdot\zb}-1\bigr)\,d\zb.
\end{equation*}

    By the density theorem \cite[Theorem~3.2]{han2023nonlocal}, with the
    kernel qualification stated above, this Fourier representation
    extends to every \(v\in\cS_{w_\del}^{\nub}(\R^d)\):
    \begin{equation}\label{eq:grad-Fourier}
        \widehat{\fG^{\nub}_{w_\del,\R^d} v}(\xib) = \lab_{w_\del}^{\nub}(\xib) \widehat{v}(\xib),\quad \text{a.e.}\ \xib\in \R^d.
    \end{equation}

The upper bound below holds with \(C=4d\sqrt{2}\pi\) by
\cite[Lemma~2.4]{han2023nonlocal}.
For the lower bound, when \(d=1\) we assume in addition that the radial
profile of \(w\) has a nonnegative, nonincreasing representative on
\((0,\infty)\). Under these assumptions,
\cite[Lemma~2.1]{han2024compactnessa} gives a constant \(c=c(w,d)>0\) such that
\begin{equation}\label{eq:grad-Fourier-symbol-bounds}
      c\min\{1,\abs{\xib}\}
      \le\abs{\lab_w^{\nub}(\xib)}
      \le C\abs{\xib},\quad\xib\in\R^d.
\end{equation}

When \(\abs{\xib}\le 1\), we have a finer estimate as shown in the following lemma.
\begin{lemma}
\label{lem:small-frequency-symbol}
Let \(w\) satisfy \Cref{assu:halfspace-kernel}
and let \(\lab_w^{\nub}\) be its Fourier symbol.
Then there exists a constant \(C>0\) such that
\[
    \left|
        |\lab_{w}^{\nub}(\xib)|^{2}
        -4\pi^{2}|\xib|^{2}
    \right|
    \leq C|\xib|^{4},
    \qquad |\xib|\leq 1.
\]
\end{lemma}

\begin{proof}
The real and imaginary parts of \(\lab_w^{\nub}(\xib)\) are 
\[
\Re \lab_w^{\nub}(\xib) = -\int_{\R^d}  \chi_{\nub}(\zb)\frac{\zb}{|\zb|} w(\zb) (1-\cos(2\pi \xib\cdot\zb))\, d\zb,
\quad
\Im \lab_w^{\nub}(\xib) = \int_{\R^d}  \chi_{\nub}(\zb)\frac{\zb}{|\zb|} w(\zb) \sin(2\pi \xib\cdot\zb)\, d\zb.
\]
    
Using \(\int_{\R^d} \abs{\zb}w(\zb)d\zb = 2d\), together with
\(|1-\cos t|\le t^2/2\) and \(|\sin t-t|\le|t|^3/6\), one obtains
\[
 \abs{\Re\lab_w^{\nub}(\xib)}\le C|\xib|^2,
 \quad
 \abs{\Im\lab_w^{\nub}(\xib)-2\pi\xib}\le C|\xib|^3.
\]
Indeed, the second and third moments of \(w\) are bounded by its first
moment, since \(w\) vanishes outside \(B_1(\bm0)\).
Using the upper bound in \eqref{eq:grad-Fourier-symbol-bounds}, we conclude that
\[ 
 \left|\abs{\lab_w^{\nub}(\xib)}^2-4\pi^2|\xib|^2\right| \le \abs{\Re\lab_w^{\nub}(\xib)}^2
       +\abs{\Im\lab_w^{\nub}(\xib)-2\pi\xib}
       \bigl(\abs{\Im\lab_w^{\nub}(\xib)}
             +2\pi|\xib|\bigr) \le C|\xib|^4.
\]
This finishes the proof.
\end{proof}

Now we are ready to compare the type (I) and type (II) energies through their Fourier symbols. The symbols above correspond to operators of different orders:
$\lambda_{\ga_\del}$ scales as $\del^{-2}$, whereas
$\lab_{w_\del}^{\nub}$ scales as $\del^{-1}$.
The energy multipliers $\lambda_{\ga_\del}$ and
$|\lab_{w_\del}^{\nub}|^2$ both scale as $\del^{-2}$.
By \Cref{lem:bond-symbol-estimates,lem:small-frequency-symbol},
their unscaled forms both
approximate $4\pi^2|\xib|^2$ with an error bounded by $C|\xib|^4$
for $|\xib|\le1$. The following comparison holds at all
frequencies.

\begin{lemma}\label{lem:symbol-comparison}
Let $\hat\ga$ satisfy \Cref{assu:kernel} and let $w$ satisfy
\Cref{assu:halfspace-kernel}. If $\hat\ga\in L^1(\R^d)$, then there
exists a constant $C=C(\hat\ga,w,d)>0$ such that
\begin{equation*}
    \lambda_{\hat\ga}(\xib)
    \le C|\lab_w^{\nub}(\xib)|^2,
    \qquad \xib\in\R^d.
\end{equation*}
Here $\lambda_{\hat\ga}$ denotes the Fourier symbol associated with
$\hat\ga$.
\end{lemma}

For $\hat\ga_\del(\zb)=\del^{-d-2}\hat\ga(\zb/\del)$,
\eqref{eq:bond-symbol-scaling} and
\eqref{eq:scaling-Fourier-symbol} give
\begin{equation*}
    \lambda_{\hat\ga_\del}(\xib)
    \le C|\lab_{w_\del}^{\nub}(\xib)|^2,
    \qquad \xib\in\R^d,
\end{equation*}
with $C$ independent of $\del$. 
\begin{lemma}
\label{lem:halfspace-energy-comparison}
Let \(\hat\ga\) satisfy \Cref{assu:kernel} and
\(\hat\ga\in L^1(\R^d)\). Define
\[
    \hat\ga_\del(\zb):=\del^{-d-2}\hat\ga(\zb/\del).
\]
Then there exists \(C=C(\hat{\ga},w,d)>0\) such that
\begin{equation*}
    \abs{v}_{\cS_{\hat\ga_\del}(\Om)}
    \le C
    \abs{v}_{\cV_\del},
    \qquad \forall v\in\cV_\del.
\end{equation*} 
\end{lemma} 

\begin{proof}
By \Cref{lem:halfspace-zero-extension-isometry},
\(E_0v\in\cS_{w_\del}^{\nub}(\R^d)\).
The scaled form of \Cref{lem:symbol-comparison} and the Fourier energy
identities give
\[
\begin{aligned}
 \abs{E_0v}_{\cS_{\hat\ga_\del}(\R^d)}^2
 &=\int_{\R^d}\lambda_{\hat\ga_\del}(\xib)
       \abs{\widehat{E_0v}(\xib)}^2\,d\xib\\
 &\le C\int_{\R^d}\abs{\lab_{w_\del}^{\nub}(\xib)}^2
       \abs{\widehat{E_0v}(\xib)}^2\,d\xib
 =C\abs{v}_{\cV_\del}^2.
\end{aligned}
\]
The constant is independent of \(\del\) by the scaling identities.
Since \(E_0v=0\) outside \(\Om\) and \(\hat\ga_\del\) vanishes outside
\(B_\del(\bm0)\), its bond seminorm on \(\R^d\) equals
\(\abs{v}_{\cS_{\hat\ga_\del}(\Om)}\).
Taking square roots proves the result.
\end{proof}

\begin{proof}[Proof of \Cref{prop:halfspace-boundary-layer}]
    Since \(2\del<\rho_1\), the horizon condition
\(\del<4\rho_1\) in \Cref{prop:nonlocal-boundary-layer} is satisfied.
    Let \(\hat\ga:=\frac{d+2}{\abs{B_1(\bm 0)}}\1_{B_1(\bm 0)}\) and \(\hat\ga_\del(\zb):=\del^{-d-2}\hat\ga(\zb/\del)\). Then \(\hat\ga\) satisfies \Cref{assu:kernel} and  \(\hat\ga\in L^1(\R^d)\). The conclusion then follows from the type (I) boundary layer estimate \Cref{prop:nonlocal-boundary-layer} applied to \(\hat\ga_\del\) and the energy comparison in \Cref{lem:halfspace-energy-comparison}:
    \[
        \norm{v}_{L^2(\Omega^-_{a\del})} \le C_a \del \abs{v}_{\cS_{\hat\ga_\del}(\Om)} \le C_a \del \abs{v}_{\cV_\del}.
    \]
\end{proof}

\section{Boundary-layer Estimates}

We cite the following classical results from \cite[Theorem~1.5.1.2 and Corollary~1.5.1.6]{grisvard2011elliptic} to define traces of Sobolev functions on boundaries of domains. These results also provide a characterization of functions with vanishing traces. 

\begin{theorem}
\label{thm:normal-traces}
Let $m\ge0$ be an integer, and let
$\Omega\subset\mathbb R^d$ be a bounded $C^{m,1}$ domain. Denote the outward unit normal of \(\Om\) by $\nb$. Let \(\Tr\) denote the operator given by \(\Tr u=u|_{\partial\Om}\) when \(u\) is a smooth function on \(\bar\Om\). Then the mapping 
\[
u\mapsto \left(\Tr u, \Tr\frac{\partial u}{\partial\nb}, \cdots, \Tr\frac{\partial^m u}{\partial\nb^m}\right)
\]
which is defined for \(u\in C^{m,1}(\bar\Om)\), has a unique continuous extension as an operator
\[
    H^{m+1}(\Omega) \to \prod_{j=0}^m H^{m+\frac12-j}(\partial\Omega).
\]
This operator has a continuous right inverse. For \(j=0,\cdots,m\), we write \(\Tr\frac{\partial^j u}{\partial\nb^j}=\Tr \partial_{\nb}^j u\) with the convention that \(\partial_{\nb}^0 u=u\). For \(k\in \{m,m+1\}\), \(u\in H^{k}_0(\Om)\) if and only if \(u\in H^k(\Om)\) and 
\[
    \Tr \partial_{\nb}^j u = 0, \quad 0\le j\le k-1,
\]
where
$H_0^k(\Omega)=\overline{C_c^\infty(\Omega)}^{\,H^k(\Omega)}$.
\end{theorem}

The next lemma shows a quantitative estimate of Sobolev norms of functions in a boundary layer in terms of the width of the layer, provided that certain vanishing trace conditions hold. A special case when \(m=0\) and \(\Om\) is a polyhedral domain is given in \cite[Lemma~10]{meng2025asymptotically}. When \(m=0\), there are no trace conditions imposed and the lemma is a simple consequence of \cite[Lemma~5.1]{pakhnin2013operator}. Recall that \(\Om_\del^-\) and \(\Om_\del^+\) are defined in \eqref{eq:notation-collars}.
\begin{lemma}\label{lem:boundary-layer-est}
    Let \(m\ge 0\) be an integer and assume that \(\Om\) is a bounded 
    \(C^{m,1}\) domain. Then there exists a constant
    \(\rho_0\in (0,1)\) depending only on \(\Om\) with the following property: for any
    \(\del_0\in (0,\rho_0)\) there is
    \(C=C(m,d,\Om,\del_0)>0\) 
    such that for any \(\del\in (0,\del_0]\) and any
    \(u\in H^{m+1}(\Om^-_{\del_0})\) satisfying  
    \begin{equation}
\label{eq:zero-normal-traces}
    \Tr(\normal^j u)=0,
    \qquad
    j=0,\ldots,m-1,
\end{equation}
one has
\begin{equation}\label{eq:boundary-layer-est}
    \forall \abs{\be}\le m,\quad \norm{\p^\be u}_{L^2(\Omm)} \le C \del^{m+\frac{1}{2}-|\be|}\norm{u}_{H^{m+1}(\Om^-_{\del_0})}.
\end{equation}
The same assertion holds with
\(\Om^-_\del\) and \(\Om^-_{\del_0}\) replaced by
\(\Om^+_\del\) and \(\Om^+_{\del_0}\), respectively.
\end{lemma}

\begin{proof}
Without loss of generality, we only prove the estimate for \(\Om^-_\del\). Throughout the proof,
$C$ depends only on $m$, $d$, $\Om$ and \(\del_0\), and may change
from line to line.

First we prove the case $m=0$.
By \cite[Lemma~5.1]{pakhnin2013operator}, there exist $\rho_0\in (0,1)$ and $C>0$ depending on \(\Om\) such that
\begin{equation}\label{eq:global-boundary-strip}
    \|v\|_{L^2(\Om^-_{\vep})}^2
    \le C\vep\|v\|_{L^2(\Om)}\|v\|_{H^1(\Om)}\le C\vep \|v\|_{H^1(\Om)}^2,
    \qquad
    v\in H^1(\Om),\quad 0<\vep<\rho_0.
\end{equation}
Although the cited result is stated for $C^1$ domains,
its proof applies to Lipschitz domains: the changes of
variables to rectifying Lipschitz boundary charts are bi-Lipschitz,
and composition with such maps is bounded from $H^1$ to $H^1$.

Fix $0<\del_0\le\rho_0$ and
$u\in H^1(\Om^-_{\del_0})$. Choose a cut-off function
$\eta\in C^\infty([0,\infty))$ satisfying
\[
    0\le\eta\le1,\qquad
    \eta=1\ \text{on }[0,1/2],\qquad
    \eta=0\ \text{on }[3/4,\infty),
\]
and define
\[
    v(\xb)=
    \begin{cases}
        \eta\!\left(
            \dfrac{\dist(\xb,\partial\Om)}{\del_0}
        \right)u(\xb),
        & \xb\in\Om^-_{\del_0},\\[1ex]
        0,
        & \xb\in\Om\setminus\Om^-_{\del_0}.
    \end{cases}
\]
Then $v\in H^1(\Om)$ and
\[
    \|v\|_{L^2(\Om)}
    \le \|u\|_{L^2(\Om^-_{\del_0})},
    \qquad
    \|v\|_{H^1(\Om)}
    \le C\del_0^{-1}
        \|u\|_{H^1(\Om^-_{\del_0})}.
\]
For $0<\del<\del_0/2$, we have $v=u$ on $\Om^-_\del$.
Consequently, \eqref{eq:global-boundary-strip} gives
\[
    \|u\|_{L^2(\Om^-_\del)}^2
    \le C\frac{\del}{\del_0}
        \|u\|_{H^1(\Om^-_{\del_0})}^2.
\]
For $\del_0/2\le\del\le\del_0$, the same estimate follows from
\[
    \|u\|_{L^2(\Om^-_\del)}^2
    \le \|u\|_{L^2(\Om^-_{\del_0})}^2
    \le \frac{2\del}{\del_0}
        \|u\|_{H^1(\Om^-_{\del_0})}^2.
\]
Thus the case \(m=0\) is proved as
\begin{equation*}
    \|u\|_{L^2(\Om^-_\del)}^2
    \le C\frac{\del}{\del_0}
        \|u\|_{H^1(\Om^-_{\del_0})}^2,
    \qquad 0<\del\le\del_0.
\end{equation*}

Now we consider the case $m\ge 1$. If $|\beta|=m$, the required estimate follows directly from the case $m=0$ applied to $\partial^\be  u\in H^1(\Om^-_{\del_0})$:
\[
    \|\partial^\be u\|_{L^2(\Om^-_\del)}^2
    \le C\frac{\del}{\del_0}
        \|\partial^\be u\|_{H^1(\Om^-_{\del_0})}^2
    \le C\frac{\del}{\del_0}
        \|u\|_{H^{m+1}(\Om^-_{\del_0})}^2.
\]
Fix now $|\beta|<m$.
Given \eqref{eq:zero-normal-traces},
\Cref{thm:normal-traces} yields
\begin{equation}
\label{eq:zero-cartesian-traces}
    \Tr(\partial^\al u)=0,
    \qquad |\al|<m.
\end{equation}
Similar to the proof of \cite[Lemma~5.1]{pakhnin2013operator}, we may reduce the estimate on a bounded domain to a hypercube using the standard way of a finite cover of \(C^{m,1}\) open sets and a partition of unity. Therefore, it suffices to prove that \eqref{eq:boundary-layer-est} holds on \(\Ga\times (0,\del)\) for \(\Ga\subset \R^{d-1}\) and \(\del\in (0,\del_0)\), i.e., if \(\partial^j_{x_d} u(\xb',0)=0\) for \(j=0,\cdots,m-1\), then for \(\be=(\beta',\be_d)\) with \(|\be|<m\),
\[
    \norm{\partial^{\be'}_{\xb'}\partial^{\be_d}_{x_d} u}_{L^2(\Ga\times (0,\del))} \le C\del^{m+\frac{1}{2}-|\be|} \|u\|_{H^{m+1}(\Ga\times (0,\del_0))}.
\]
We claim that for \(g\in H^k(0,\del)\) with \(g(0)=g'(0)=\cdots=g^{(k-1)}(0)=0\), 
\[
\norm{g}_{L^2(0,\del)}\le C\del^{k}\norm{g^{(k)}}_{L^2(0,\del)}.
\]
Indeed, by Taylor's theorem with integral remainder,
\[
    g(r)=\int_0^r\frac{(r-t)^{k-1}}{(k-1)!}g^{(k)}(t)\,dt,
\]
and Cauchy--Schwarz, followed by integration in \(r\), gives
\[
    \int_0^\del \abs{g(r)}^2\,dr
    \le C\int_0^\del \left(\int_0^r (r-t)^{2k-2}\,dt\right)
        \left(\int_0^r \abs{g^{(k)}(t)}^2\,dt\right)\,dr
    \le C\del^{2k}\int_0^\del \abs{g^{(k)}(t)}^2\,dt,
\]
which proves the claim. Let \(g(t):=\partial^{\be'}_{\xb'}\partial^{\be_d}_{x_d} u(\xb',t)\). Then \(g^{(j)}(t)=\partial^{\be'}_{\xb'}\partial^{\be_d+j}_{x_d} u(\xb',t)\) and \(g^{(j)}(0)=0\) for \(j=0,\ldots,m-|\be|-1\) by \eqref{eq:zero-cartesian-traces}. Applying the claim with \(k=m-|\be|\) gives
\[
    \int_{0}^\del \left|\partial^{\be'}_{\xb'}\partial^{\be_d}_{x_d} u(\xb',t)\right|^2 dt
    \le C\del^{2(m-|\be|)} \int_0^\del \left|\partial^{\be'}_{\xb'}\partial^{\be_d+m-|\be|}_{x_d} u(\xb',t)\right|^2 dt.
\]
Integrating over \(\xb'\in\Ga\) and applying the proved estimate for \(|\al|=m\) where \(\al=(\be', \be_d+m-|\be|)\) gives
\[
    \|\partial^{\be'}_{\xb'}\partial^{\be_d}_{x_d} u\|_{L^2(\Ga\times (0,\del))}^2 \le C\del^{2(m-|\be|)} \|\partial^\al u\|_{L^2(\Ga\times (0,\del_0))}^2\le C\del^{2(m-|\be|)+1} \|u\|_{H^{m+1}(\Ga\times (0,\del_0))}^2.
\]
This finishes the proof.
\end{proof}

\begin{remark}
    There are a few important special cases of
    \Cref{lem:boundary-layer-est} which we will use in our analysis.
    \begin{enumerate}
        \item When \(m=0\), the zero normal trace condition \eqref{eq:zero-normal-traces} is empty, \(\Om\) is Lipschitz, and \(\abs{\be}=0\). Then \eqref{eq:boundary-layer-est} reads as \begin{equation*}
            \norm{u}_{L^2(\Omm)} \le C \del^{\frac{1}{2}}\norm{u}_{H^{1}(\Om^-_{\del_0})},\quad\forall u\in H^1(\Om^-_{\del_0}).
        \end{equation*}
        \item When \(m=1\), \eqref{eq:zero-normal-traces} reduces to the usual zero trace condition \(u=0\) on \(\p\Om\), and \(\Om\in C^{1,1}\). Taking \(\be=0\) in \eqref{eq:boundary-layer-est} gives \begin{equation*}
            \norm{u}_{L^2(\Omm)} \le C \del^{\frac{3}{2}}\norm{u}_{H^{2}(\Om^-_{\del_0})},\quad\forall u\in H^2(\Om^-_{\del_0}) \text{ with } \Tr u = 0.
        \end{equation*}
        \item For \(m\ge 1\), taking \(|\be|=1\) in \eqref{eq:boundary-layer-est} gives \begin{equation*}
            \norm{\nabla u}_{L^2(\Omm, \R^d)} \le C \del^{m-\frac{1}{2}}\norm{u}_{H^{m+1}(\Om^-_{\del_0})},\quad\forall u\in H^{m+1}(\Om^-_{\del_0}) \text{ satisfying \eqref{eq:zero-normal-traces} }.
        \end{equation*}
        \item For \(m\ge 2\), taking \(|\be|=2\) in \eqref{eq:boundary-layer-est} gives \begin{equation*}
            \norm{\nabla^2 u}_{L^2(\Omm, \R^{d\times d})} \le C \del^{m-\frac{3}{2}}\norm{u}_{H^{m+1}(\Om^-_{\del_0})},\quad\forall u\in H^{m+1}(\Om^-_{\del_0}) \text{ satisfying \eqref{eq:zero-normal-traces} }.
        \end{equation*}
    \end{enumerate}
\end{remark}

Combining the estimates on both sides yields the following corollary.
\begin{corollary}\label{cor:local-boundary-double-layer-est}
    Let \(m\ge1\) be an integer, let \(a,b>0\) be fixed, and let
    \(\Om\subset \mathbb R^d\) have \(C^{m,1}\) boundary. Let \(\del_0\)
    be the fixed collar width from \Cref{lem:boundary-layer-est}, and suppose
    that \(\del>0\) satisfies only \(b\del\le\del_0\). Suppose that \(S\) is
    piecewise defined by \(S^-\in H^{m+1}(\Om)\) and
    \(S^+\in H^{m+1}(\Om^+_{\del_0})\), and that the normal traces of \(S\) vanish up to
order \(m-1\) on both sides, namely
\begin{equation*}
    \Tr(\normal^j S^-)
    =\Tr(\normal^j S^+)=0,
    \qquad
    j=0,\ldots,m-1 .
\end{equation*}
Then there exists a constant 
\(C=C(m,d,\Om,\del_0,a,b)>0\), independent of \(\del\), such that
\[
\|\nabla^k S\|_{L^2((\partial\Om)^{a,b}_\del)}\le
    C\del^{m+\frac{1}{2}-k}
    \left(\norm{S^-}_{H^{m+1}(\Om)}+\norm{S^+}_{H^{m+1}(\Om^+_{\del_0})}\right),\quad m\ge k.
\]
\end{corollary}

\bibliographystyle{abbrv} 

\bibliography{mylibbetterbibtex}

\end{document}